\documentclass[12pt]{article}%
\usepackage{amsmath}
\usepackage{amsfonts}
\usepackage{amssymb}
\usepackage{graphicx}%
\newtheorem{theorem}{Theorem}

\newtheorem{lemma}[theorem]{Lemma}

\newenvironment{proof}[1][Proof]{\noindent\textbf{#1.} }{\ \rule{0.5em}{0.5em}}

\def\lfhook#1{\setbox0=\hbox{#1}{\ooalign{\hidewidth
\lower1.5ex\hbox{'}\hidewidth\crcr\unhbox0}}}
\begin{document}

\title{Asymptotic analysis by the saddle point method of the Anick-Mitra-Sondhi model}
\author{Diego Dominici\thanks{ddomin1@uic.edu} \ and Charles
Knessl\thanks{knessl@uic.edu}\\Department of Mathematics, Statistics and Computer Science \\University of Illinois at Chicago (M/C 249) \\851 South Morgan Street \\Chicago, IL 60607-7045, USA}
\maketitle

\begin{abstract}
We consider a fluid queue where the input process consists of $N$ identical
sources that turn \textbf{on} and \textbf{off }at exponential waiting times.
The server works at the constant rate $c$ and an \textbf{on} source generates
fluid at unit rate. This model was first formulated and analyzed by Anick,
Mitra and Sondhi in \cite{MR84a:68020}. We obtain an alternate representation
of the joint steady state distribution of the buffer content and the number of
\textbf{on} sources. This is given as a contour integral that we then analyze
in the limit $N\rightarrow\infty.$ We give detailed asymptotic results for the
joint distribution, as well as the associated marginal and conditional
distributions. In particular, simple conditional limits laws are obtained.
These shows how the buffer content behaves conditioned on the number of active
sources and vice versa. Numerical comparisons show that our asymptotic results
are very accurate even for $N=20.$

\end{abstract}

\section{Introduction}

In traditional queueing theory one usually considers a system with one or more
servers, at which customers arrive to receive some kind of service. As all
servers may be busy, the customers may have to wait for service in a queue.
Since there is uncertainty about the actual arrival times of the customers
and/or their service requirements, queues are typically modeled
stochastically, as random processes.

In the last twenty years, models have appeared in which a continuous quantity,
referred to as a \textit{fluid}, takes the role of the number of customers or
queue length. In these models, fluid flows into a reservoir according to some
stochastic process. The server may be viewed as a tap at the bottom of the
reservoir, allowing fluid to flow out. The rate at which this happens is often
constant, but may also be stochastic. Since the fluid reservoir takes the role
of the traditional customer queue, it is often referred to as a \textit{fluid
queue}. A third term which is often encountered is \textit{fluid buffer},
stressing the fact that the storage of fluid is temporary, to prevent loss of
fluid at times when the input rate exceeds the output rate. Such models are
used as approximations to discrete queuing models of manufacturing systems,
high-speed data networks, transportation systems, in the theory of dams, etc.
The distribution of the fluid level or buffer content, the average buffer
content, the overflow probability (in case of a finite buffer) and the output
process are the main quantities of interest in these fluid models.

In the last decade, the literature on queueing theory has paid considerable
attention to \textit{Markov-modulated fluid models} \cite{MR1357855},
\cite{fluid}. In these models, a fluid buffer is either filled or depleted, or
both, at rates which are determined by the current state of a
\textit{background Markov process}, also called \textit{Markovian random
environment}.

\subsection{The Markov-modulated fluid model}

In this section we describe a general model of fluid entering and leaving a
single buffer system \cite{MR97k:60252}. Let $X(t)$ denote the amount of fluid
at time $t$ in the buffer. Furthermore, let $Z(t)$ be a continuous-time Markov
process. $Z(t)$ will be said to evolve \textquotedblleft in the
background\textquotedblright. We will assume that $Z(t)$ has a finite state
space $\mathcal{N}$ with $\mathcal{N}=\left\{  0,1,\ldots,N\right\}  $.

The buffer content $X(t)$ is regulated (or driven) by $Z(t)$ in such a way
that the \textit{net input rate} into the buffer (i.e., the rate of change of
its content) is $\mathrm{d}\left[  Z(t)\right]  $. The function $\mathrm{d}%
\left(  \cdot\right)  $ is called the \textit{drift function}. When the buffer
capacity is infinite, the dynamics of $X(t)$ are given by%
\begin{equation}
\frac{dX}{dt}=\left\{
\begin{array}
[c]{c}%
\mathrm{d}\left[  Z(t)\right]  ,\quad X(t)>0\\
\max\left\{  \mathrm{d}\left[  Z(t)\right]  ,0\right\}  ,\quad X(t)=0
\end{array}
\right.  . \label{fluid infinity}%
\end{equation}
The condition at $X(t)=0$ ensures that the process $X(t)$ does not become
negative. When the buffer capacity is finite, say $K,$ the dynamics are given
by%
\[
\frac{dX}{dt}=\left\{
\begin{array}
[c]{c}%
\mathrm{d}\left[  Z(t)\right]  ,\quad0<X(t)<K\\
\max\left\{  \mathrm{d}\left[  Z(t)\right]  ,0\right\}  ,\quad X(t)=0\\
\min\left\{  \mathrm{d}\left[  Z(t)\right]  ,0\right\}  ,\quad X(t)=K
\end{array}
\right.  .
\]
The condition at $X(t)=K$ prevents the buffer content from exceeding $K.$

One of the main reasons why Markov-modulated fluid models have attracted so
much attention is that they are relevant for modelling telecommunications
networks. In the literature, more attention has been paid to models in which
the buffer capacity is infinitely large, because this case is easier to
analyze. For most practical situations in telecommunications the infinite
buffer case is a good approximation for the finite buffer case, since overflow
of the buffer is assumed to be extremely rare. For highly loss-sensitive
traffic, allowed loss fractions in the order of $10^{-6}$ to $10^{-10}$ are
typical \cite{MR98h:60131}.

Unlike classical queueing models which assume renewal arrivals, Markov
modulated fluid models can capture autocorrelations in arrival processes. The
continuous nature of the\ fluid also makes them more tractable analytically.
Many results have been obtained for a variety of fluid queueing systems
\cite{MR97j:60165}, \cite{MR2003c:60147}, \cite{MR2001i:60155},
\cite{MR2000f:60143}, \cite{MR2002e:60150}, \cite{MR1932878}, \cite{MR1926146}%
, \cite{MR1932126}, \cite{MR2002b:60167}, \cite{MR2002f:60182},
\cite{MR2000d:60155}, \cite{MR2002j:60174}.

In many applications, the process $Z(t)$ evolves as a finite state
\textit{birth-death process. }Birth-death processes have wide applications in
many practical problems \cite{MR36:913}. They can be regarded as continuous
time analogs of random walks.

The paper which has become the main reference for birth-death fluid models was
written by Anick, Mitra and Sondhi \cite{MR84a:68020}. We shall refer to their
model as the AMS model. The AMS model describes an infinitely large fluid
buffer which is fed by $N$ identical exponential \textbf{on}-\textbf{off}
sources and emptied by an output channel with constant capacity $c.$ An
\textbf{on} source turns \textbf{off} at rate $1$ and an \textbf{off} source
turns \textbf{on} at rate $\lambda.$ Thus, the net input is regulated by a
specific birth-death process $Z(t)$ with state space $\mathcal{N},$ birth rate
$\lambda_{k}=\lambda\left(  N-k\right)  ,$ death rate $\mu_{k}=k$ and drift
function $\mathrm{d}(k)=k-c.$ The rates are conditioned on $Z(t)=k,\ 0\leq
k\leq N.$

In a way, the AMS model was a generalization of \ the earlier model of Kosten
\cite{MR58:3798a}, \cite{MR58:3798b}, where the limiting case $N\rightarrow
\infty,\ \lambda\rightarrow0$ (with $\lambda N$ fixed) is considered.
According to the authors in \cite{MR84a:68020}, their model and its variants
had also been proposed in other previous papers \cite{Chu}, \cite{hashida},
\cite{rudin}.

In \cite{kosten1}, Kosten generalizes the AMS model by considering the same
problem for a \textit{Multi Group Finite Source} (MGFS) system, consisting of
$m$ groups of identical sources and a single shared buffer. In
\cite{MR88b:90057} Kosten considers a central processor, working at uniform
speed, that receives information at uniform rates from a number of sources and
a buffer which stores the information that cannot be handled directly. He
shows that the stationary probability of overflow $G(x)$ satisfies $G(x)\sim
Ce^{-\alpha x}$ as $x\rightarrow\infty.$ For fairly complex aggregates of
sources, a method is developed (called the decomposition method) to determine
the exponent $\alpha$. The coefficient $C$ is determined by simulation.

Daigle and Landford \cite{R11:FractalKa-JSAC99} used the AMS model to study
packet voice communication systems. Tucker \cite{tucker} also used a similar
model, but the server capacity was defined as an integer number of information
units and the model had finite buffer size. Results on the buffer content
distribution, cell loss and delay were given, along with simulation comparisons.

The feature\ common to the above models is that they all assume a continuous
time distribution for the \textbf{on-off} intervals. On the other hand, Li
\cite{Li} introduced a discrete time model with a finite number of
\textbf{on-off }sources and geometric distributions for the \textbf{on-off
}intervals. He assumed that in one time unit only one \textbf{on-off }source
can change state. The channel capacity was assumed to be an integer number of
sources and the buffer size could be either zero (burst switching-clipping
case) or infinity (packet switching case).

In \cite{MR89h:60166} Mitra shows how to deal with states $i\in\mathcal{N}$
for which $\mathrm{d}(i)=0$. He also generalizes the AMS model in two ways.
Again, he considers a buffer which receives input from $N$ i.i.d. exponential
sources, called producers. However, now the output is also assumed to be
Markov-modulated, by $M$ i.i.d. exponential consumers. A second generalization
is that he considers the case of finite buffer capacity.

In the manufacturing literature Wijngaard \cite{wijngaard} has analyzed the
model treated in \cite{MR89h:60166} for the case of one producing and one
consuming machine. Sevast{$^{\prime}$}yanov \cite{Sevastyanov} has given an
innovative approximation for a production line model of many stages, which
incorporates the exact solution of a pair of stages. All of this work assumes
finite capacity buffers. Another extension of the AMS model is presented in
\cite{MR97i:60133}, where the authors allow the possibility of two sources of
input, one slow source of a fluid and another source that generates packets.

Some birth-death fluid queues with finite $\mathcal{N}$ \ have also been
proposed as approximative models, when the modulating process in the original
model is a multidimensional Markov chain, see \cite{MR91k:94028}. In
\cite{kogan} the authors allow the modulating process to be an arbitrary
birth-death process on a finite state space, but require the drift function
\textrm{d}$(i)$ to have a particular sign structure. An analysis without any
restriction on the sign structure of the drift vector has been given in the
survey paper \cite{birth}.

Until recently, authors have mostly focused on \textit{stationary}
probabilities. The transient analysis of \ stochastic fluid models, i.e., the
analysis of the fluid distribution in the buffer at an arbitrary time $t$, is
a complex and computationally intensive task. There are only a few papers
dedicated to the transient behavior of fluid queues with Markov input. Narayan
and Kulkarni \cite{MR1962652} derive explicit expressions for the Laplace
transform of the joint distribution of the first time the buffer becomes empty
and the state of the Markov process at that time.

The Laplace transform has been often used to evaluate the transient behavior
of fluid flow models. In \cite{MR95m:60139} Ren and Kobayashi studied the
transient distribution of the buffer content for the AMS model. The same
authors deal with the case of multiple types of inputs in \cite{Ren}. These
studies have been extended to the Markov modulated input rate model by Tanaka,
et. al., in \cite{tanaka}.

Another approach to transient behavior is due to Sericola \cite{sericola}.
This leads to a numerically stable recursive method for computing transient
and first passage distributions, and applies to general rate and drift functions.

The equilibrium probability that the buffer content exceeds $x,$ $\Pr\left[
X(\infty)>x\right]  ,$ in the AMS model has been studied asymptotically, as
$N\rightarrow\infty,$ by Weiss \cite{MR87k:60243} and Morrison
\cite{MR90c:60067}. The exact expression in \cite{MR84a:68020} is convenient
for calculating this probability for moderate values of $N.$ However, when $N$
is large it is useful to have asymptotic approximations$.$

Weiss \cite{MR87k:60243} (see also \cite{MR96i:60029}) analyzed the large $N$
asymptotics using the theory of large deviations. This yields an approximation
of the form $\Pr\left[  X(\infty)>NB\right]  \approx\exp\left[  -\mathrm{I}%
(B)N\right]  ,$ where $\mathrm{I}(B)$ is characterized as a solution of a
variational problem. Overflow asymptotics of this type for more general
Markov-modulated fluid models are derived in \cite{MR97k:60257}. Duffield
\cite{duffield} considers time-dependent asymptotics for the overflow probability.

Morrison \cite{MR90c:60067} investigated the AMS model starting from the exact
expression in \cite{MR84a:68020}, using the Euler-Maclaurin formula. He
constructed different asymptotic approximations to $\Pr\left[  X(\infty
)>x\right]  $ for the three scales $x=O(1/N),\quad x=O(1)$ \ and \ $x=O(N).$

The heavy traffic analysis of the AMS model, where $N\rightarrow\infty$ and
$\frac{c}{N}=\frac{\lambda}{\lambda+1}+O\left(  N^{-1/2}\right)  ,$ was
studied by Knessl and Morrison \cite{MR91m:60183}. In this limit, computing
the joint steady-state distribution of the buffer content and the number of
\textbf{on} sources can be reduced to solving a partial differential equation,
which is forward parabolic in a part of the domain and backward parabolic in
the remaining portion.

In this paper we analyze the joint distribution of the buffer content and
number of \textbf{on} sources in the AMS model asymptotically, as the number
of sources $N\rightarrow\infty.$ Our results extend those obtained in
\cite{MR90c:60067}.

In Section 2 we describe the AMS model in more detail and in Section 3 we
derive an integral representation of the solution, from which the spectral
representation obtained by the authors in \cite{MR84a:68020} can be easily
derived. In Sections 4-9 we obtain asymptotic expansions from the exact
solution by using the saddle point method. We have to consider eleven relevant
regions of the two-dimensional state space. In Section 10 we summarize our
results and numerically compare our approximations with the exact solution. In
Section 11 we analyze the conditional and marginal distributions.

\section{The Anick-Mitra-Sondhi (AMS) model}

In this model a data-handling switch receives messages from $N$ mutually
independent sources, which independently and asynchronously alternate between
\textbf{on}\ and \textbf{off}\ states. The \textbf{on}\ and \textbf{off}%
\ periods are exponentially distributed for each source. These two
distributions, while not necessarily identical, are common to all sources.

The unit of time is selected to be the average \textbf{on}\ period and the
average \textbf{off}\ period is denoted by $\frac{1}{\lambda}$. Thus, if we
denote by On (Off) the random variables measuring the elapsed time a source is
in the \textbf{on} (\textbf{off) }state, we have%
\[
\Pr\left[  \text{On}\leq T\right]  =\int\limits_{0}^{T}e^{-s}ds,\quad
\Pr\left[  \text{Off}\leq T\right]  =\int\limits_{0}^{T}\lambda e^{-\lambda
s}ds
\]
and%
\begin{equation}
\mathrm{E}\left[  \Pr\left(  \text{a source is \textbf{on}}\right)  \right]
=\frac{\lambda}{\lambda+1}. \label{EisOn}%
\end{equation}
An \textbf{on}\ source will generate packets of information at a uniform rate
that we assume, without loss of generality, to be $1$ packet per unit of time.
While a source is \textbf{off, }it generates no packets. Thus, when $k$
sources are \textbf{on}\ simultaneously, the instantaneous receiving rate at
the switch is $k$. The switch stores or buffers the incoming information that
is in excess of the maximum transmission rate $c$ of an output channel (thus,
$c$ is also the ratio of the output channel capacity to an \textbf{on}%
\ source's transmission rate). We assume that $c$ is not an integer.

Let $Z(t)$ be the number of \textbf{on}\ sources at time $t$, and $X(t)$ be
the buffer content. As long as the buffer is not empty, the instantaneous rate
of change of the buffer content is $Z(t)-c$. Once the buffer is empty, it
remains so as long as $Z(t)<c$. We assume that the buffer capacity is
infinite, so we need the following stability condition
\begin{equation}
\frac{\lambda}{\lambda+1}<\gamma<1,\quad\gamma=\frac{c}{N} \label{stability}%
\end{equation}
to ensure the existence of a steady state distribution. Note that $c<N$
implies that the buffer may be non-empty and $c>\frac{\lambda}{\lambda+1}N$
means that the output rate exceeds the average input rate.

Here we consider only the equilibrium probabilities, defined by%
\[
F_{k}(x)=\underset{t\rightarrow\infty}{\lim}\Pr\left[  Z(t)=k,\quad X(t)\leq
x\right]  ,\quad0\leq k\leq N,\quad x>0.
\]
The balance equations are \cite{MR84a:68020}
\begin{equation}
(k-c)\frac{\partial F_{k}}{\partial x}=\lambda\left[  N-(k-1)\right]
\,\,F_{k-1}+(k+1)\,\,F_{k+1}-\left[  \lambda(N-k)+k\right]  \,F_{k},\quad0\leq
k\leq N, \label{AMS}%
\end{equation}
with%
\begin{equation}
F_{k}(x)\equiv0,\quad k\notin\lbrack0,N]. \label{boundary}%
\end{equation}
Moreover, if the number of sources \textbf{on}\ at any time exceeds $c$, then
the buffer content increases and the buffer cannot stay empty. Hence,%
\begin{equation}
F_{k}(0)=0,\quad\left\lfloor c\right\rfloor +1\leq k\leq N. \label{zero}%
\end{equation}
Also,%
\begin{equation}
F_{k}(\infty)=\frac{1}{(1+\lambda)^{N}}\binom{N}{k}\lambda^{k},\quad0\leq
k\leq N \label{inf}%
\end{equation}
since $F_{k}(\infty)$ is the probability that $k$ out of $N$ sources are
\textbf{on}\ simultaneously and from (\ref{EisOn}) this follows a binomial
distribution with parameter $\frac{\lambda}{\lambda+1}$.

In \cite{MR84a:68020} Anick, Mitra and Sondhi derived the spectral
representation of the solution to (\ref{AMS})-(\ref{inf}). They wrote the
equation (\ref{AMS}) in matrix form%
\begin{equation}
\mathbf{DF}^{\prime}(x)=\mathbf{M~F(}x\mathbf{),}\quad x\geq0 \label{matrix}%
\end{equation}
where
\[
\mathbf{F(}x\mathbf{)=}\left[  F_{0}(x),F_{1}(x),\ldots,F_{N}(x)\right]
^{T},
\]
$\mathbf{D}$ is an $\left(  N+1\right)  \times\left(  N+1\right)  $ diagonal
matrix with entries%
\[
\mathbf{D}_{ij}=\left(  j-c\right)  \delta_{ij},\quad0\leq i,j\leq N
\]
and $\mathbf{M}$ is an $\left(  N+1\right)  \times\left(  N+1\right)  $
tridiagonal matrix with entries
\[
\mathbf{M}_{ij}=\left\{
\begin{array}
[c]{c}%
j,\quad j=i+1,\quad0\leq i\leq N-1\\
-\left[  \lambda(N-j)+j\right]  ,\quad j=i,\quad0\leq i\leq N\quad\\
\lambda\left(  N-j\right)  ,\quad j=i-1,\quad1\leq i\leq N\\
0,\quad\left\vert i-j\right\vert \geq2
\end{array}
\right.  .
\]
Using (\ref{matrix}), they represented the solution $\mathbf{F(}x\mathbf{)}$
as%
\begin{equation}
\mathbf{F(}x\mathbf{)=F(}\infty\mathbf{)+}%
{\displaystyle\sum\limits_{i=0}^{N-\left\lfloor c\right\rfloor -1}}
a_{i}e^{z_{i}x}\mathbf{\phi}_{i} \label{spectral0}%
\end{equation}
where $z_{i}$ is a negative eigenvalue and $\mathbf{\phi}_{i}$ is the
associated right eigenvector of the matrix $\mathbf{D}^{-1}\mathbf{M}$%
\[
\mathbf{D}^{-1}\mathbf{M\phi}_{i}=z_{i}\mathbf{\phi}_{i}.
\]
They showed that (\ref{AMS}) and (\ref{zero}) imply conditions on the
derivatives of $F_{N}(x)$ at $x=0,$
\begin{equation}
\frac{d^{n}F_{N}}{dx^{n}}(0)=0,\quad0\leq n\leq N-\left\lfloor c\right\rfloor
-1. \label{derivati}%
\end{equation}
From (\ref{derivati}) they found the coefficients $a_{i}$ in (\ref{spectral0})
with the help of Vandermonde determinants. Asymptotic and numerical results
for the overflow probability of a predefined buffer backlog were also presented.

\section{Exact solution}

Here we shall solve the problem (\ref{AMS})-(\ref{inf}) by a different method
than that used in \cite{MR84a:68020}. This will lead to an alternate form of
the solution that shall prove useful for asymptotic analysis.

Looking for a solution to (\ref{AMS}) of the form%
\[
F_{k}(x)=e^{\theta x}h_{k}(\theta)
\]
we get the following difference equation for $h_{k}(\theta)$%
\begin{equation}
\lambda\left(  N-k+1\right)  \,\,h_{k-1}+(k+1)\,\,h_{k+1}-\left[
\lambda(N-k)+k+(k-c)\theta\right]  \,h_{k}=0,\quad0\leq k\leq N. \label{diff}%
\end{equation}
To solve (\ref{diff}) we represent $h_{k}(\theta)$ as an integral
\begin{equation}
h_{k}(\theta)=\frac{1}{2\pi i}%
{\displaystyle\oint\limits_{\mathcal{C}}}
\frac{H(\theta,w)}{w^{N-k+1}}dw \label{h}%
\end{equation}
where $\mathcal{C}$ is an small loop around the origin in the complex plane
and $H(\theta,w)$ is an analytic function of $w$ inside $\mathcal{C}$. Using
(\ref{h}) in (\ref{diff}) and integrating by parts gives an ODE for
$H(\theta,w)$%
\begin{equation}
\frac{\partial H}{\partial w}=\frac{(w+\gamma\theta-1-\theta)N}{w^{2}%
+(\lambda-\theta-1)w-\lambda}H \label{diffH}%
\end{equation}
whose solution is%
\begin{equation}
H(\theta,w)=\left[  1-R_{1}(\theta)w\right]  ^{NV(\theta)}\left[
1+R_{2}(\theta)w\right]  ^{N\left[  1-V(\theta)\right]  } \label{H}%
\end{equation}
with%
\begin{equation}
V(\theta)=\frac{1}{2}\left[  1+\frac{\left(  2\gamma-1\right)  \theta
-\lambda-1}{\Delta(\theta)}\right]  ,\quad\Delta(\theta)=\sqrt{(\theta
+1-\lambda)^{2}+4\lambda} \label{V}%
\end{equation}%
\begin{equation}
R_{1}(\theta)=\frac{1}{2\lambda}\left[  \Delta(\theta)+(\lambda-1-\theta
)\right]  ,\quad R_{2}(\theta)=\frac{1}{2\lambda}\left[  \Delta(\theta
)-(\lambda-1-\theta)\right]  . \label{R12}%
\end{equation}
Using the binomial theorem we can also represent $h_{k}(\theta)$ as%
\begin{equation}
h_{k}(\theta)=%
{\displaystyle\sum\limits_{i=0}^{N-k}}
\binom{NV(\theta)}{i}\binom{N\left[  1-V(\theta)\right]  }{N-k-i}\left[
-R_{1}(\theta)\right]  ^{i}\left[  R_{2}(\theta)\right]  ^{N-k-i}
\label{binomial}%
\end{equation}
which shows that $h_{k}(\theta)$ is a polynomial in $\theta$ of degree $N-k,$
\begin{align*}
h_{N}(\theta)  &  =1,\quad h_{N-1}(\theta)=\frac{1+\left(  1-\gamma\right)
\theta}{\lambda}N,\\
h_{N-2}(\theta)  &  =\frac{\left(  \gamma-1\right)  \theta^{2}+\left(
\phi-2\right)  \theta-1}{2\lambda^{2}}N+\frac{\left[  1+\left(  1-\gamma
\right)  \theta\right]  ^{2}}{2\lambda^{2}}N^{2},\\
&  \vdots\\
h_{N-k}(\theta)  &  =c_{1}(\theta)N+\cdots+c_{k-1}(\theta)N^{k-1}%
+\frac{\left[  1+\left(  1-\gamma\right)  \theta\right]  ^{k}}{\lambda^{k}%
k!}N^{k},\quad1\leq k\leq N.
\end{align*}
The polynomials $h_{k}(\theta)$ are related to the Krawtchouk polynomials
\cite{MR51:8724}. The functions $c_{k}(\theta)$ are entire (polynomial)
functions of $\theta,$ that may be computed from (\ref{binomial}).

Since $H$ is analytic inside $\mathcal{C}$, we have%
\[
h_{N+1}(\theta)=\frac{1}{2\pi i}%
{\displaystyle\oint\limits_{\mathcal{C}}}
H(\theta,w)dw=0.
\]
Imposing the boundary condition (\ref{boundary}) yields $h_{k}(\theta)=0$ for
$k\leq-1.$ Thus,
\[
\frac{\partial^{k}H}{\partial w^{k}}\left(  \theta,0\right)  =0,\quad k\geq
N+1
\]
and therefore $H(\theta,w)$ must be a polynomial in $w$ of degree $\leq N$. In
view of (\ref{H}) we must have
\begin{equation}
V(\theta_{j})=\frac{j}{N},\quad0\leq j\leq N \label{Veig}%
\end{equation}
which gives the equation for the eigenvalues $\theta_{j}.$ Since from
(\ref{inf}) we need $F_{k}(x)\rightarrow$ $F_{k}(\infty)$ for $x\rightarrow
\infty,$ we consider only the negative eigenvalues. Solving for $\theta_{j}$
in (\ref{Veig}) we get
\begin{equation}
\theta_{j}=-\sigma\left(  \frac{j}{N}\right)  <0,\quad0\leq j\leq
N-\left\lfloor c\right\rfloor -1 \label{sigmaeig}%
\end{equation}
where%
\begin{equation}
\sigma(x)=\frac{\rho+2(\lambda-1)(1-x)x+(1-2x)\sqrt{\rho^{2}+4\lambda(1-x)x}%
}{2(\gamma-x)(1-\gamma-x)} \label{sigma}%
\end{equation}
and%
\begin{equation}
\rho=\gamma-\lambda+\lambda\gamma\label{rho}%
\end{equation}
with $0<\rho<1$ from (\ref{stability}). Combining the above results and using
(\ref{inf}), we find the spectral representation of the solution%
\begin{equation}
F_{k}(x)=\frac{1}{(1+\lambda)^{N}}\binom{N}{k}\lambda^{k}+%
{\displaystyle\sum\limits_{j=0}^{N-\left\lfloor c\right\rfloor -1}}
a_{j}e^{\theta_{j}x}h_{k}(\theta_{j}) \label{spectral}%
\end{equation}
where the constants $a_{j}$ are to be determined. To find $a_{j}$ we employ
arguments different from those in \cite{MR84a:68020}$.$

Taking the Laplace transform%
\[
\mathcal{F}_{k}(\vartheta)=%
{\displaystyle\int\limits_{0}^{\infty}}
e^{-\vartheta x}F_{k}(x)dx
\]
in (\ref{AMS}) and using the boundary conditions (\ref{zero}) for
$k\geq\left\lfloor c\right\rfloor +1$ gives%
\begin{equation}
(k-c)\vartheta\mathcal{F}_{k}=\lambda\left[  N-(k-1)\right]  \,\mathcal{F}%
_{k-1}+(k+1)\,\,\mathcal{F}_{k+1}-\left[  \lambda(N-k)+k\right]
\,\mathcal{F}_{k},\quad\left\lfloor c\right\rfloor +1\leq k\leq N.
\label{laplacediff}%
\end{equation}
Equation (\ref{laplacediff}) has two independent solutions, $h_{k}(\vartheta)$
and $h_{k}^{\ast}(\vartheta),$ with $h_{k}(\vartheta)$ given by (\ref{h}) and
\begin{equation}
h_{k}^{\ast}(\vartheta)=\frac{1}{2\pi i}%
{\displaystyle\oint\limits_{\mathcal{C}}}
H^{\ast}(\vartheta,w)\frac{1}{w^{k+1}}dw \label{h2}%
\end{equation}
where $\mathcal{C}$ is a small loop around the origin and
\[
H^{\ast}(\vartheta,w)=\left[  1-\lambda R_{2}(\vartheta)w\right]
^{NV(\theta)}\left[  1+\lambda R_{1}(\vartheta)w\right]  ^{N\left[
1-V(\vartheta)\right]  }.
\]
Using the binomial theorem we get from (\ref{h2})
\begin{equation}
h_{k}^{\ast}(\vartheta)=\lambda^{k}%
{\displaystyle\sum\limits_{i=0}^{k}}
\binom{NV(\vartheta)}{i}\binom{N\left[  1-V(\vartheta)\right]  }{k-i}\left[
-R_{2}(\vartheta)\right]  ^{i}\left[  R_{1}(\vartheta)\right]  ^{k-i}.
\label{binomial2}%
\end{equation}
From (\ref{binomial2}) we have%
\begin{align*}
h_{0}^{\ast}(\vartheta)  &  =1,\quad h_{1}^{\ast}(\vartheta)=\left(
\lambda-\gamma\vartheta\right)  N,\\
h_{2}^{\ast}(\vartheta)  &  =\frac{1}{2}\left[  -\lambda^{2}+\left(
\gamma\lambda+\lambda-\gamma\right)  \vartheta-\gamma\vartheta^{2}\right]
N+\frac{\left(  \lambda-\gamma\vartheta\right)  ^{2}}{2}N^{2},\\
&  \vdots\\
h_{k}^{\ast}(\vartheta)  &  =c_{1}^{\ast}(\vartheta)N+\cdots+c_{k-1}^{\ast
}(\vartheta)N^{k-1}+\frac{\left(  \lambda-\gamma\vartheta\right)  ^{k}}%
{k!}N^{k},\quad1\leq k\leq N
\end{align*}
where the $c_{k}^{\ast}(\vartheta)$ are polynomials that can be identified
from (\ref{binomial2}) using (\ref{V}) and (\ref{R12}). Therefore, the general
solution to (\ref{laplacediff}) is given by%
\[
\mathcal{F}_{k}(\vartheta)=A_{1}(\vartheta)h_{k}(\vartheta)+A_{2}\left(
\vartheta\right)  h_{k}^{\ast}(\vartheta),\quad\left\lfloor c\right\rfloor
+1\leq k\leq N
\]
with $A_{1}(\vartheta),$ $A_{2}\left(  \vartheta\right)  $ still undetermined.
We note that $h_{k}$ and $h_{k}^{\ast}$ are \textit{entire} functions of
$\vartheta.$ Since $h_{N+1}\left(  \vartheta\right)  =0$ but $h_{N+1}^{\ast
}(\vartheta)\neq0,$ we see from (\ref{boundary}) that we need $A_{2}\left(
\vartheta\right)  =0,$ so that%
\begin{equation}
\mathcal{F}_{k}(\vartheta)=A_{1}(\vartheta)h_{k}(\vartheta). \label{laplace3}%
\end{equation}

Taking the Laplace transform in (\ref{spectral}) yields%
\begin{equation}
\mathcal{F}_{k}(\vartheta)=\frac{1}{(1+\lambda)^{N}}\binom{N}{k}\lambda
^{k}\frac{1}{\vartheta}+%
{\displaystyle\sum\limits_{j=0}^{N-\left\lfloor c\right\rfloor -1}}
a_{j}h_{k}(\theta_{j})\frac{1}{\vartheta-\theta_{j}}, \label{laplace}%
\end{equation}
from which we see that the only singularities of $\mathcal{F}_{k}(\vartheta)$
are simple poles at $\vartheta=0$ \ and \ $\vartheta=\theta_{j},\quad0\leq
j\leq N-\left\lfloor c\right\rfloor -1.$ \ Since $h_{k}(\vartheta)$ is entire,
$A_{1}(\vartheta)$ must have poles at these points also. We thus rewrite
(\ref{laplace3}) as
\[
\mathcal{F}_{k}(\vartheta)=B(\vartheta)\frac{1}{\vartheta}\left[
{\displaystyle\prod\limits_{j=0}^{N-\left\lfloor c\right\rfloor -1}}
\frac{1}{\theta_{j}-\vartheta}\right]  h_{k}(\vartheta)
\]
where $B(\vartheta)$\ is an entire function. Using the inversion formula, we
can represent \ $F_{k}(x)$ \ in the form%
\begin{equation}
F_{k}(x)=\frac{1}{2\pi i}%
{\displaystyle\int\limits_{\mathrm{Br}^{+}}}
e^{x\vartheta}B(\vartheta)\frac{1}{\vartheta}\left[
{\displaystyle\prod\limits_{j=0}^{N-\left\lfloor c\right\rfloor -1}}
\frac{1}{\theta_{j}-\vartheta}\right]  h_{k}(\vartheta)d\vartheta\label{F1}%
\end{equation}
where $\mathrm{Br}^{+}$ is a vertical contour on which $\operatorname{Re}%
(\vartheta)>0$. \ From (\ref{laplace}) the residue of $\mathcal{F}%
_{k}(\vartheta)$ at \ $\vartheta=0$ \ is equal to%
\begin{equation}
\frac{1}{(1+\lambda)^{N}}\binom{N}{k}\lambda^{k} \label{Res0}%
\end{equation}
and from (\ref{H}) we know that%
\[
H(0,w)=\left(  1+\frac{w}{\lambda}\right)  ^{N}%
\]
which implies that%
\begin{equation}
h_{k}(0)=\binom{N}{k}\lambda^{k-N}. \label{h(0)}%
\end{equation}
Thus, from (\ref{F1}), we get%
\[
B(0)=\left(  \frac{\lambda}{\lambda+1}\right)  ^{N}%
{\displaystyle\prod\limits_{j=0}^{N-\left\lfloor c\right\rfloor -1}}
\theta_{j}.
\]
We determine $B(\vartheta)$ by examining the limit $\vartheta\rightarrow
\infty$, which correspond to considering the boundary condition at $x=0.$

From (\ref{V}) and (\ref{R12}) we can easily obtain the asymptotic expressions%
\begin{align}
V(\vartheta)  &  \sim\gamma-\phi\frac{1}{\vartheta},\quad\vartheta
\rightarrow\infty\nonumber\\
R_{1}(\vartheta)  &  \sim\frac{1}{\vartheta}+(\lambda-1)\frac{1}{\vartheta
^{2}},\quad\vartheta\rightarrow\infty\label{VR1}\\
R_{2}(\vartheta)  &  \sim\frac{1}{\lambda}\vartheta+\frac{1-\lambda}{\lambda
}+\frac{1}{\vartheta},\quad\vartheta\rightarrow\infty.\nonumber
\end{align}
Setting $w=\frac{u}{\vartheta}$ in (\ref{h}) and using (\ref{VR1}) yields%
\begin{align}
h_{k}(\vartheta)  &  \sim\vartheta^{N-k}\frac{1}{2\pi i}%
{\displaystyle\oint\limits_{\mathcal{C}}}
\left(  1+\frac{u}{\lambda}\right)  ^{(1-\gamma)N}\frac{du}{u^{N-k+1}%
}\label{hinf}\\
&  =\left(  \frac{\vartheta}{\lambda}\right)  ^{N-k}\binom{\left(
1-\gamma\right)  N}{N-k},\quad\vartheta\rightarrow\infty.\nonumber
\end{align}
We also have%
\begin{equation}%
{\displaystyle\prod\limits_{j=0}^{N-\left\lfloor c\right\rfloor -1}}
\frac{1}{\theta_{j}-\vartheta}\sim\frac{(-1)^{N-\left\lfloor c\right\rfloor }%
}{\vartheta^{N-\left\lfloor c\right\rfloor }},\quad\vartheta\rightarrow\infty.
\label{Pinf}%
\end{equation}
Combining (\ref{hinf}) and (\ref{Pinf}) we get%
\begin{equation}
\mathcal{F}_{k}(\vartheta)\sim B(\vartheta)\frac{(-1)^{N-\left\lfloor
c\right\rfloor }}{\lambda^{N-k}}\binom{\left(  1-\gamma\right)  N}%
{N-k}\vartheta^{\left\lfloor c\right\rfloor -k-1},\quad\vartheta
\rightarrow\infty. \label{Lasym}%
\end{equation}
Since (\ref{zero}) implies that
\[
\underset{\vartheta\rightarrow\infty}{\lim}\left[  \vartheta\mathcal{F}%
_{k}(\vartheta)\right]  =0,\quad\left\lfloor c\right\rfloor +1\leq k\leq N,
\]
we obtain%
\[
B(\vartheta)=o\left(  \vartheta^{k-\left\lfloor c\right\rfloor }\right)
,\quad\vartheta\rightarrow\infty,\quad\left\lfloor c\right\rfloor +1\leq k\leq
N.
\]
Setting $k=\left\lfloor c\right\rfloor +1$ we see that $B(\vartheta)$ is an
entire function that is $o\left(  \vartheta\right)  $ as $\vartheta
\rightarrow\infty.$ By the generalized Liouville theorem \cite{MR7:53e},
$B(\vartheta)$ must be a constant%
\[
B(\vartheta)=B(0)=\left(  \frac{\lambda}{\lambda+1}\right)  ^{N}%
{\displaystyle\prod\limits_{j=0}^{N-\left\lfloor c\right\rfloor -1}}
\theta_{j}.
\]
Therefore we can write%
\begin{equation}
F_{k}(x)=\left(  \frac{\lambda}{\lambda+1}\right)  ^{N}\frac{1}{2\pi i}%
{\displaystyle\int\limits_{\mathrm{Br}^{+}}}
e^{x\theta}\frac{1}{\vartheta}\left[
{\displaystyle\prod\limits_{j=0}^{N-\left\lfloor c\right\rfloor -1}}
\frac{\theta_{j}}{\theta_{j}-\vartheta}\right]  h_{k}(\vartheta)d\vartheta.
\label{Fexact}%
\end{equation}

Closing $\mathrm{Br}^{+}$ in the half-plane $\operatorname{Re}(\vartheta
)\leq0$ we recover the spectral representation%
\[
F_{k}(x)=\frac{1}{(1+\lambda)^{N}}\binom{N}{k}\lambda^{k}+%
{\displaystyle\sum\limits_{j=0}^{N-\left\lfloor c\right\rfloor -1}}
a_{j}e^{\theta_{j}x}h_{k}(\theta_{j})
\]
with
\[
a_{j}=-\left(  \frac{\lambda}{\lambda+1}\right)  ^{N}%
{\displaystyle\prod\limits_{\substack{i=0\\i\neq j}}^{N-\left\lfloor
c\right\rfloor -1}}
\frac{\theta_{i}}{\theta_{i}-\theta_{j}}.
\]
In conclusion, we have proved the following theorem.

\begin{theorem}
The solution of (\ref{AMS})-(\ref{inf}) is given by%
\begin{equation}
F_{k}(x)=\left(  \frac{\lambda}{\lambda+1}\right)  ^{N}\frac{1}{2\pi i}%
{\displaystyle\int\limits_{\mathrm{Br}^{+}}}
e^{x\vartheta}\frac{1}{\vartheta}\left[
{\displaystyle\prod\limits_{j=0}^{N-\left\lfloor c\right\rfloor -1}}
\frac{\theta_{j}}{\theta_{j}-\vartheta}\right]  h_{k}(\vartheta)d\vartheta
\label{Exact}%
\end{equation}
where $\mathrm{Br}^{+}$ is a vertical contour on which $\operatorname{Re}%
(\vartheta)>0,$
\[
h_{k}(\vartheta)=\frac{1}{2\pi i}%
{\displaystyle\oint\limits_{\mathcal{C}}}
\frac{H(\vartheta,w)}{w^{N-k+1}}dw,
\]
$\mathcal{C}$ \ is a small loop around the origin in the complex $w-$plane
\[
H(\vartheta,w)=\left[  1-R_{1}(\vartheta)w\right]  ^{NV(\vartheta)}\left[
1+R_{2}(\vartheta)w\right]  ^{N\left[  1-V(\vartheta)\right]  },
\]%
\[
V(\vartheta)=\frac{1}{2}\left[  1+\frac{\left(  2\gamma-1\right)
\vartheta-\lambda-1}{\Delta(\vartheta)}\right]  ,\quad\Delta(\vartheta
)=\sqrt{(\vartheta+1-\lambda)^{2}+4\lambda},
\]%
\[
R_{1}(\vartheta)=\frac{1}{2\lambda}\left[  \Delta(\vartheta)+(\lambda
-1-\vartheta)\right]  ,\quad R_{2}(\vartheta)=\frac{1}{2\lambda}\left[
\Delta(\vartheta)-(\lambda-1-\vartheta)\right]  ,
\]%
\[
\theta_{j}=-\sigma\left(  \frac{j}{N}\right)  <0,\quad0\leq j\leq
N-\left\lfloor c\right\rfloor -1,
\]
and%
\[
\sigma(x)=\frac{\rho+2(\lambda-1)(1-x)x+(1-2x)\sqrt{\rho^{2}+4\lambda(1-x)x}%
}{2(\gamma-x)(1-\gamma-x)}.
\]

\end{theorem}

We have thus re-derived, using different arguments, the spectral
representation%
\begin{equation}
F_{k}(x)=\frac{1}{(1+\lambda)^{N}}\binom{N}{k}\lambda^{k}-\left(
\frac{\lambda}{\lambda+1}\right)  ^{N}%
{\displaystyle\sum\limits_{j=0}^{N-\left\lfloor c\right\rfloor -1}}
\left[
{\displaystyle\prod\limits_{\substack{i=0\\i\neq j}}^{N-\left\lfloor
c\right\rfloor -1}}
\frac{\theta_{i}}{\theta_{i}-\theta_{j}}\right]  e^{\theta_{j}x}h_{k}%
(\theta_{j}) \label{spectral1}%
\end{equation}
obtained by the authors in \cite{MR84a:68020}. We also derived the integral
representation (\ref{Exact}) that will prove more useful than (\ref{spectral1}%
) for asymptotic analysis. We note that the coefficients in the sum in
(\ref{spectral1}) alternate in sign and this makes it very difficult to obtain
asymptotic results for $F_{k}(x)$ from the spectral representation.

We study asymptotic properties of $F_{k}(x)$ as $N\rightarrow\infty$ for
various ranges of $k$ and $x.$ We shall use (\ref{Exact}), approximate the
integrand for $N$ large and then evaluate the integral over $\mathrm{Br}^{+}$
asymptotically. We will introduce the scaled variables $\left(  y,z\right)
=\left(  \frac{x}{N},\frac{k}{N}\right)  $ and consider the limit
$N\rightarrow\infty,$ with $y$ and $z$ fixed. In Section \ref{main approx} we
will derive asymptotic approximations to $F_{k}(x)$ that are valid in most of
the strip $\mathfrak{D=}\left\{  \left(  y,z\right)  :y\geq0,\quad0\leq
z\leq1\right\}  ,$ with the exception of a few boundary regions and one curve
in $\overset{\circ}{\mathfrak{D}}$ (the interior of $\mathfrak{D).}$These
other ranges are considered in Sections \ref{section corner} -
\ref{section y=0}.

\section{The main approximation}

\label{main approx}

We first evaluate the product in (\ref{Exact}) for $N\rightarrow\infty$.

\begin{lemma}
\label{lemma2}Let \
\[
P(\vartheta)=%
{\displaystyle\prod\limits_{j=0}^{N-\left\lfloor c\right\rfloor -1}}
\frac{\theta_{j}}{\theta_{j}-\vartheta}.
\]
Then, for \ $N\rightarrow\infty\ $and a fixed $\vartheta>\theta_{0},$ we have%
\begin{equation}
P(\vartheta)\sim\widetilde{P}(\vartheta;N)\equiv\sqrt{\frac{-\theta_{0}%
}{\vartheta-\theta_{0}}}\exp\left[  N\mu(\vartheta)\right]  ,\quad\label{P2}%
\end{equation}
with%
\begin{gather}
\mu(\vartheta)=-\frac{1}{2}\ln\left[  \gamma(1-\gamma)\vartheta+\rho\right]
\label{mu}\\
+\frac{\vartheta\left(  1-2\gamma\right)  +\lambda+1}{2\Delta(\vartheta)}%
\ln\left[  \frac{\left[  \gamma(1-\gamma)\vartheta+\rho\right]  \left[
\lambda-1-\vartheta+\Delta(\vartheta)\right]  }{\left(  \lambda-1-\vartheta
\right)  \rho+\left(  \lambda+1\right)  ^{2}\gamma(1-\gamma)+\left(
1-\lambda\right)  \gamma(1-\gamma)\vartheta+\Delta(\vartheta)\delta}\right]
\nonumber
\end{gather}
and%
\begin{equation}
\theta_{0}=-\frac{\rho}{\gamma\left(  1-\gamma\right)  }<0,\quad\delta=\left(
1-\gamma\right)  ^{2}\lambda+\gamma^{2}>0, \label{delta}%
\end{equation}%
\[
\rho=\gamma\lambda+\gamma-\lambda.
\]
We note from (\ref{mu}) that $\mu(0)=0.$
\end{lemma}

\begin{proof}
Since we assumed that $c$ is not an integer, we introduce the fractional part
of $c$ defined by%
\begin{equation}
\alpha=c-\left\lfloor c\right\rfloor ,\quad0<\alpha<1 \label{alpha}%
\end{equation}
which allows us to write $N-\left\lfloor c\right\rfloor $ as $(1-\gamma
)N+\alpha.$ Since
\[
\ln\left[  P(\vartheta)\right]  =\ln\left[
{\displaystyle\prod\limits_{k=0}^{(1-\gamma)N+\alpha-1}}
\frac{\sigma\left(  \frac{k}{N}\right)  }{\vartheta+\sigma\left(  \frac{k}%
{N}\right)  }\right]  =%
{\displaystyle\sum\limits_{k=0}^{(1-\gamma)N+\alpha-1}}
\ln\left[  \frac{\sigma\left(  \frac{k}{N}\right)  }{\vartheta+\sigma\left(
\frac{k}{N}\right)  }\right]
\]
we can apply the Euler-Maclaurin formula \cite{MR80d:00030} to get%
\begin{align*}
\ln\left[  P(\vartheta)\right]   &  =\frac{1}{2}\left\{  \ln\left[
\frac{\sigma\left(  0\right)  }{\vartheta+\sigma\left(  0\right)  }\right]
+\ln\left[  \frac{\sigma\left(  1-\gamma+\frac{\alpha-1}{N}\right)
}{\vartheta+\sigma\left(  1-\gamma+\frac{\alpha-1}{N}\right)  }\right]
\right\} \\
&  +%
{\displaystyle\int\limits_{0}^{(1-\gamma)N+\alpha-1}}
\ln\left[  \frac{\sigma\left(  \frac{k}{N}\right)  }{\vartheta+\sigma\left(
\frac{k}{N}\right)  }\right]  dk+o(1).
\end{align*}
Changing variables from $k$ to $x=\frac{k}{N}$ we obtain%
\begin{equation}%
{\displaystyle\int\limits_{0}^{(1-\gamma)N+\alpha-1}}
\ln\left[  \frac{\sigma\left(  \frac{k}{N}\right)  }{\vartheta+\sigma\left(
\frac{k}{N}\right)  }\right]  dk=\frac{1}{\varepsilon}%
{\displaystyle\int\limits_{0}^{1-\gamma-(1-\alpha)\varepsilon}}
\ln\left[  \frac{\sigma\left(  x\right)  }{\vartheta+\sigma\left(  x\right)
}\right]  dx \label{I1.3}%
\end{equation}
with%
\begin{equation}
\varepsilon=\frac{1}{N}. \label{epsilon}%
\end{equation}
If we denote the indefinite integral in (\ref{I1.3}) by%
\[
\mathrm{I}_{1}(z)=%
{\displaystyle\int\limits^{z}}
\ln\left[  \frac{\sigma\left(  x\right)  }{\vartheta+\sigma\left(  x\right)
}\right]  dx
\]
then we have, as $\varepsilon\rightarrow0^{+},$
\begin{equation}%
{\displaystyle\int\limits_{0}^{1-\gamma-(1-\alpha)\varepsilon}}
\ln\left[  \frac{\sigma\left(  x\right)  }{\vartheta+\sigma\left(  x\right)
}\right]  dx\sim\mathrm{I}_{1}\left[  (1-\gamma)^{-}\right]  -(1-\alpha
)\mathrm{I}_{1}^{\prime}\left[  (1-\gamma)^{-}\right]  \varepsilon
-\mathrm{I}_{1}(0),\quad\label{44}%
\end{equation}
where $\mathrm{I}_{1}\left[  (1-\gamma)^{-}\right]  $ means $\underset
{z\ \uparrow\ 1-\gamma}{\lim}\mathrm{I}_{1}(z)$. From (\ref{sigma}) we get%
\begin{equation}
\sigma(0)=\frac{\rho}{\gamma(1-\gamma)},\quad\sigma\left(  1-\gamma
-\varepsilon\right)  \sim\frac{\phi}{\varepsilon},\quad\varepsilon
\rightarrow0^{+} \label{sigma(0)}%
\end{equation}
with%
\begin{equation}
\phi=\gamma+\lambda-\gamma\lambda. \label{phi}%
\end{equation}
Therefore,%
\begin{equation}
\mathrm{I}_{1}^{\prime}(1-\gamma)=\underset{x\ \uparrow\ 1-\gamma}{\lim}%
\ln\left[  \frac{\sigma\left(  x\right)  }{\vartheta+\sigma\left(  x\right)
}\right]  =\ln\left(  1\right)  =0. \label{limit1}%
\end{equation}

To find $\mathrm{I}_{1}(z),$ we use integration by parts%
\begin{equation}
\mathrm{I}_{1}(z)=-%
{\displaystyle\int\limits^{z}}
\ln\left[  \frac{\vartheta}{\sigma\left(  x\right)  }+1\right]  dx=-z\ln
\left[  \frac{\vartheta}{\sigma\left(  z\right)  }+1\right]  -\vartheta%
{\displaystyle\int\limits^{z}}
x\frac{\sigma^{\prime}\left(  x\right)  }{\left[  \vartheta+\sigma\left(
x\right)  \right]  \sigma\left(  x\right)  }dx. \label{I1.1}%
\end{equation}
Since (\ref{Veig}) implies that $V\left[  -\sigma(x)\right]  =x,$ we can write
(\ref{I1.1}) as%
\[
\mathrm{I}_{1}(z)=-z\ln\left[  \frac{\vartheta}{\sigma\left(  z\right)
}+1\right]  -\vartheta%
{\displaystyle\int\limits^{z}}
\frac{V\left[  -\sigma(x)\right]  }{\left[  \vartheta+\sigma\left(  x\right)
\right]  \sigma\left(  x\right)  }\sigma^{\prime}\left(  x\right)  dx
\]
or, changing the integrating variable from $x$ to $\sigma,$%
\[
\mathrm{I}_{1}(z)=-z\ln\left[  \frac{\vartheta}{\sigma\left(  z\right)
}+1\right]  -\vartheta%
{\displaystyle\int\limits^{\sigma(z)}}
\frac{V\left(  -\sigma\right)  }{\left(  \vartheta+\sigma\right)  \sigma
}d\sigma
\]%
\[
=-z\ln\left[  \frac{\vartheta}{\sigma\left(  z\right)  }+1\right]  -\vartheta%
{\displaystyle\int\limits^{\sigma(z)}}
\frac{1}{2}\left[  1-\frac{\left(  2\gamma-1\right)  \sigma+\lambda+1}%
{\sqrt{(-\sigma+1-\lambda)^{2}+4\lambda}}\right]  \frac{1}{\left(
\vartheta+\sigma\right)  \sigma}d\sigma
\]%
\begin{align}
&  =-z\ln\left[  \frac{\vartheta}{\sigma\left(  z\right)  }+1\right]
-\frac{\vartheta}{2}%
{\displaystyle\int\limits^{\sigma(z)}}
\frac{1}{\left(  \vartheta+\sigma\right)  \sigma}d\sigma\label{I1.4}\\
&  +\frac{\vartheta}{2}%
{\displaystyle\int\limits^{\sigma(z)}}
\frac{\left(  2\gamma-1\right)  \sigma+\lambda+1}{\sqrt{(-\sigma
+1-\lambda)^{2}+4\lambda}}\frac{1}{\left(  \vartheta+\sigma\right)  \sigma
}d\sigma.\nonumber
\end{align}
The first integral in (\ref{I1.4}) is
\[
\int\frac{1}{\left(  \vartheta+\sigma\right)  \sigma}d\sigma=\frac
{1}{\vartheta}\ln\left[  \frac{\sigma}{\sigma+\vartheta}\right]  .
\]
The second integral in (\ref{I1.4}) is quite complicated and after some
calculation we find that%
\begin{gather*}
\int\frac{\left(  2\gamma-1\right)  \sigma+\lambda+1}{\sqrt{(-\sigma
+1-\lambda)^{2}+4\lambda}}\frac{1}{\left(  \vartheta+\sigma\right)  \sigma
}d\sigma\\
=\frac{1}{2\vartheta}\ln\left\{  \frac{\left[  \left(  \lambda+1\right)
^{2}+\left(  \lambda-1\right)  \sigma-\left(  \lambda+1\right)  \Delta
(-\sigma)\right]  ^{2}}{4\lambda\sigma^{2}}\right\} \\
+\frac{\vartheta\left(  1-2\gamma\right)  +\lambda+1}{\vartheta\Delta
(\vartheta)}\ln\left[  2\frac{\left(  \lambda-1-\vartheta\right)
\sigma+\left(  \lambda+1\right)  ^{2}+\left(  1-\lambda\right)  \vartheta
+\Delta(\vartheta)\Delta(-\sigma)}{\vartheta+\sigma}\right]  .
\end{gather*}
Hence,
\begin{equation}
\mathrm{I}_{1}(z)=\left(  z-\frac{1}{2}\right)  \ln\left[  \frac{\sigma\left(
z\right)  }{\vartheta+\sigma\left(  z\right)  }\right]  +Q_{1}\left[
\sigma(z)\right]  \label{I1.2}%
\end{equation}
with%
\begin{align}
Q_{1}(\sigma)  &  =\frac{1}{4}\ln\left\{  \frac{\left[  \left(  \lambda
+1\right)  ^{2}+\left(  \lambda-1\right)  \sigma-\left(  \lambda+1\right)
\Delta(-\sigma)\right]  ^{2}}{4\lambda\sigma^{2}}\right\} \label{Q}\\
&  +\frac{\vartheta\left(  1-2\gamma\right)  +\lambda+1}{2\Delta(\vartheta
)}\ln\left[  2\frac{\left(  \lambda-1-\vartheta\right)  \sigma+\left(
\lambda+1\right)  ^{2}+\left(  1-\lambda\right)  \vartheta+\Delta
(\vartheta)\Delta(-\sigma)}{\vartheta+\sigma}\right]  .\nonumber
\end{align}

As $z\uparrow1-\gamma$ we have%
\[
Q_{1}\left[  \sigma(z)\right]  \sim-\frac{1}{4}\ln\left(  \lambda\right)
+\frac{\vartheta\left(  1-2\gamma\right)  +\lambda+1}{2\Delta(\vartheta)}%
\ln\left\{  2\left[  \lambda-1-\vartheta+\Delta(\vartheta)\right]  \right\}
\]
and taking (\ref{limit1}) into account we get%
\[
\mathrm{I}_{1}\left[  (1-\gamma)^{-}\right]  =-\frac{1}{4}\ln\left(
\lambda\right)  +\frac{\vartheta\left(  1-2\gamma\right)  +\lambda+1}%
{2\Delta(\vartheta)}\ln\left\{  2\left[  \lambda-1-\vartheta+\Delta
(\vartheta)\right]  \right\}  .
\]
Using (\ref{sigma(0)}) in (\ref{Q}) we have%
\begin{gather*}
Q_{1}\left[  \sigma(0)\right]  =\frac{1}{4}\ln\left(  \frac{\rho^{2}}{\lambda
}\right)  +\frac{\vartheta\left(  1-2\gamma\right)  +\lambda+1}{2\Delta
(\vartheta)}\\
\times\ln\left[  2\frac{\left(  \lambda-1-\vartheta\right)  \rho+\left(
\lambda+1\right)  ^{2}\gamma(1-\gamma)+\left(  1-\lambda\right)
\gamma(1-\gamma)\vartheta+\Delta(\vartheta)\delta}{\gamma(1-\gamma
)\vartheta+\rho}\right]  .
\end{gather*}
From (\ref{sigma(0)}) and (\ref{I1.2}) we get%
\begin{align*}
\mathrm{I}_{1}(0)  &  =-\frac{1}{2}\ln\left[  \frac{\rho}{\gamma
(1-\gamma)\vartheta+\rho}\right]  +Q_{1}\left[  \sigma(0)\right] \\
&  =\frac{1}{2}\ln\left[  \gamma(1-\gamma)\vartheta+\rho\right]  -\frac{1}%
{4}\ln\left(  \lambda\right)  +\frac{\vartheta\left(  1-2\gamma\right)
+\lambda+1}{2\Delta(\vartheta)}\\
&  \times\ln\left[  2\frac{\left(  \lambda-1-\vartheta\right)  \rho+\left(
\lambda+1\right)  ^{2}\gamma(1-\gamma)+\left(  1-\lambda\right)
\gamma(1-\gamma)\vartheta+\Delta(\vartheta)\delta}{\gamma(1-\gamma
)\vartheta+\rho}\right]  .
\end{align*}
Combining these results in (\ref{44}) yields%
\begin{gather*}%
{\displaystyle\int\limits_{0}^{1-\gamma-(1-\alpha)\varepsilon}}
\ln\left[  \frac{\sigma\left(  x\right)  }{\vartheta+\sigma\left(  x\right)
}\right]  dx\sim\frac{\vartheta\left(  1-2\gamma\right)  +\lambda+1}%
{2\Delta(\vartheta)}\ln\left\{  2\left[  \lambda-1-\vartheta+\Delta
(\vartheta)\right]  \right\} \\
-\frac{1}{2}\ln\left[  \gamma(1-\gamma)\vartheta+\rho\right]  -\frac
{\vartheta\left(  1-2\gamma\right)  +\lambda+1}{2\Delta(\vartheta)}\\
\times\ln\left[  2\frac{\left(  \lambda-1-\vartheta\right)  \rho+\left(
\lambda+1\right)  ^{2}\gamma(1-\gamma)+\left(  1-\lambda\right)
\gamma(1-\gamma)\vartheta+\Delta(\vartheta)\delta}{\gamma(1-\gamma
)\vartheta+\rho}\right]
\end{gather*}%
\begin{align*}
&  =-\frac{1}{2}\ln\left[  \gamma(1-\gamma)\vartheta+\rho\right]
+\frac{\vartheta\left(  1-2\gamma\right)  +\lambda+1}{2\Delta(\vartheta)}\\
&  \times\ln\left[  \frac{\left[  \gamma(1-\gamma)\vartheta+\rho\right]
\left[  \lambda-1-\vartheta+\Delta(\vartheta)\right]  }{\left(  \lambda
-1-\vartheta\right)  \rho+\left(  \lambda+1\right)  ^{2}\gamma(1-\gamma
)+\left(  1-\lambda\right)  \gamma(1-\gamma)\vartheta+\Delta(\vartheta)\delta
}\right]  .
\end{align*}

From (\ref{sigma(0)}) we see that as $N\rightarrow\infty$
\[
\frac{1}{2}\left\{  \ln\left[  \frac{\sigma\left(  0\right)  }{\vartheta
+\sigma\left(  0\right)  }\right]  +\ln\left[  \frac{\sigma\left(
1-\gamma+\frac{\alpha-1}{N}\right)  }{\vartheta+\sigma\left(  1-\gamma
+\frac{\alpha-1}{N}\right)  }\right]  \right\}  \sim\frac{1}{2}\ln\left[
\frac{\rho}{\gamma(1-\gamma)\vartheta+\rho}\right]  .
\]
Therefore, we conclude that%
\begin{gather*}
\ln\left[  P(\vartheta)\right]  =\frac{1}{2}\ln\left[  \frac{\rho}%
{\gamma(1-\gamma)\vartheta+\rho}\right]  -\frac{1}{2\varepsilon}\ln\left[
\gamma(1-\gamma)\vartheta+\rho\right]  +\frac{1}{\varepsilon}\frac
{\vartheta\left(  1-2\gamma\right)  +\lambda+1}{2\Delta(\vartheta)}\\
\times\ln\left[  \frac{\left[  \gamma(1-\gamma)\vartheta+\rho\right]  \left[
\lambda-1-\vartheta+\Delta(\vartheta)\right]  }{\left(  \lambda-1-\vartheta
\right)  \rho+\left(  \lambda+1\right)  ^{2}\gamma(1-\gamma)+\left(
1-\lambda\right)  \gamma(1-\gamma)\vartheta+\Delta(\vartheta)\delta}\right]
+o(1).
\end{gather*}
Exponentiating this yields the formula for $\widetilde{P}(\vartheta;N)$.
\end{proof}

We note from (\ref{mu}) that%
\[
\mu\left(  \vartheta\right)  \rightarrow\frac{\gamma}{\delta},\quad
\vartheta\downarrow\theta_{0}%
\]
and that $\widetilde{P}(\vartheta;N)$ has a branch point singularity at
$\vartheta=\theta_{0}.$

We evaluate next $h_{k}(\vartheta)$ in (\ref{Exact}) for $N\rightarrow\infty$
and $k=O(N)$ with $\vartheta$ fixed. We introduce the new variables%
\begin{equation}
y=\frac{x}{N},\quad z=\frac{k}{N},\quad y\geq0,\quad0\leq z\leq1 \label{(y,z)}%
\end{equation}
and consider the asymptotic approximation to $F_{k}(x)$ inside the domain%
\begin{equation}
\mathfrak{D}=\left\{  (y,z):y\geq0,\quad0\leq z\leq1\right\}  . \label{D}%
\end{equation}

Writing $h_{k}(\vartheta)$ in terms of $z$ we have%
\begin{align}
h_{k}(\vartheta)  &  =\frac{1}{2\pi i}%
{\displaystyle\oint\limits_{\mathcal{C}}}
\frac{\left[  1-R_{1}(\vartheta)w\right]  ^{NV(\vartheta)}\left[
1+R_{2}(\vartheta)w\right]  ^{N\left[  1-V(\vartheta)\right]  }}{w^{N(1-z)+1}%
}dw\nonumber\\
&  =\frac{1}{2\pi i}%
{\displaystyle\oint\limits_{\mathcal{C}}}
\frac{1}{w}\exp\left[  N\eta(w,\vartheta,z)\right]  dw \label{eta1}%
\end{align}
with%
\begin{equation}
\eta(w,\vartheta,z)=V(\vartheta)\ln\left[  1-R_{1}(\vartheta)w\right]
+\left[  1-V(\vartheta)\right]  \ln\left[  1+R_{2}(\vartheta)w\right]
-(1-z)\ln(w). \label{eta}%
\end{equation}
We first assume that $\gamma<z<1$ and compute the integral (\ref{eta1}) as
$N\rightarrow\infty$ by the saddle point method \cite{MR89d:41049}. Thus, we
locate the saddle points $W(\vartheta,z)$ by solving the equation
\[
\left.  \frac{\partial}{\partial w}\eta(w,\vartheta,z)\right\vert
_{w=W(\vartheta,z)}=0,
\]
or%
\begin{equation}
zW^{2}+\left[  \left(  \gamma-z\right)  \vartheta+z\lambda-z-\lambda\right]
W+\left(  1-z\right)  \lambda=0 \label{EqW}%
\end{equation}
and find two saddle points $W_{+}\left(  \vartheta,z\right)  $ and
$W_{-}\left(  \vartheta,z\right)  $ defined by
\begin{equation}
W_{\pm}\left(  \vartheta,z\right)  =\frac{\vartheta+1-\lambda}{2}+\frac{1}%
{2z}\left[  \lambda-\gamma\vartheta\pm\sqrt{D(\vartheta,z)}\right]  \label{W}%
\end{equation}
with%
\begin{equation}
D(\vartheta,z)=\rho^{2}+\left[  2\left(  \lambda+1\right)  \rho+2\phi
\vartheta\right]  (z-\gamma)+\left[  \left(  \lambda+1\right)  ^{2}+2\left(
1-\lambda\right)  \vartheta+\vartheta^{2}\right]  (z-\gamma)^{2}. \label{Disc}%
\end{equation}

We next study the motion of the saddle points $W_{\pm}\left(  \vartheta
,z\right)  ,$ and also the branch points $-\frac{1}{R_{2}(\vartheta)}$ and
$\frac{1}{R_{1}(\vartheta)}$ of $\eta(w,\vartheta,z),$ as $\vartheta$
decreases from $\infty$ to $\theta_{0},$ with $z>\gamma.$ First we observe
from (\ref{R12}) that for all real $\vartheta,$\ $R_{1}(\vartheta)\ $and
$R_{2}(\vartheta)>0$ and thus%
\[
-\frac{1}{R_{2}(\vartheta)}<0<\frac{1}{R_{1}(\vartheta)}.
\]
For large $\vartheta$ we have%
\[
\frac{1}{R_{1}}\sim\vartheta,\quad W_{+}\sim\frac{z-\gamma}{z}\vartheta,\quad
W_{-}\sim\frac{1-z}{z-\gamma}\frac{\lambda}{\vartheta}\quad\vartheta
\rightarrow\infty
\]
and hence for all $\vartheta>0$%
\begin{equation}
-\frac{1}{R_{2}}<0<W_{-}<W_{+}<\frac{1}{R_{1}}. \label{range1}%
\end{equation}
When $\vartheta=0$ the branch point $\frac{1}{R_{1}(\vartheta)}$ coalesces
with the saddle point $W_{+}$
\[
\frac{1}{R_{1}(0)}=1,\quad W_{+}(0,z)=1,\quad W_{-}(0,z)=\frac{\lambda}%
{z}\left(  1-z\right)
\]
and we get%
\begin{equation}
-\frac{1}{R_{2}(0)}<0<W_{-}(0,z)<W_{+}(0,z)=\frac{1}{R_{1}(0)}. \label{range2}%
\end{equation}
For $\theta_{0}<\vartheta<0,$ the saddle point $W_{+}$ moves to the right of
the branch point $\frac{1}{R_{1}(\vartheta)}$%
\begin{equation}
-\frac{1}{R_{2}}<0<W_{-}<\frac{1}{R_{1}}<W_{+}. \label{range3}%
\end{equation}
Finally, when $\vartheta=\theta_{0}$ with $\gamma<z<\frac{\gamma^{2}}{\delta}$
the branch point $\frac{1}{R_{1}(\vartheta)}$ coalesces with the saddle point
$W_{-}$%
\[
\frac{1}{R_{1}(\theta_{0})}=\frac{\lambda}{\gamma}\left(  1-\gamma\right)
,\quad W_{+}(\theta_{0},z)=\frac{\left(  1-z\right)  \lambda}{\left(
1-\gamma\right)  z},\quad W_{-}(\theta_{0},z)=\frac{\lambda}{\gamma}\left(
1-\gamma\right)
\]
and we have%
\begin{equation}
-\frac{1}{R_{2}(\theta_{0})}<0<W_{-}(\theta_{0},z)=\frac{1}{R_{1}(\theta_{0}%
)}<W_{+}(\theta_{0},z). \label{range4}%
\end{equation}

Note that when $\vartheta=\vartheta^{\ast},$ where
\begin{equation}
\vartheta^{\ast}(z)=\frac{\lambda(z-1)-z+2\sqrt{z\lambda(1-z)}}{z-\gamma
}<\theta_{0},\quad\gamma<z<\frac{\gamma^{2}}{\delta}, \label{thetastar}%
\end{equation}
the discriminant $D(\vartheta,z)$ vanishes and the saddle points coalesce.
Thus,%
\[
-\frac{1}{R_{2}(\vartheta^{\ast})}<0<\frac{1}{R_{1}(\vartheta^{\ast})}%
<W_{-}(\vartheta^{\ast},z)=W_{+}(\vartheta^{\ast},z)=\sqrt{\frac
{\lambda\left(  1-z\right)  }{z}}.
\]
Furthermore,
\begin{equation}
\frac{d}{dz}D(\vartheta^{\ast},z)=\frac{\partial D}{\partial\vartheta}%
\frac{d\vartheta^{\ast}}{dz}+\frac{\partial D}{\partial z}=0,\quad\text{at
\ }\vartheta=\vartheta^{\ast} \label{derivdiscr}%
\end{equation}
since%
\begin{gather*}
\frac{\partial D}{\partial\vartheta}(\vartheta^{\ast},z)=4\left(
z-\gamma\right)  \sqrt{\lambda z\left(  1-z\right)  },\\
\frac{d\vartheta^{\ast}}{dz}=\frac{\phi\sqrt{\lambda z\left(  1-z\right)
}+\lambda\left(  2z\gamma-z-\gamma\right)  }{\sqrt{\lambda z\left(
1-z\right)  }\left(  z-\gamma\right)  ^{2}},\\
\frac{\partial D}{\partial z}(\vartheta^{\ast},z)=-\frac{4}{z-\gamma}\left[
\phi\sqrt{\lambda z\left(  1-z\right)  }+\lambda\left(  2z\gamma
-z-\gamma\right)  \right]  .
\end{gather*}
Hence, $D(\vartheta,z)$ has a double zero at $\vartheta=\vartheta^{\ast}.$

From (\ref{range1})-(\ref{range4}) we conclude that for $\vartheta>\theta_{0}$
with $z>\gamma$ we must deform the contour $\mathcal{C}$ to the steepest
descent contour through the saddle $W_{-}$ and therefore
\begin{equation}
h_{k}(\vartheta)\sim\frac{1}{\sqrt{2\pi N}}\frac{1}{W_{-}}\exp\left[
N\eta(W_{-},\vartheta,z)\right]  \left[  \left.  \frac{\partial^{2}\eta
}{\partial w^{2}}\right\vert _{w=W_{-}}\right]  ^{-\frac{1}{2}},\quad
\vartheta>\theta_{0},\ \gamma<z<1.\ \label{hsaddle}%
\end{equation}
Note that the steepest descent directions at $w=W_{-}$ are $\arg\left(
w-W_{-}\right)  =\pm\frac{\pi}{2}$ and we can shift $\mathcal{C}$ into another
circle that goes through this saddle point.

Using (\ref{eta}) and (\ref{EqW}) we find that
\begin{equation}
\frac{\partial^{2}\eta}{\partial w^{2}}\left(  W_{\pm},\vartheta,z\right)
=\frac{A(\vartheta,z)W_{\pm}+B(\vartheta,z)}{\left(  W_{\pm}\right)
^{2}\left[  \left(  \vartheta\gamma-\lambda\right)  W_{\pm}+\lambda\right]
^{2}} \label{etaww}%
\end{equation}
with%
\begin{gather*}
A(\vartheta,z)=\lambda^{2}(z-1)\left(  z\lambda+z-\lambda\right)  +\left[
(-2\lambda^{2}+\lambda-\gamma\lambda^{2}-\gamma)z^{2}+\lambda\left(
4\gamma\lambda+2\lambda-\gamma\right)  z-3\gamma\lambda^{2}\right]
\vartheta\\
+\left(  z-\gamma\right)  \left[  \left(  2\gamma\lambda+\lambda
-2\gamma\right)  z-3\gamma\lambda\right]  \vartheta^{2}-\gamma\left(
z-\gamma\right)  ^{2}\vartheta^{3}%
\end{gather*}%
\[
B(\vartheta,z)=\lambda(1-z)\left\{  \left(  z\lambda+z-\lambda\right)
\lambda+\left[  2\gamma\lambda+\left(  \gamma-\lambda-\gamma\lambda\right)
z\right]  \vartheta+\gamma\left(  z-\gamma\right)  \vartheta^{2}\right\}  .
\]
$\allowbreak$

Using (\ref{P2}) and (\ref{hsaddle}) in (\ref{Exact}) we have\
\begin{align}
F_{k}(x)  &  \sim\left(  \frac{\lambda}{\lambda+1}\right)  ^{N}\frac{1}%
{\sqrt{2\pi N}}\frac{1}{2\pi i}%
{\displaystyle\int\limits_{\mathrm{Br}^{+}}}
\frac{1}{\vartheta}\frac{1}{W_{-}}\sqrt{\frac{-\theta_{0}}{\vartheta
-\theta_{0}}}\left[  \left.  \frac{\partial^{2}\eta}{\partial w^{2}%
}\right\vert _{w=W_{-}}\right]  ^{-\frac{1}{2}}\label{F2}\\
&  \times\exp\left\{  N\left[  y\vartheta+\mu(\vartheta)+\eta(W_{-}%
,\vartheta,z)\right]  \right\}  d\vartheta,\quad z>\gamma.\nonumber
\end{align}
To compute the above integral as $N\rightarrow\infty,$ we shall again use the
saddle point method. We obtain $F_{k}(x)\sim G_{1}(y,z)$ with
\begin{align}
G_{1}(y,z)  &  =\frac{1}{2\pi N}\left(  \frac{\lambda}{\lambda+1}\right)
^{N}\frac{1}{\Theta}\frac{1}{W_{-}(\Theta,z)}\sqrt{\frac{-\theta_{0}}%
{\Theta-\theta_{0}}}\sqrt{\frac{1}{\eta_{ww}(W_{-},\Theta,z)}}\label{G1}\\
\times &  \sqrt{\frac{1}{\Psi_{\vartheta\vartheta}(y,W_{-},\Theta,z)}}%
\exp\left[  N\Psi(y,W_{-},\Theta,z)\right]  ,\quad\Theta>0,\quad
z>\gamma\nonumber
\end{align}
where
\begin{equation}
\Psi(y,w,\vartheta,z)=y\vartheta+\mu(\vartheta)+\eta\left(  w,\vartheta
,z\right)  \label{Psi}%
\end{equation}
and $\Theta(y,z)$ is defined to be the solution of the equation%
\begin{equation}
0=\frac{\partial\Psi}{\partial\vartheta}(y,W_{-},\Theta,z)=y+\mu^{\prime
}(\Theta)+\frac{\partial\eta}{\partial\vartheta}\left[  W_{-}\left(
\Theta,z\right)  ,\Theta,z\right]  \label{Theta}%
\end{equation}
where we have used
\[
\frac{\partial\eta}{\partial w}\left[  W_{-}\left(  \Theta,z\right)
,\Theta,z\right]  =0.
\]
For each value of $\Theta,$ (\ref{Theta}) defines implicitly a curve in the
$(y,z)$-plane.

Since (\ref{F2}) has a pole at $\vartheta=0,$ (\ref{G1}) is only valid if $y$
and $z$ are such that $\Theta>0,$ i.e., the saddle lies to the right of the
pole. Then $\mathrm{Br}^{+}$ can be deformed to the steepest descent contour
through $\vartheta=\Theta.$ We also note that when $\Theta=0,$ (\ref{G1}) is
singular. We shall find an asymptotic approximation valid when the saddle is
close to $\vartheta=0$ in Section \ref{transition}. From (\ref{Theta}) we see
that $\Theta=0$ corresponds to the curve%
\begin{equation}
y=Y_{0}(z)=\frac{z-\gamma}{\lambda+1}-\frac{\rho}{\left(  \lambda+1\right)
^{2}}\ln\left(  \frac{z\lambda+z-\lambda}{\rho}\right)  ,\quad\gamma\leq
z\leq1. \label{Y0}%
\end{equation}

For $y>Y_{0}(z)$ the saddle point $\Theta$ lies on the negative real axis, to
the left of the pole at $\vartheta=0.$ We re-write (\ref{Exact}) as%
\begin{equation}
F_{k}(x)-\frac{1}{(1+\lambda)^{N}}\binom{N}{k}\lambda^{k}=\left(
\frac{\lambda}{\lambda+1}\right)  ^{N}\frac{1}{2\pi i}%
{\displaystyle\int\limits_{\mathrm{Br}^{-}}}
e^{x\theta}\frac{1}{\vartheta}\left[
{\displaystyle\prod\limits_{j=0}^{N-\left\lfloor c\right\rfloor -1}}
\frac{\theta_{j}}{\theta_{j}-\vartheta}\right]  h_{k}(\vartheta)d\vartheta
\label{Fnegative}%
\end{equation}
where $\theta_{0}<\operatorname{Re}(\vartheta)<0$ on $\mathrm{Br}^{-}.$ Note
that the residue at \ $\vartheta=0$ was already calculated in (\ref{Res0}) and
corresponds to $F_{k}(\infty)$. Using (\ref{P2}) and (\ref{hsaddle}) in
(\ref{Fnegative}) and approximating the integral by the saddle point method,
we obtain
\begin{gather}
F_{k}(x)-\frac{1}{(1+\lambda)^{N}}\binom{N}{k}\lambda^{k}\sim\frac{1}{2\pi
N}\left(  \frac{\lambda}{\lambda+1}\right)  ^{N}\frac{1}{\Theta}\frac{1}%
{W_{-}(\Theta,z)}\sqrt{\frac{-\theta_{0}}{\Theta-\theta_{0}}}\sqrt{\frac
{1}{\eta_{ww}(W_{-},\Theta,z)}}\label{Fsaddle2}\\
\times\sqrt{\frac{1}{\Psi_{\vartheta\vartheta}(y,W_{-},\Theta,z)}}\exp\left[
N\Psi(y,W_{-},\Theta,z)\right]  ,\quad\theta_{0}<\Theta<0,\quad\gamma
<z<1\nonumber
\end{gather}
with $\Psi(y,w,\vartheta,z)$ and $\Theta(y,z)$ defined by (\ref{Psi}) and
(\ref{Theta}) respectively.

From (\ref{Theta}) we find that the value $\Theta=\theta_{0}$ corresponds to
the curve%
\begin{equation}
y=Y_{1}(z)=\frac{\left[  \gamma(1-\gamma)\right]  ^{2}\rho}{\delta^{2}%
}\left\{  \frac{\gamma^{2}-\delta z}{\gamma(1-\gamma)\rho}+\ln\left[
\frac{\gamma(1-\gamma)\rho}{\gamma^{2}-\delta z}\right]  -1\right\}
,\quad\gamma\leq z<\frac{\gamma^{2}}{\delta}. \label{Y1}%
\end{equation}
Our approximation to the integral in (\ref{Exact}) assumed that $\vartheta
>\theta_{0}.$ However, (\ref{Fsaddle2}) passes smoothly through the curve
$y=Y_{1}(z)$ since from (\ref{etaww}) we have%
\[
\frac{\partial^{2}\eta}{\partial w^{2}}\left(  W_{-},\vartheta,z\right)
\sim\frac{\gamma^{2}\left[  \left(  1-\gamma\right)  ^{2}\lambda\gamma
^{2}+\delta^{2}\left(  \frac{\gamma^{2}}{\delta}-z\right)  \right]  }%
{\lambda^{2}(1-\gamma)^{2}\delta^{2}\left(  \vartheta-\theta_{0}\right)
},\quad\vartheta\rightarrow\theta_{0}.
\]
The discriminant $D(\Theta,z)$ vanishes when%
\begin{equation}
y=Y^{\ast}(z)=-\mu^{\prime}(\vartheta^{\ast})-\frac{\partial\eta}%
{\partial\vartheta}\left[  \sqrt{\frac{\lambda\left(  1-z\right)  }{z}%
},\vartheta^{\ast},z\right]  ,\quad\gamma<z<\frac{\gamma^{2}}{\delta}
\label{Ymax}%
\end{equation}
where $\vartheta^{\ast}(z)$ was defined in (\ref{thetastar}). But, as we
showed in (\ref{derivdiscr}), the derivative of the discriminant also vanishes
along the curve $y=Y^{\ast}(z)$ so that (\ref{Fsaddle2}) can be smoothly
continued through this curve by simply replacing $W_{-}$ by $W_{+}.$

We summarize our results below.

\begin{figure}[t]
\begin{center}
\rotatebox{0} {\resizebox{15cm}{!}{\includegraphics{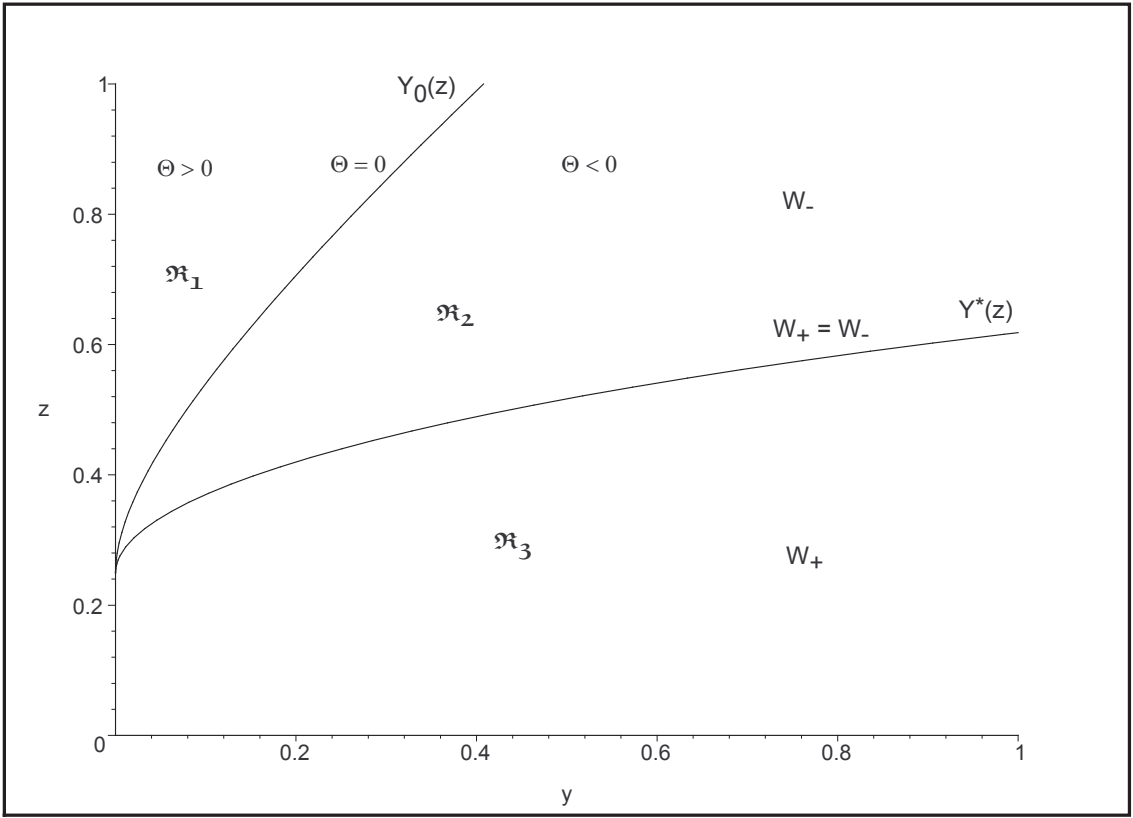}}}
\end{center}
\caption{A sketch of the regions $\mathfrak{R}_{1}, \mathfrak{R}_{2}$ and
$\mathfrak{R}_{3}$.}%
\label{figR1R2R3}%
\end{figure}

\begin{theorem}
\label{theorem1}Defining $y=\frac{x}{N}$ and $z=\frac{k}{N},$ for
$(y,z)\in\overset{\circ}{\mathfrak{D}}=\left\{  (y,z):y>0,\quad0<z<1\right\}
$ we have

\begin{enumerate}
\item
\begin{align}
F_{k}(x)  &  \sim G_{1}(y,z)\equiv\frac{1}{2\pi N}\left(  \frac{\lambda
}{\lambda+1}\right)  ^{N}\frac{1}{\Theta}\frac{1}{W_{-}(\Theta,z)}\sqrt
{\frac{-\theta_{0}}{\Theta-\theta_{0}}}\sqrt{\frac{1}{\eta_{ww}(W_{-}%
,\Theta,z)}}\label{FR1}\\
\times &  \sqrt{\frac{1}{\Psi_{\vartheta\vartheta}(y,W_{-},\Theta,z)}}%
\exp\left[  N\Psi(y,W_{-},\Theta,z)\right]  ,\quad(y,z)\in\mathfrak{R}%
_{1}\nonumber
\end{align}%
\begin{equation}
\mathfrak{R}_{1}=\left\{  (y,z):0<y<Y_{0}(z),\quad\gamma<z<1\right\}  ,
\label{R1}%
\end{equation}%
\[
Y_{0}(z)=\frac{z-\gamma}{\lambda+1}-\frac{\rho}{\left(  \lambda+1\right)
^{2}}\ln\left(  \frac{z\lambda+z-\lambda}{\rho}\right)  ,\quad\gamma\leq
z\leq1
\]
where
\[
\eta(w,\vartheta,z)=V(\vartheta)\ln\left[  1-R_{1}(\vartheta)w\right]
+\left[  1-V(\vartheta)\right]  \ln\left[  1+R_{2}(\vartheta)w\right]
-(1-z)\ln(w)
\]%
\[
W_{\pm}\left(  \vartheta,z\right)  =\frac{\vartheta+1-\lambda}{2}+\frac{1}%
{2z}\left[  \lambda-\gamma\vartheta\pm\sqrt{D(\vartheta,z)}\right]
\]%
\[
D(\vartheta,z)=\rho^{2}+\left[  2\left(  \lambda+1\right)  \rho+2\phi
\vartheta\right]  (z-\gamma)+\left[  \left(  \lambda+1\right)  ^{2}+2\left(
1-\lambda\right)  \vartheta+\vartheta^{2}\right]  (z-\gamma)^{2}.
\]%
\[
\Psi(y,w,\vartheta,z)=y\vartheta+\mu(\vartheta)+\eta\left(  w,\vartheta
,z\right)
\]
and $\Theta(y,z)$ is defined to be the solution of the equation%
\[
y+\mu^{\prime}(\Theta)+\frac{\partial\eta}{\partial\vartheta}\left[
W_{-}\left(  \Theta,z\right)  ,\Theta,z\right]  =0
\]%
\begin{gather*}
\mu(\vartheta)=-\frac{1}{2}\ln\left[  \gamma(1-\gamma)\vartheta+\rho\right] \\
+\frac{\vartheta\left(  1-2\gamma\right)  +\lambda+1}{2\Delta(\vartheta)}%
\ln\left[  \frac{\left[  \gamma(1-\gamma)\vartheta+\rho\right]  \left[
\lambda-1-\vartheta+\Delta(\vartheta)\right]  }{\left(  \lambda-1-\vartheta
\right)  \rho+\left(  \lambda+1\right)  ^{2}\gamma(1-\gamma)+\left(
1-\lambda\right)  \gamma(1-\gamma)\vartheta+\Delta(\vartheta)\delta}\right]  .
\end{gather*}

\item
\begin{equation}
F_{k}(x)-\frac{1}{(1+\lambda)^{N}}\binom{N}{k}\lambda^{k}\sim G_{1}%
(y,z),\quad(y,z)\in\mathfrak{R}_{2} \label{FR2}%
\end{equation}%
\begin{equation}
\mathfrak{R}_{2}=\left\{  (y,z):Y_{0}(z)<y\leq Y^{\ast}(z),\quad
\gamma<z<1\right\}  \label{R2}%
\end{equation}
with%
\[
Y^{\ast}(z)=-\mu^{\prime}(\vartheta^{\ast})-\frac{\partial\eta}{\partial
\vartheta}\left[  \sqrt{\frac{\lambda\left(  1-z\right)  }{z}},\vartheta
^{\ast},z\right]  ,\quad\gamma<z<\frac{\gamma^{2}}{\delta}%
\]%
\[
\vartheta^{\ast}(z)=\frac{\lambda(z-1)-z+2\sqrt{z\lambda(1-z)}}{z-\gamma
},\quad\gamma<z<\frac{\gamma^{2}}{\delta}.
\]

\item
\begin{equation}
F_{k}(x)-\frac{1}{(1+\lambda)^{N}}\binom{N}{k}\lambda^{k}\sim G_{2}%
(y,z),\quad(y,z)\in\mathfrak{R}_{3} \label{FR3}%
\end{equation}%
\begin{equation}
\mathfrak{R}_{3}=\left\{  (y,z):y>Y^{\ast}(z),\quad\gamma<z<\frac{\gamma^{2}%
}{\delta}\right\}  \cup\left\{  (y,z):0<y,\quad0<z\leq\gamma\right\}
\label{R3}%
\end{equation}
where%
\begin{align}
G_{2}(y,z)  &  =\frac{1}{2\pi N}\left(  \frac{\lambda}{\lambda+1}\right)
^{N}\frac{1}{\Theta^{+}}\frac{1}{W_{+}(\Theta^{+},z)}\sqrt{\frac{-\theta_{0}%
}{\Theta^{+}-\theta_{0}}}\sqrt{\frac{1}{\eta_{ww}(W_{+},\Theta^{+},z)}%
}\label{G2}\\
\times &  \sqrt{\frac{1}{\Psi_{\vartheta\vartheta}(y,W_{+},\Theta^{+},z)}}%
\exp\left[  N\Psi(y,W_{+},\Theta^{+},z)\right]  .\nonumber
\end{align}
and $\Theta^{+}(y,z)$ is defined to be the solution of the equation%
\begin{equation}
y+\mu^{\prime}(\Theta^{+})+\frac{\partial\eta}{\partial\vartheta}\left[
W_{+}\left(  \Theta^{+},z\right)  ,\Theta^{+},z\right]  =0. \label{Thetaplus}%
\end{equation}

\end{enumerate}
\end{theorem}

The three regions $\mathfrak{R}_{j}$ are sketched in Figure \ref{figR1R2R3}.
Note that as we pass from $\mathfrak{R}_{1}$ to $\mathfrak{R}_{2}$ there is a
\textquotedblleft phase transition\textquotedblright\ in the asymptotics.
However, the solution is smooth as we pass from $\mathfrak{R}_{2}$ to
$\mathfrak{R}_{3}.$

Theorem \ref{theorem1} ceases to be valid along the curve $y=Y_{0}%
(z),\quad\gamma\leq z\leq1,$ because the function $G_{1}(y,z)$ is singular
there. We shall find a corner layer solution near the point $(0,\gamma)$ in
Section \ref{section corner}, a transition layer solution along $y=Y_{0}%
(z),\quad\gamma<z<1$ in Section \ref{transition} and a corner layer solution
near the point $(Y_{0}(1),1)$ in Section \ref{corner z=1}.

From (\ref{W}) we have%
\begin{align}
W_{-}\left(  \vartheta,z\right)   &  \sim\frac{\lambda}{1+\left(
1-\gamma\right)  \vartheta}\left(  1-z\right)  ,\quad z\rightarrow
1,\quad\vartheta>\theta_{0}\label{Wlimit}\\
W_{+}\left(  \vartheta,z\right)   &  \sim\frac{\lambda-\gamma\vartheta}%
{z},\quad z\rightarrow0,\quad\vartheta<\theta_{0}\nonumber
\end{align}
and from (\ref{eta}) we get
\begin{align}
\eta_{ww}\left(  W_{-},\vartheta,z\right)   &  \sim\frac{\left[  1+\left(
1-\gamma\right)  \vartheta\right]  ^{2}}{\lambda^{2}\left(  1-z\right)
},\quad z\rightarrow1,\quad\vartheta>\theta_{0}\label{etalimit}\\
\eta_{ww}\left(  W_{+},\vartheta,z\right)   &  \sim\frac{z^{3}}{\left(
\lambda-\gamma\vartheta\right)  ^{2}},\quad z\rightarrow0,\quad\vartheta
<\theta_{0}.\nonumber
\end{align}
Hence,
\begin{align*}
G_{1}(y,z)  &  =O\left(  \frac{1}{\sqrt{1-z}}\right)  ,\quad z\rightarrow1\\
G_{2}(y,z)  &  =O\left(  \frac{1}{\sqrt{z}}\right)  ,\quad z\rightarrow0
\end{align*}
and our approximation develops singularities as $z\rightarrow0$ and
$z\rightarrow1,$ which corresponds to $k\approx0$ and $k\approx N.$ We shall
find a boundary layer solution near $z=0$ in Section \ref{sectionz=0} and a
boundary layer solution near $z=1$ in Section \ref{section z=1}. In these
cases we will need to re-examine the expansion of (\ref{eta1}).

Finally, from (\ref{Theta}) we see that $\Theta\rightarrow\infty$ as
$y\rightarrow0\ $if $\gamma<z<1.$ We shall find a boundary layer approximation
near $y=0$ with $z>\gamma$ in Section \ref{section y=0}. To analyze the case
of small $y$ we will need an expansion of $P(\vartheta)$ in Lemma \ref{lemma2}
valid for $\vartheta$ large. Note that Theorem \ref{theorem1} does apply for
$y\rightarrow0$ as long as $0<z<\gamma.$ In this range we can show that
$\Theta^{+}$ in (\ref{Thetaplus}) remains finite. In particular, the boundary
masses $F_{k}(0)$ for $k=Nz$ with $0<z<\gamma$ can be computed by setting
$y=0$ in (\ref{FR3}). We also note that in this range $F_{k}(0)\sim
F_{k}(\infty)$ and the exponentially small (as $N\rightarrow\infty)$
difference can be estimated from (\ref{FR3}).

\section{The corner layer at $(0,\gamma)$ (Region I)}

$\label{section corner}$

We next consider the vicinity of the corner point $(y,z)=(0,\gamma).$ We first
evaluate the product $P(\vartheta)$ for $N\rightarrow\infty$ and
$\vartheta=O(N).$

\begin{lemma}
\label{lemma5}Let $\vartheta=NS,\quad S>0,\quad S=O(1)$ and \
\[
P(NS)=\mathcal{P}(S)=%
{\displaystyle\prod\limits_{j=0}^{N-\left\lfloor c\right\rfloor -1}}
\frac{\theta_{j}}{\theta_{j}-NS}.
\]
Then $\mathcal{P}(S)\sim\widehat{P}(S;N)$ as $N\rightarrow\infty,$ where
\begin{align}
\widehat{P}(S;N)  &  =\frac{1}{\sqrt{2\pi N}}\sqrt{\frac{\rho}{\phi
\gamma\left(  1-\gamma\right)  }}\left(  \frac{\phi}{S}\right)  ^{\alpha}%
\exp\left\{  -N\left(  1-\gamma\right)  \ln\left[  \left(  1-\gamma\right)
SN\right]  -N\gamma\ln(\gamma)\right\} \label{Paprox}\\
&  \times\Gamma\left(  \frac{\phi}{S}+1-\alpha\right)  \exp\left\{  \frac
{\phi}{S}\ln\left[  \frac{\gamma}{\phi\left(  1-\gamma\right)  N}\right]
+\frac{2\phi-\rho-1}{S}\right\}  .\nonumber
\end{align}
Here $\Gamma\left(  \cdot\right)  $ is the Gamma function, $\alpha
=c-\left\lfloor c\right\rfloor $ is the fractional part of $c$, $\phi
=\gamma+\lambda-\gamma\lambda$ and $\rho=\gamma\lambda+\gamma-\lambda.$
\end{lemma}

\begin{proof}
We write%
\[
\mathcal{P}(S)=%
{\displaystyle\prod\limits_{j=0}^{(1-\gamma)N+\alpha-1}}
\frac{\vartheta_{j}}{\vartheta_{j}-SN}=%
{\displaystyle\prod\limits_{j=0}^{(1-\gamma)N+\alpha-1}}
\frac{-\sigma\left(  \frac{j}{N}\right)  }{-\sigma\left(  \frac{j}{N}\right)
-SN}%
\]%
\[
=%
{\displaystyle\prod\limits_{j=0}^{(1-\gamma)N+\alpha-1}}
\frac{\sigma\left(  \frac{j}{N}\right)  \left(  1-\gamma-\frac{j}{N}\right)
}{\left[  \left(  1-\gamma\right)  N-j\right]  S+\sigma\left(  \frac{j}%
{N}\right)  \left(  1-\gamma-\frac{j}{N}\right)  }.
\]
Changing index from $j$ to $n=(1-\gamma)N+\alpha-1-j$ and defining
\begin{equation}
\sigma_{n}=\sigma\left(  1-\gamma+\frac{\alpha-1-n}{N}\right)  \label{sigman}%
\end{equation}
we have%
\[
\mathcal{P}(S)=%
{\displaystyle\prod\limits_{n=0}^{(1-\gamma)N+\alpha-1}}
\frac{\left(  \frac{n+1-\alpha}{N}\right)  \sigma_{n}}{\left(  n+1-\alpha
\right)  S+\left(  \frac{n+1-\alpha}{N}\right)  \sigma_{n}}%
\]%
\begin{align}
&  =%
{\displaystyle\prod\limits_{n=0}^{(1-\gamma)N+\alpha-1}}
\left[  \frac{1}{\phi}\left(  \frac{n+1-\alpha}{N}\right)  \sigma_{n}\right]
\times%
{\displaystyle\prod\limits_{n=0}^{(1-\gamma)N+\alpha-1}}
\frac{\phi}{\phi+(n+1-\alpha)S}\label{3products}\\
&  \times%
{\displaystyle\prod\limits_{n=0}^{(1-\gamma)N+\alpha-1}}
\frac{\phi+(n+1-\alpha)S}{(n+1-\alpha)S+\left(  \frac{n+1-\alpha}{N}\right)
\sigma_{n}}\equiv\mathcal{P}_{1}(S)\mathcal{P}_{2}(S)\mathcal{P}%
_{3}(S).\nonumber
\end{align}
We evaluate individually the three products in (\ref{3products}) beginning
with $\mathcal{P}_{1}(S).$

Setting $\varepsilon=\frac{1}{N}$ and using the Euler-MacLaurin summation
formula, we have%
\begin{align*}
&  \ln\left[
{\displaystyle\prod\limits_{n=0}^{(1-\gamma)N+\alpha-1}}
\varepsilon\left(  n+1-\alpha\right)  \sigma_{n}\right]  =%
{\displaystyle\sum\limits_{n=0}^{(1-\gamma)N+\alpha-1}}
\ln\left[  \varepsilon\left(  n+1-\alpha\right)  \sigma_{n}\right] \\
&  =\frac{1}{2}\ln\left[  \left(  1-\gamma\right)  \sigma_{(1-\gamma
)N+\alpha-1}\right]  +\frac{1}{2}\ln\left[  \varepsilon\left(  1-\alpha
\right)  \sigma_{0}\right] \\
&  +%
{\displaystyle\int\limits_{0}^{(1-\gamma)N+\alpha-1}}
\ln\left\{  \varepsilon\left(  n+1-\alpha\right)  \sigma\left[  1-\gamma
+\left(  \alpha-1-n\right)  \varepsilon\right]  \right\}  dn+o(1).
\end{align*}
Changing variables to $x=1-\gamma+\left(  \alpha-1-n\right)  \varepsilon$ , we
get%
\begin{align*}
&
{\displaystyle\int\limits_{0}^{(1-\gamma)N+\alpha-1}}
\ln\left\{  \varepsilon\left(  n+1-\alpha\right)  \sigma\left[  1-\gamma
+\left(  \alpha-1-n\right)  \varepsilon\right]  \right\}  dn\\
&  =\frac{1}{\varepsilon}%
{\displaystyle\int\limits_{0}^{1-\gamma+\left(  \alpha-1\right)  \varepsilon}}
\ln\left[  \left(  1-\gamma-x\right)  \sigma\left(  x\right)  \right]  dx.
\end{align*}
Defining
\[
\mathrm{I}_{2}(z)=%
{\displaystyle\int\limits^{z}}
\ln\left[  \left(  1-\gamma-x\right)  \sigma\left(  x\right)  \right]  dx,
\]
we have%
\[%
{\displaystyle\int\limits_{0}^{1-\gamma+\left(  \alpha-1\right)  \varepsilon}}
\ln\left[  \left(  1-\gamma-x\right)  \sigma\left(  x\right)  \right]
dx\sim\mathrm{I}_{2}\left(  1-\gamma\right)  -(1-\alpha)\mathrm{I}_{2}%
^{\prime}\left(  1-\gamma\right)  \ \varepsilon-\mathrm{I}_{2}(0),\quad
\varepsilon\rightarrow0.
\]
From (\ref{sigma}) we have%
\[
\left(  1-\gamma-x\right)  \sigma\left(  x\right)  =\frac{\rho+2(\lambda
-1)(1-x)x+(1-2x)\sqrt{\rho^{2}+4\lambda(1-x)x}}{2(\gamma-x)}%
\]
and evaluating this at \ $x=1-\gamma$ we get%
\[
\mathrm{I}_{2}^{\prime}\left(  1-\gamma\right)  =\ln(\phi).
\]
After some calculations, we find that%
\begin{align}
\mathrm{I}_{2}(z)  &  =(\gamma-z)\ln(2)-z-\ln(1-z)\nonumber\\
&  +\frac{\gamma}{2}\ln\left[  \rho^{2}-2\lambda\left(  2z\gamma
-\gamma-z\right)  +\phi K_{1}(z)\right] \nonumber\\
&  +\frac{1}{2}(\gamma-1)\ln\left[  \rho^{2}+2\lambda\left(  2z\gamma
-\gamma-z+1\right)  +\phi K_{1}(z)\right] \label{I2}\\
&  +z\ln\left[  \frac{\rho+2z\left(  1-z\right)  \left(  \lambda-1\right)
+(1-2z)K_{1}(z)}{\gamma-z}\right] \nonumber\\
&  +\frac{1}{2}\ln\left[  \rho^{2}+2\lambda\left(  1-z\right)  +\rho
K_{1}(z)\right] \nonumber
\end{align}
where%
\[
K_{1}(z)=\sqrt{\rho^{2}+4\lambda z\left(  1-z\right)  }.
\]
From (\ref{I2}) we obtain%

\[
\mathrm{I}_{2}(0)=2\gamma\ln(2)+\frac{1}{2}(\gamma+1)\ln(\lambda
+1)+\frac{\gamma}{2}\ln\left(  \delta\right)  +\gamma\ln(\gamma)
\]
and%

\[
\mathrm{I}_{2}(1-\gamma)=2\gamma\ln(2)+\frac{1}{2}(\gamma+1)\ln(\lambda
+1)+\frac{\gamma}{2}\ln\left(  \delta\right)  +\gamma-1.
\]
Thus,%
\[
\mathrm{I}_{2}(1-\gamma)-\mathrm{I}_{2}(0)=\gamma-1-\gamma\ln(\gamma)
\]
and we have%
\[
\frac{1}{\varepsilon}%
{\displaystyle\int\limits_{0}^{1-\gamma+\left(  \alpha-1\right)  \varepsilon}}
\ln\left[  \left(  1-\gamma-x\right)  \sigma\left(  x\right)  \right]
dx\sim\left[  \gamma-1-\gamma\ln(\gamma)\right]  \frac{1}{\varepsilon}+\left(
\alpha-1\right)  \ln(\phi),\quad\varepsilon\rightarrow0.
\]
From (\ref{sigma}) and (\ref{sigman}) we get%
\begin{align*}
&  \frac{1}{2}\left\{  \ln\left[  \left(  1-\gamma\right)  \sigma
_{(1-\gamma)N+\alpha-1}\right]  +\ln\left[  \left(  1-\alpha\right)
\varepsilon\sigma_{0}\right]  \right\} \\
&  =\frac{1}{2}\ln\left[  \left(  1-\gamma\right)  \sigma\left(  0\right)
\right]  +\frac{1}{2}\ln\left\{  \left(  1-\alpha\right)  \varepsilon
\sigma\left[  1-\gamma+\left(  \alpha-1\right)  \varepsilon\right]  \right\}
\sim\frac{1}{2}\ln\left(  \frac{\rho\phi}{\gamma}\right)  ,\quad
\varepsilon\rightarrow0.
\end{align*}
Therefore we conclude that%
\begin{align*}
&  \ln\left[
{\displaystyle\prod\limits_{n=0}^{(1-\gamma)N+\alpha-1}}
\left(  \frac{n+1-\alpha}{N}\right)  \sigma_{n}\right] \\
&  \sim\left[  \gamma-1-\gamma\ln(\gamma)\right]  N+\left(  \alpha-1\right)
\ln(\phi)+\frac{1}{2}\ln\left(  \frac{\rho\phi}{\gamma}\right)  ,\quad
N\rightarrow\infty
\end{align*}
and then%
\begin{align}
&  \mathcal{P}_{1}(S)=%
{\displaystyle\prod\limits_{n=0}^{(1-\gamma)N+\alpha-1}}
\frac{1}{\phi}\left(  \frac{n+1-\alpha}{N}\right)  \sigma_{n}=\phi
^{-(1-\gamma)N-\alpha}%
{\displaystyle\prod\limits_{n=0}^{(1-\gamma)N+\alpha-1}}
\left(  \frac{n+1-\alpha}{N}\right)  \sigma_{n}\label{P1 final}\\
&  \sim\sqrt{\frac{\rho}{\gamma\phi}}\exp\left[  -(1-\gamma)\ln(\phi
)N+(\gamma-1)N-\gamma\ln(\gamma)N\right]  ,\quad N\rightarrow\infty.\nonumber
\end{align}

We next consider the second product in (\ref{3products}). Using properties of
Gamma functions \cite{MR94b:00012} we have%
\[
\mathcal{P}_{2}(S)=%
{\displaystyle\prod\limits_{n=0}^{(1-\gamma)N+\alpha-1}}
\frac{\phi}{\phi+S(n+1-\alpha)}=\left(  \frac{\phi}{S}\right)  ^{(1-\gamma
)N+\alpha}\frac{\Gamma\left(  \frac{\phi}{S}+1-\alpha\right)  }{\Gamma\left[
(1-\gamma)N+\frac{\phi}{S}+1\right]  }%
\]
and Stirling's formula \cite{MR97k:01072} gives
\begin{align}
\mathcal{P}_{2}(S)  &  \sim\frac{1}{\sqrt{2\pi N}}\left(  \frac{\phi}%
{S}\right)  ^{\alpha}\Gamma\left(  \frac{\phi}{S}+1-\alpha\right)  \left[
\left(  1-\gamma\right)  N\right]  ^{-\frac{\phi}{S}}\frac{1}{\sqrt{1-\gamma}%
}\label{P2final}\\
&  \times\exp\left\{  (1-\gamma)\ln\left(  \frac{\phi}{S}\right)
N+(1-\gamma)N-(1-\gamma)\left[  \ln\left(  1-\gamma\right)  N\right]
\right\}  .\nonumber
\end{align}

Finally we determine the asymptotic approximation to $\mathcal{P}_{3}(S).$ The
Euler-MacLaurin summation formula yields%
\[
\ln\left[  \mathcal{P}_{3}(S)\right]  =%
{\displaystyle\sum\limits_{n=0}^{(1-\gamma)N+\alpha-1}}
\ln\left[  \frac{\phi+(n+1-\alpha)S}{(n+1-\alpha)\left(  S+\varepsilon
\sigma_{n}\right)  }\right]
\]%
\begin{gather}
=\frac{1}{2}\ln\left[  \frac{\phi+(1-\gamma)SN}{(1-\gamma)\left(
SN+\sigma_{(1-\gamma)N+\alpha-1}\right)  }\right]  +\frac{1}{2}\ln\left[
\frac{\phi+(1-\alpha)S}{(1-\alpha)\left(  S+\varepsilon\sigma_{0}\right)
}\right] \label{int1}\\
+%
{\displaystyle\int\limits_{0}^{(1-\gamma)N+\alpha-1}}
\ln\left\{  \frac{\phi+(n+1-\alpha)S}{(n+1-\alpha)\left(  S+\varepsilon
\sigma\left[  1-\gamma+\left(  \alpha-1-n\right)  \varepsilon\right]  \right)
}\right\}  dn+o(1).\nonumber
\end{gather}
Changing variables in the integral in (\ref{int1}) from $n$ to $x=1-\gamma
+\left(  \alpha-1-n\right)  \varepsilon,$ we have%
\[%
{\displaystyle\int\limits_{0}^{(1-\gamma)N+\alpha-1}}
\ln\left\{  \frac{\phi+(n+1-\alpha)S}{(n+1-\alpha)\left(  S+\varepsilon
\sigma\left[  1-\gamma+\left(  \alpha-1-n\right)  \varepsilon\right]  \right)
}\right\}  dn
\]%
\begin{equation}
=\frac{1}{\varepsilon}%
{\displaystyle\int\limits_{0}^{1-\gamma+\left(  \alpha-1\right)  \varepsilon}}
\ln\left\{  \frac{\phi\varepsilon+S(1-\gamma-x)}{(1-\gamma-x)\left[
\varepsilon\sigma(x)+S\right]  }\right\}  dx. \label{int.1}%
\end{equation}
We write (\ref{int.1}) as the difference of two integrals%
\[%
{\displaystyle\int\limits_{0}^{1-\gamma+\left(  \alpha-1\right)  \varepsilon}}
\ln\left\{  \frac{\phi\varepsilon+(1-\gamma-x)S}{(1-\gamma-x)\left[
\varepsilon\sigma(x)+S\right]  }\right\}  dx
\]%
\begin{equation}
=%
{\displaystyle\int\limits_{0}^{1-\gamma+\left(  \alpha-1\right)  \varepsilon}}
\ln\left[  \frac{\phi\varepsilon}{1-\gamma-x}+S\right]  dx-%
{\displaystyle\int\limits_{0}^{1-\gamma+\left(  \alpha-1\right)  \varepsilon}}
\ln\left[  \varepsilon\sigma(x)+S\right]  dx. \label{int.2}%
\end{equation}
The first integral in (\ref{int.2}) is easily evaluated as $\varepsilon
\rightarrow0$ to give%
\begin{gather*}%
{\displaystyle\int\limits_{0}^{1-\gamma+\left(  \alpha-1\right)  \varepsilon}}
\ln\left[  \frac{\phi\varepsilon}{1-\gamma-x}+S\right]  dx=\left(
1-\gamma\right)  \ln\left(  S\right) \\
+\frac{\phi}{S}\ln\left\{  \frac{\left(  1-\gamma\right)  eS}{\left[
\phi+\left(  1-\alpha\right)  S\right]  \varepsilon}\right\}  \varepsilon
+\left(  1-\alpha\right)  \ln\left[  \frac{1-\alpha}{\phi+\left(
1-\alpha\right)  S}\right]  \varepsilon+o(\varepsilon).
\end{gather*}
$\allowbreak\allowbreak$ $\allowbreak\allowbreak$Introducing the function%
\[
\mathrm{I}_{3}(z)=%
{\displaystyle\int\limits^{z}}
\ln\left[  \varepsilon\sigma(x)+S\right]  dx,
\]
in (\ref{int.2}) and using (\ref{V}) we have%
\begin{gather*}
\mathrm{I}_{3}(z)=z\ln\left[  \varepsilon\sigma(z)+S\right]  -%
{\displaystyle\int\limits^{z}}
\frac{\varepsilon x}{\varepsilon\sigma(x)+S}\sigma^{\prime}(x)dx\\
=z\ln\left[  \varepsilon\sigma(z)+S\right]  -%
{\displaystyle\int\limits^{\sigma(z)}}
\frac{\varepsilon V\left(  -\sigma\right)  }{\varepsilon\sigma+S}d\sigma.
\end{gather*}
After some work we find that%
\begin{gather}
\mathrm{I}_{3}(z)=\left(  z-\frac{1}{2}\right)  \ln\left[  \varepsilon
\sigma(z)+S\right]  +\frac{\left(  2\gamma-1\right)  S-\left(  \lambda
+1\right)  \varepsilon}{2Q_{2}}\nonumber\\
\times\ln\left\{  2\frac{\Delta\left[  -\sigma(z)\right]  Q_{2}+\left[
\left(  \lambda-1\right)  \varepsilon-S\right]  \sigma(z)+\left(
1-\lambda\right)  S+\left(  \lambda+1\right)  ^{2}\varepsilon}{\varepsilon
\sigma(z)+S}\right\} \label{int4.1}\\
+\left(  \gamma-\frac{1}{2}\right)  \ln\left\{  \lambda-1+\sigma
(z)+\Delta\left[  -\sigma(z)\right]  \right\} \nonumber
\end{gather}
with%
\[
Q_{2}(S)=\sqrt{\left[  \left(  \lambda+1\right)  \varepsilon\right]
^{2}-2\left(  \lambda-1\right)  \varepsilon S+S^{2}}.
\]
From (\ref{int4.1}) we have, as $\varepsilon\rightarrow0$
\begin{align*}
\mathrm{I}_{3}\left[  1-\gamma+\left(  \alpha-1\right)  \varepsilon\right]
&  \sim\left(  \gamma-\frac{1}{2}\right)  \ln\left(  \frac{8\lambda}%
{S}\right)  -\frac{\phi}{S}\ln\left\{  \frac{4\varepsilon\lambda\left[
\phi+\left(  1-\alpha\right)  S\right]  }{\phi S}\right\}  \varepsilon\\
&  +\left(  1-\alpha\right)  \ln\left[  \frac{1-\alpha}{\phi+\left(
1-\alpha\right)  S}\right]  \varepsilon+\frac{\left(  2\gamma-1\right)
\left(  \lambda-1\right)  }{S}\varepsilon
\end{align*}
and%
\[
\mathrm{I}_{3}(0)\sim\left(  \frac{1}{2}-\gamma\right)  \ln\left(  \frac
{S}{8\lambda}\right)  -\frac{\phi}{S}\ln\left[  \frac{4\lambda\left(
1-\gamma\right)  }{\gamma}\right]  \varepsilon+\left[  \frac{1}{2}\left(
\lambda-1\right)  -\rho\right]  \frac{\varepsilon}{S}.
\]
We conclude that%
\[
\frac{1}{\varepsilon}%
{\displaystyle\int\limits_{0}^{1-\gamma+\left(  \alpha-1\right)  \varepsilon}}
\ln\left\{  \frac{\phi\varepsilon+S(1-\gamma-x)}{(1-\gamma-x)\left[
\varepsilon\sigma(x)+S\right]  }\right\}  dx\sim\frac{\phi}{S}\ln\left(
\frac{\gamma}{\phi}\right)  +\frac{2\phi-\rho-1}{S}.
\]
Finally we compute the remaining part of (\ref{int1})%
\begin{align*}
&  \frac{1}{2}\ln\left[  \frac{\phi\varepsilon+(1-\gamma)S}{(1-\gamma
)S+(1-\gamma)\varepsilon\sigma_{(1-\gamma)N+\alpha-1}}\right]  +\frac{1}{2}%
\ln\left[  \frac{\phi+(1-\alpha)S}{(1-\alpha)S+\left(  1-\alpha\right)
\varepsilon\sigma_{0}}\right] \\
&  =\frac{1}{2}\ln\left[  \frac{\phi\varepsilon+(1-\gamma)S}{(1-\gamma
)S+(1-\gamma)\sigma(0)\varepsilon}\right]  +\frac{1}{2}\ln\left[  \frac
{\phi+S(1-\alpha)}{(1-\alpha)S+\left(  1-\alpha\right)  \varepsilon
\sigma\left[  1-\gamma+\left(  \alpha-1\right)  \varepsilon\right]  }\right]
\\
&  \sim\frac{1}{2}\frac{\phi\gamma-\rho}{\left(  1-\gamma\right)  \gamma
S}\varepsilon-\frac{1}{2}\frac{\left(  1-\alpha\right)  \left(  \gamma
\lambda^{2}+\gamma-4\gamma\lambda-\lambda^{2}+2\lambda\right)  }{\phi\left[
\phi+\left(  1-\alpha\right)  S\right]  }\varepsilon=o(1).
\end{align*}
Hence,%
\begin{equation}
\ln\left[  \mathcal{P}_{3}(S)\right]  \sim\frac{\phi}{S}\ln\left(
\frac{\gamma}{\phi}\right)  +\frac{2\phi-\rho-1}{S},\quad N\rightarrow\infty.
\label{P3 final}%
\end{equation}

Combining (\ref{P1 final}), (\ref{P2final}) and (\ref{P3 final}) we obtain%
\[
\mathcal{P}_{1}(S)\mathcal{P}_{2}(S)\mathcal{P}_{3}(S)\sim\widehat{P}(S;N)
\]
as in (\ref{Paprox}).
\end{proof}

Next, we shall analyze the function $h_{k}(\vartheta)$ for $N\rightarrow
\infty$ and $\left(  y,z\right)  $ close to the corner $(0,\gamma).$ To do so,
we introduce the new variable $l,$ defined by%
\begin{equation}
l=k-c+\alpha,\quad-\infty<l<\infty,\quad l=O(1) \label{l}%
\end{equation}
and the new function $h^{(1)}(l,S)=h_{l+c-\alpha}(SN).$ From (\ref{h}) and
(\ref{eta}) we have%
\begin{equation}
h^{(1)}(l,S)=\frac{1}{2\pi i}%
{\displaystyle\oint\limits_{\mathcal{C}}}
\exp\left[  \eta\left(  w,SN,\gamma\right)  \right]  \frac{1}{w^{\alpha-l+1}%
}dw \label{hl.1}%
\end{equation}
Using (\ref{VR1}) we see that as $N\rightarrow\infty$ the branch points are
located at
\[
\frac{1}{R_{1}(SN)}\sim SN,\quad-\frac{1}{R_{2}(SN)}\sim-\frac{\lambda}{SN}%
\]
and that, for $N\rightarrow\infty$ with $S$ and $w$ fixed,%
\begin{align}
\exp\left[  \eta\left(  w,SN,\gamma\right)  \right]   &  \sim\exp\left\{
\left[  \left(  1-\gamma\right)  N+\frac{\phi}{S}\right]  \ln\left(  \frac
{SN}{\lambda}\right)  +\frac{\left(  1-\lambda\right)  \left(  1-\gamma
\right)  }{S}\right\} \label{etal.1}\\
&  \times\exp\left[  \frac{\lambda\left(  1-\gamma\right)  }{Sw}-\frac{\gamma
}{S}w\right]  w^{\frac{\phi}{S}}.\nonumber
\end{align}
Therefore, we replace the contour $\mathcal{C}$ in (\ref{hl.1}) by a Hankel's
loop and use (\ref{etal.1}) to obtain
\begin{align*}
h^{(1)}(l,S)  &  \sim\exp\left[  \left(  \frac{1-\gamma}{\varepsilon}%
+\frac{\phi}{S}\right)  \ln\left(  \frac{S}{\lambda\varepsilon}\right)
+\frac{\left(  1-\lambda\right)  \left(  1-\gamma\right)  }{S}\right] \\
&  \times\frac{1}{2\pi i}%
{\displaystyle\oint\limits_{\mathcal{H}}}
\exp\left[  \frac{\lambda\left(  1-\gamma\right)  }{Sw}-\frac{\gamma}%
{S}w\right]  w^{\frac{\phi}{S}+l-\alpha-1}dw.
\end{align*}
where $\mathcal{H}$ begins at $w=-\infty,$ encircles $w=0$ once in a positive
sense and returns to its starting point. We assume that the branch of
$w^{\frac{\phi}{S}+l-\alpha-1}$ takes its principal value at the point where
the contour crosses the positive real axis and is continuous everywhere,
except for $w$ real and negative \cite{MR97i:41001}.

Introducing the new variable $t=\frac{\beta}{2\gamma}\frac{1}{w},$ with%
\begin{equation}
\beta=2\sqrt{\lambda\gamma\left(  1-\gamma\right)  } \label{beta}%
\end{equation}
we get%
\begin{align*}
h^{(1)}(l,S)  &  \sim\exp\left[  \left(  \frac{1-\gamma}{\varepsilon}%
+\frac{\phi}{S}\right)  \ln\left(  \frac{S}{\lambda\varepsilon}\right)
+\frac{\left(  1-\lambda\right)  \left(  1-\gamma\right)  }{S}+\left(
\frac{\phi}{S}+l-\alpha\right)  \ln\left(  \frac{\beta}{2\gamma}\right)
\right] \\
&  \times\frac{1}{2\pi i}%
{\displaystyle\int\limits_{\mathcal{H}^{\prime}}}
\exp\left[  \frac{\beta}{2S}\left(  t-\frac{1}{t}\right)  \right]
t^{-\frac{\phi}{S}-l+\alpha-1}dt
\end{align*}
where $\mathcal{H}^{\prime}$ is another Hankel contour of the same shape as
$\mathcal{H}.$ Using an integral representation of the Bessel function
$J_{\cdot}(\cdot)$ \cite{MR96i:33010} we obtain%
\begin{equation}
h_{k}\left(  \vartheta\right)  =h^{(1)}(l,S)\sim\left(  \frac{SN}{\lambda
}\right)  ^{(1-\gamma)N+\frac{\phi}{S}}e^{(1-\phi)\frac{1}{S}}\left(
\frac{\beta}{2\gamma}\right)  ^{\frac{\phi}{S}+l-\alpha}J_{l-\alpha+\frac
{\phi}{S}}\left(  \frac{\beta}{S}\right)  . \label{h1}%
\end{equation}

Using (\ref{Paprox}) and (\ref{h1}) we obtain the following theorem.

\begin{theorem}
\label{corner} Let $k=c-\alpha+l=\left\lfloor c\right\rfloor +l,\quad
-\infty<l<\infty,\quad l=O(1)$ and $x=\frac{\chi}{N},\quad\chi=O(1).$ Then
$F_{k}(x)\sim F_{l}^{(1)}(\chi)$ as $N\rightarrow\infty$ where%
\begin{align}
F_{l}^{(1)}(\chi)  &  =\frac{1}{\sqrt{2\pi N}}\sqrt{\frac{\rho}{\gamma
\phi\left(  1-\gamma\right)  }}\left(  \frac{\beta}{2\gamma}\right)
^{l-\alpha}\exp\left[  N\Phi(\gamma)\right] \nonumber\\
&  \times\frac{1}{2\pi i}%
{\displaystyle\int\limits_{\mathrm{Br}^{+}}}
e^{\chi S}\frac{1}{S}\Gamma\left(  \frac{\phi}{S}+1-\alpha\right)
J_{l-\alpha+\frac{\phi}{S}}\left(  \frac{\beta}{S}\right) \label{Fcorner}\\
&  \times\exp\left\{  \frac{\phi}{S}\ln\left[  \frac{\beta S}{2\lambda
\phi\left(  1-\gamma\right)  }\right]  +\frac{\phi-\rho}{S}\right\}  \left(
\frac{\phi}{S}\right)  ^{\alpha}dS,\nonumber
\end{align}
$\mathrm{Br}^{+}$ is a vertical contour in the complex plane with
$\operatorname{Re}(S)>0,$ $J_{\cdot}(\cdot)$ is the Bessel function,
$\Gamma\left(  \cdot\right)  $ is the Gamma function ,%
\[
\beta=2\sqrt{\lambda\gamma\left(  1-\gamma\right)  }%
\]
and%
\begin{equation}
\Phi(z)=-z\ln(z)-(1-z)\ln\left(  1-z\right)  +z\ln\left(  \lambda\right)
-\ln\left(  \lambda+1\right)  . \label{Phi}%
\end{equation}

\end{theorem}

\section{The transition layer at$\ y=Y_{0}(z)$ (Region II)}

\label{transition}

As we observed previously in Section \ref{main approx}, the approximation
(\ref{G1}) is not valid along the curve $y=Y_{0}(z)$ given by%
\[
Y_{0}(z)=\frac{z-\gamma}{\lambda+1}-\frac{\rho}{\left(  \lambda+1\right)
^{2}}\ln\left(  \frac{z\lambda+z-\lambda}{\rho}\right)  ,\quad\gamma\leq
z\leq1.
\]
The condition $y\approx Y_{0}(z)$ corresponds to the saddle $\Theta$ being
close to the pole at $\vartheta=0$ in (\ref{Exact}). To obtain the expansion
of (\ref{G1}) for $y\approx Y_{0}(z),$ we first approximate the integrand in
(\ref{F2}) for small $\vartheta,$ using%
\[
W_{-}\sim\frac{\lambda\left(  1-z\right)  }{z},\quad\eta_{ww}(W_{-}%
,\vartheta,z)\sim\frac{z^{3}}{\lambda^{2}\left(  1-z\right)  }%
\]%
\[
\eta(W_{-},\vartheta,z)+\mu\left(  \vartheta\right)  \sim\ln\left[  \left(
1-z\right)  ^{z-1}\lambda^{z-1}z^{-z}\right]  -Y_{0}(z)\vartheta+\frac
{Y_{2}(z)}{2}\vartheta^{2}%
\]%
\begin{equation}
Y_{2}(z)=\frac{2\zeta}{\left(  \lambda+1\right)  ^{4}}\ln\left(
\frac{z\lambda+z-\lambda}{\rho}\right)  -\frac{\left(  z-\gamma\right)
\left[  2\rho\zeta+3\left(  \lambda+1\right)  \zeta\left(  z-\gamma\right)
+\allowbreak\left(  \lambda-1\right)  \left(  \lambda+1\right)  ^{2}\left(
z-\gamma\right)  ^{2}\right]  }{\left(  z\lambda+z-\lambda\right)  ^{2}\left(
\lambda+1\right)  ^{3}} \label{Y2}%
\end{equation}%
\begin{equation}
\zeta=\lambda^{2}\left(  \gamma-1\right)  +2\lambda-\gamma. \label{zeta}%
\end{equation}
Using the above in (\ref{F2}) yields%
\[
F_{k}(x)\sim\frac{1}{\sqrt{2\pi N}}\sqrt{\frac{1}{z\left(  1-z\right)  }}%
\exp\left[  \Phi(z)N\right]  \frac{1}{2\pi i}%
{\displaystyle\int\limits_{\mathrm{Br}^{+}}}
\frac{1}{\vartheta}\exp\left\{  \left[  y\vartheta-Y_{0}(z)\vartheta
+\frac{Y_{2}(z)}{2}\vartheta^{2}\right]  N\right\}  d\vartheta
\]
with $\Phi(z)$ given by (\ref{Phi}). Changing variables to
\begin{equation}
\vartheta=\frac{v}{\sqrt{Y_{2}(z)N}},\quad V(y,z)=\left[  y-Y_{0}(z)\right]
\sqrt{\frac{N}{Y_{2}(z)}} \label{V.1}%
\end{equation}
we obtain \cite{MR2001c:00002}%
\begin{align*}
&  \frac{1}{2\pi i}%
{\displaystyle\int\limits_{\mathrm{Br}^{+}}}
{}\exp\left\{  \left[  y-Y_{0}(z)\right]  \vartheta N+\frac{Y_{2}(z)}%
{2}\vartheta^{2}N\right\}  \frac{1}{\vartheta}d\vartheta\\
&  =\frac{1}{2\pi i}%
{\displaystyle\int\limits_{\mathrm{Br}^{\prime}}}
{}\exp\left(  vV+\frac{v^{2}}{2}\right)  \frac{1}{v}dv=\frac{1}{\sqrt{2\pi}}%
{\displaystyle\int\limits_{-\infty}^{V}}
\exp\left(  -\frac{t^{2}}{2}\right)  dt
\end{align*}
where $\mathrm{Br}^{\prime}$ is another vertical contour, on which
$\operatorname{Re}(v)>0.$ Therefore,%
\begin{equation}
F_{k}(x)=F^{(2)}(V,z)\sim\frac{1}{\sqrt{2\pi N}}\sqrt{\frac{1}{z\left(
1-z\right)  }}\exp\left[  \Phi(z)N\right]  \frac{1}{\sqrt{2\pi}}%
{\displaystyle\int\limits_{-\infty}^{V}}
\exp\left(  -\frac{t^{2}}{2}\right)  dt. \label{F2.1}%
\end{equation}
This is valid for $y-Y_{0}(z)=O\left(  N^{-\frac{1}{2}}\right)  $ with
$\gamma<z<1$ and gives the transition between regions $\mathfrak{R}_{1}$ and
$\mathfrak{R}_{2}$ in Theorem \ref{theorem1}.

\section{The boundary $z=0$ (Region III)}

$\label{sectionz=0}$

We shall next analyze the part of $\partial\mathfrak{D}$ corresponding to
$z=0.$ We start by proving the following lemma.

\begin{lemma}
Let $\ k=O(1).$ Then as $N\rightarrow\infty$%
\begin{equation}
h_{k}(\vartheta)\sim\frac{N^{k}}{k!}\left(  \lambda-\gamma\vartheta\right)
^{k}\left[  -R_{1}(\vartheta)\right]  ^{V(\vartheta)N}\left[  R_{2}%
(\vartheta)\right]  ^{\left[  1-V(\vartheta)\right]  N},\quad\Theta^{+}\left(
0,0\right)  <\vartheta<\theta_{0} \label{hz=0}%
\end{equation}
where $\Theta^{+}\left(  y,z\right)  $ is defined in (\ref{Thetaplus}).
\end{lemma}

\begin{proof}
From (\ref{H}) we have for $N\rightarrow\infty,$ with $w=O(N),$
\begin{equation}
H\left(  \vartheta,w\right)  =w^{N}\left[  -R_{1}(\vartheta)\right]
^{V(\vartheta)N}\left[  R_{2}(\vartheta)\right]  ^{\left[  1-V(\vartheta
)\right]  N}\exp\left[  \frac{\left(  \lambda-\gamma\vartheta\right)  }%
{w}N+O\left(  \frac{N}{w^{2}}\right)  \right]  . \label{h0.1}%
\end{equation}
We replace $\mathcal{C}$ by $\mathcal{C}^{\prime}=\frac{\mathcal{C}}{N}$ and
change variables from $w$ to $u=\frac{w}{N}$ in (\ref{h}) to get
\[
h_{k}(\vartheta)=N^{k-N}\frac{1}{2\pi i}%
{\displaystyle\oint\limits_{\mathcal{C}}}
\frac{H(\vartheta,Nu)}{u^{N-k+1}}du.
\]
Using (\ref{h0.1}) we have%
\[
h_{k}(\vartheta)\sim N^{k}\left[  -R_{1}(\vartheta)\right]  ^{V(\vartheta
)N}\left[  R_{2}(\vartheta)\right]  ^{\left[  1-V(\vartheta)\right]  N}%
\frac{1}{2\pi i}%
{\displaystyle\oint\limits_{\mathcal{C}}}
\exp\left[  \frac{\left(  \lambda-\gamma\vartheta\right)  }{u}\right]
\frac{1}{u^{1-k}}du
\]
which evaluates to give (\ref{hz=0}).
\end{proof}

Writing (\ref{Exact}) as
\[
F_{k}(x)-F_{k}(\infty)=\left(  \frac{\lambda}{\lambda+1}\right)  ^{N}\frac
{1}{2\pi i}%
{\displaystyle\int\limits_{\mathrm{Br}^{-}}}
e^{x\vartheta}\frac{1}{\vartheta}\left[
{\displaystyle\prod\limits_{j=0}^{N-\left\lfloor c\right\rfloor -1}}
\frac{\theta_{j}}{\theta_{j}-\vartheta}\right]  h_{k}(\vartheta)d\vartheta
\]
and using (\ref{P2}) and (\ref{hz=0}) we have%
\begin{gather*}
F_{k}(x)-F_{k}(\infty)\sim\frac{N^{k}}{k!}\left(  \frac{\lambda}{\lambda
+1}\right)  ^{N}\frac{1}{2\pi i}%
{\displaystyle\int\limits_{\mathrm{Br}^{-}}}
\frac{1}{\vartheta}\sqrt{\frac{-\theta_{0}}{\vartheta-\theta_{0}}}\left(
\lambda-\gamma\vartheta\right)  ^{k}\\
\times\exp\left\{  y\vartheta N+\mu(\vartheta)N+V(\vartheta)\ln\left[
-R_{1}(\vartheta)\right]  N+\left[  1-V(\vartheta)\right]  \ln\left[
R_{2}(\vartheta)\right]  N\right\}  d\vartheta.
\end{gather*}
Use of the saddle point method yields%
\begin{gather}
F_{k}(x)-F_{k}(\infty)\equiv F_{k}^{(3)}(y)-F_{k}(\infty)\sim\frac{N^{k}}%
{k!}\left(  \frac{\lambda}{\lambda+1}\right)  ^{N}\frac{1}{\sqrt{2\pi N}%
}\label{F3.5}\\
\times\frac{1}{\Theta_{0}}\sqrt{\frac{-\theta_{0}}{\Theta_{0}-\theta_{0}}%
}\left(  \lambda-\gamma\Theta_{0}\right)  ^{k}\exp\left[  N\Psi_{0}\left(
y,\Theta_{0}\right)  \right]  \left[  \frac{\partial^{2}\Psi_{0}}%
{\partial\vartheta^{2}}\left(  y,\Theta_{0}\right)  \right]  ^{-\frac{1}{2}%
}\nonumber
\end{gather}
where
\[
\Psi_{0}\left(  y,\vartheta\right)  =y\vartheta+\mu(\vartheta)+V(\vartheta
)\ln\left[  -R_{1}(\vartheta)\right]  +\left[  1-V(\vartheta)\right]
\ln\left[  R_{2}(\vartheta)\right]
\]
and $\Theta_{0}(y)$ is the solution to the equation%
\begin{equation}
\frac{\partial\Psi_{0}}{\partial\vartheta}\left(  y,\Theta_{0}\right)  =0.
\label{Theta0}%
\end{equation}
We note from (\ref{Theta0}) that%
\begin{align*}
\Theta_{0}(y)  &  \rightarrow\Theta^{+}(0,0),\quad y\rightarrow0\\
\Theta_{0}(y)  &  \rightarrow\theta_{0},\quad y\rightarrow\infty
\end{align*}
and thus we shall limit $\Theta_{0}$ to the range $\Theta^{+}(0,0)\leq
\Theta_{0}<\theta_{0}$ in (\ref{F3.5}). Expansion (\ref{F3.5}) remains valid
for $y\rightarrow0$ and can be used to evaluate $F_{k}(0)-F_{k}(\infty)$ for
$k=O(1).$

\section{The boundary $z=1$ $\label{section z=1}$}

We examine the scale $k=N-O(1),$ which corresponds to $1-z=O(N^{-1}).$ The
analysis is different for three ranges of $y.$

\subsection{The boundary layer at $\ z=1,\quad0<y<Y_{0}(1)$ (Region IV)
$\label{z=1.1}$}

\begin{lemma}
\label{lemma 8}Let \ $k=N-j,\quad j=O(1),\quad\vartheta>0,\quad\vartheta
=O(1).$ Then, as $N\rightarrow\infty$%
\begin{equation}
h_{k}(\vartheta)\sim\frac{N^{j}}{j!}\left[  \frac{1+\left(  1-\gamma\right)
\vartheta}{\lambda}\right]  ^{j}. \label{hz=1.1}%
\end{equation}

\end{lemma}

\begin{proof}
From (\ref{h}) we have%
\[
h_{k}(\vartheta)=\frac{1}{2\pi i}%
{\displaystyle\oint\limits_{\mathcal{C}}}
\frac{H(\vartheta,w)}{w^{j+1}}dw.
\]
Changing variables to $\ u=wN$ we get%
\begin{align*}
h_{k}(\vartheta)  &  =N^{j}\frac{1}{2\pi i}%
{\displaystyle\oint\limits_{\mathcal{C}^{\prime}}}
\exp\left\{  NV(\vartheta)\ln\left[  1-R_{1}(\vartheta)\frac{u}{N}\right]
\right\} \\
&  \times\exp\left\{  N\left[  1-V(\vartheta)\right]  \ln\left[
1+R_{2}(\vartheta)\frac{u}{N}\right]  \right\}  \frac{du}{u^{j+1}}%
\end{align*}
\
\[
\sim N^{j}\frac{1}{2\pi i}%
{\displaystyle\oint\limits_{\mathcal{C}^{\prime}}}
\exp\left[  \frac{1+\left(  1-\gamma\right)  \vartheta}{\lambda}u\right]
\frac{du}{u^{j+1}},\quad N\rightarrow\infty
\]
where $\mathcal{C}^{\prime}=\frac{\mathcal{C}}{N}$ is a small loop around
$u=0.$ Hence, the result follows by expanding $\exp\left[  \frac{1+\left(
1-\gamma\right)  \vartheta}{\lambda}u\right]  $ in Taylor series.
\end{proof}

Using (\ref{P2}) and (\ref{hz=1.1}) in (\ref{Exact}) we have%
\[
F_{k}(x)\sim\frac{N^{j}}{j!}\left(  \frac{\lambda}{\lambda+1}\right)
^{N}\frac{1}{2\pi i}%
{\displaystyle\int\limits_{\mathrm{Br}^{+}}}
\frac{1}{\vartheta}\sqrt{\frac{-\theta_{0}}{\vartheta-\theta_{0}}}\left[
\frac{1+\left(  1-\gamma\right)  \vartheta}{\lambda}\right]  ^{j}\exp\left[
y\vartheta N+\mu(\vartheta)N\right]  d\vartheta.
\]
We evaluate the integral above for $N\rightarrow\infty$ by the saddle point
method and obtain%
\begin{align}
F_{k}(x)  &  \equiv F_{j}^{(4)}(y)\sim\frac{N^{j}}{j!}\frac{1}{\sqrt{2\pi N}%
}\left(  \frac{\lambda}{\lambda+1}\right)  ^{N}\frac{1}{\Theta_{1}}\sqrt
{\frac{-\theta_{0}}{\Theta_{1}-\theta_{0}}}\left[  \frac{1+\left(
1-\gamma\right)  \Theta_{1}}{\lambda}\right]  ^{j}\label{Fz-=1.1}\\
&  \times\exp\left[  y\Theta_{1}N+\mu(\Theta_{1})N\right]  \frac{1}{\sqrt
{\mu^{\prime\prime}(\Theta_{1})}}\nonumber
\end{align}
where $\Theta_{1}(y)$ is the solution to%
\begin{equation}
y+\mu^{\prime}(\Theta_{1})=0. \label{Eqtheta1}%
\end{equation}

\subsection{The corner layer at $\left(  Y_{0}(1),1\right)  $ (Region V)
\label{corner z=1}}

From (\ref{Eqtheta1}) we see that$\ \Theta_{1}(y)\rightarrow0$ \ as
\ $y\rightarrow Y_{0}(1).$ Therefore, we shall find a solution valid in the
neighborhood of the point $\left(  Y_{0}(1),1\right)  .$ From (\ref{hz=1.1})
we have as $N\rightarrow\infty,$ with $j=O(1)$%
\[
h_{k}(\vartheta)\sim\frac{1}{j!}\left(  \frac{N}{\lambda}\right)  ^{j}%
,\quad\vartheta\rightarrow0
\]
and from (\ref{mu}) we get%
\begin{align}
\mu(\vartheta)  &  \sim-\frac{\rho\left[  \ln\left(  \rho\right)  -1\right]
+1}{\left(  \lambda+1\right)  ^{2}}\vartheta+\left[  \frac{\left(
1-\gamma\right)  \left(  \rho-4\lambda+1\right)  }{2\left(  \lambda+1\right)
^{3}}-\frac{\zeta\ln\left(  \rho\right)  }{\left(  \lambda+1\right)  ^{4}%
}\right]  \vartheta^{2}\label{mu=0}\\
&  =-Y_{0}(1)\vartheta+\frac{1}{2}Y_{2}(1)\vartheta^{2},\quad\vartheta
\rightarrow0.\nonumber
\end{align}
Using (\ref{P2}) and the above in (\ref{Fexact}) we obtain%
\[
F_{k}(x)\sim\frac{1}{j!}\left(  \frac{N}{\lambda}\right)  ^{j}\left(
\frac{\lambda}{\lambda+1}\right)  ^{N}\frac{1}{2\pi i}%
{\displaystyle\int\limits_{\mathrm{Br}^{+}}}
\frac{1}{\vartheta}\exp\left[  y\vartheta N-Y_{0}(1)\vartheta N+\frac{1}%
{2}Y_{2}(1)\vartheta^{2}N\right]  d\vartheta.
\]
Defining%
\[
V(y,1)=\left[  y-Y_{0}(1)\right]  \sqrt{\frac{N}{Y_{2}(1)}}%
\]
we obtain, for $j=O(1)$ and $y-Y_{0}(1)=O\left(  N^{-\frac{1}{2}}\right)  $%
\begin{equation}
F_{k}(x)\equiv F_{j}^{\left(  5\right)  }(V)\sim\frac{1}{j!}\left(  \frac
{N}{\lambda}\right)  ^{j}\left(  \frac{\lambda}{\lambda+1}\right)  ^{N}%
\frac{1}{\sqrt{2\pi}}%
{\displaystyle\int\limits_{-\infty}^{V(y,1)}}
\exp\left(  -\frac{t^{2}}{2}\right)  dt. \label{Tranz=1}%
\end{equation}

\subsection{The boundary layer at $\ z=1,\quad y>Y_{0}(1)$ (Region VI)}

We shall examine the portion of the boundary \ $z=1$ for $y>Y_{0}(1).$ To do
so, we write $F_{k}(x)$ as%
\[
F_{k}(x)-F_{k}(\infty)=\left(  \frac{\lambda}{\lambda+1}\right)  ^{N}\frac
{1}{2\pi i}%
{\displaystyle\int\limits_{\mathrm{Br}^{-}}}
e^{x\vartheta}\frac{1}{\vartheta}\left[
{\displaystyle\prod\limits_{j=0}^{N-\left\lfloor c\right\rfloor -1}}
\frac{\theta_{j}}{\theta_{j}-\vartheta}\right]  h_{k}(\vartheta)d\vartheta
\]
and, repeating the same calculation done in Section \ref{z=1.1}, we conclude
that
\begin{equation}
F_{k}(x)-F_{k}(\infty)\equiv F_{j}^{\left(  6\right)  }(y)-\left(
\frac{\lambda}{1+\lambda}\right)  ^{N}\binom{N}{j}\lambda^{-j}\sim F_{j}%
^{(4)}(y). \label{F6}%
\end{equation}

\section{The boundary \ $y=0$ $\label{section y=0}$}

In this Section we analyze the part of the boundary where the boundary
condition (\ref{zero}) applies. Theorem \ref{theorem1} is not valid for small
$y$ since for $\gamma<z\leq1$ the saddle $\Theta(y,z)\rightarrow\infty$ as
$y\rightarrow0^{+},$ but our approximation to the integrand in (\ref{Exact})
was only valid for $\vartheta=O(1).$ We will scale $y=x/N$ and use Lemma
\ref{lemma5}, which applies for $\vartheta=SN=O(N).$ We note that
$y\vartheta=xS.$ We will need to examine separately the cases $\gamma<z<1$ and
$1-z=O(\varepsilon).$

\subsection{The boundary layer at $y=0,\quad\gamma<z<1$ (Region VII)}

\begin{lemma}
Let \ $k=Nz,\quad\gamma<z<1,\quad\vartheta=SN,\quad S>0,\quad S=O(1).$ Then,%
\begin{align}
h_{k}(\vartheta)  &  \sim\frac{1}{\sqrt{2\pi N}}\exp\left\{  \left(
1-z\right)  \ln\left[  \frac{\left(  z-\gamma\right)  S}{\lambda
(1-z)}N\right]  N+\left(  1-\gamma\right)  \ln\left(  \frac{1-\gamma}%
{z-\gamma}\right)  N\right\} \label{hy=0}\\
&  \times\sqrt{\frac{\left(  1-\gamma\right)  }{\left(  1-z\right)  \left(
z-\gamma\right)  }}\exp\left[  \frac{\phi}{S}\ln\left(  \frac{1-\gamma
}{z-\gamma}\right)  +\frac{(1-z)\left(  1-\lambda\right)  }{S}\right]
.\nonumber
\end{align}

\end{lemma}

\begin{proof}
Changing variables in (\ref{h}) from $w$ to \ $u=wSN,$ we have%
\begin{align*}
h_{k}(\vartheta)  &  =\frac{1}{2\pi i}%
{\displaystyle\oint\limits_{\mathcal{C}^{\prime}}}
\exp\left\{  V(SN)\ln\left[  1-R_{1}(SN)\frac{u}{SN}\right]  N-(1-z)\ln\left(
\frac{u}{SN}\right)  N\right\} \\
&  \times\exp\left\{  \left[  1-V(SN)\right]  \ln\left[  1+R_{2}(SN)\frac
{u}{SN}\right]  N\right\}  \frac{du}{u}%
\end{align*}%
\begin{align*}
&  \sim\frac{1}{2\pi i}%
{\displaystyle\oint\limits_{\mathcal{C}^{\prime}}}
\exp\left[  \left(  1-z\right)  \ln\left(  \frac{SN}{u}\right)  N+\left(
1-\gamma\right)  \ln\left(  1+\frac{u}{\lambda}\right)  N\right] \\
&  \times\exp\left[  \frac{\phi}{S}\ln\left(  1+\frac{u}{\lambda}\right)
+\frac{u\left(  1-\gamma\right)  \left(  1-\lambda\right)  }{\left(
\lambda+u\right)  S}\right]  \frac{du}{u}.
\end{align*}
Use of the saddle point method gives, for $N\rightarrow\infty$ with $S$ and
$z$ fixed,%
\begin{align*}
h_{k}(\vartheta)  &  \sim\frac{\lambda}{\sqrt{2\pi N}}\frac{1}{u_{1}}%
\exp\left[  \left(  1-z\right)  \ln\left(  \frac{SN}{u_{1}}\right)  N+\left(
1-\gamma\right)  \ln\left(  1+\frac{u_{1}}{\lambda}\right)  N\right] \\
&  \times\exp\left[  \frac{\phi}{S}\ln\left(  1+\frac{u_{1}}{\lambda}\right)
+\frac{u_{1}\left(  1-\gamma\right)  \left(  1-\lambda\right)  }{\left(
\lambda+u_{1}\right)  S}\right]  \sqrt{\frac{\left(  1-\gamma\right)  \left(
1-z\right)  }{\left(  z-\gamma\right)  ^{3}}}%
\end{align*}
where the saddle point is located at%
\[
u_{1}=\frac{\lambda(1-z)}{z-\gamma}.
\]

\end{proof}

Using (\ref{Paprox}) and (\ref{hy=0}) in (\ref{Exact}), we get%
\begin{gather}
F_{k}(x)\sim\frac{1}{2\pi N}\left(  \frac{\lambda}{\lambda+1}\right)
^{N}\sqrt{\frac{\rho}{\phi\gamma\left(  1-z\right)  \left(  z-\gamma\right)
}}\exp\left[  -\gamma\ln(\gamma)N\right] \nonumber\\
\times\frac{1}{2\pi i}%
{\displaystyle\int\limits_{\mathrm{Br}^{\prime}}}
\exp\left\{  xSN-\left(  1-\gamma\right)  \ln\left[  \left(  z-\gamma\right)
SN\right]  N+\left(  1-z\right)  \ln\left[  \frac{\left(  z-\gamma\right)
S}{\lambda\left(  1-z\right)  }N\right]  N\right\} \label{F7.5}\\
\times\exp\left\{  \frac{\phi}{S}\ln\left[  \frac{\gamma}{\phi\left(
z-\gamma\right)  N}\right]  +\frac{2\phi-\rho-1+(1-z)\left(  1-\lambda\right)
}{S}\right\} \nonumber\\
\times\frac{1}{S}\Gamma\left(  \frac{\phi}{S}+1-\alpha\right)  \left(
\frac{\phi}{S}\right)  ^{\alpha}dS.\nonumber
\end{gather}
where $\mathrm{Br}^{\prime}$ is a vertical contour on which $\operatorname{Re}%
(S)>0.$ Computing the integral in (\ref{F7.5}) as $N\rightarrow\infty$ by the
saddle point method, we find that the integrand has a saddle at%
\begin{equation}
S^{\ast}(x,z)=\frac{z-\gamma}{x} \label{Sstar}%
\end{equation}
and hence the leading term for (\ref{F7.5}) is
\begin{gather*}
F_{k}(x)\sim\left(  \frac{1}{2\pi N}\right)  ^{\frac{3}{2}}\left(
\frac{\lambda}{\lambda+1}\right)  ^{N}\sqrt{\frac{\rho}{\phi\gamma\left(
1-z\right)  }}\frac{1}{z-\gamma}\exp\left[  -\gamma\ln(\gamma)N\right] \\
\times\exp\left\{  \left(  z-\gamma\right)  N-\left(  1-\gamma\right)
\ln\left[  \frac{\left(  z-\gamma\right)  ^{2}}{x}N\right]  N+\left(
1-z\right)  \ln\left[  \frac{\left(  z-\gamma\right)  ^{2}}{\lambda\left(
1-z\right)  x}N\right]  N\right\} \\
\times\exp\left\{  \frac{x\phi}{z-\gamma}\ln\left[  \frac{\gamma}{\phi\left(
z-\gamma\right)  N}\right]  +\frac{2\phi-\rho-1+(1-z)\left(  1-\lambda\right)
}{z-\gamma}x\right\} \\
\times\Gamma\left(  \frac{x\phi}{z-\gamma}+1-\alpha\right)  \left(
\frac{x\phi}{z-\gamma}\right)  ^{\alpha}.
\end{gather*}
After some simplification we get%
\begin{gather}
F_{k}(x)\equiv F^{(7)}(x,z)\sim\left(  \frac{1}{2\pi N}\right)  ^{\frac{3}{2}%
}\sqrt{\frac{\rho}{\phi\gamma\left(  1-z\right)  }}\frac{1}{z-\gamma}%
\Gamma\left(  \frac{x\phi}{z-\gamma}+1-\alpha\right)  \left(  \frac{x\phi
}{z-\gamma}\right)  ^{\alpha}\nonumber\\
\times\exp\left\{  (z-\gamma)\ln\left[  \frac{xe}{\left(  z-\gamma\right)
^{2}N}\right]  N+z\ln(\lambda)N-\left(  1-z\right)  \ln\left(  1-z\right)
N-\ln(\lambda+1)N-\gamma\ln(\gamma)N\right\} \label{Fx=0}\\
\times\exp\left\{  \frac{x\phi}{z-\gamma}\ln\left[  \frac{\gamma}{\phi\left(
z-\gamma\right)  N}\right]  +\left[  \frac{2\lambda\left(  1-\gamma\right)
}{z-\gamma}+\left(  \lambda-1\right)  \right]  x\right\}  .\nonumber
\end{gather}
Note that (\ref{Fx=0}) is singular as $z\rightarrow\gamma.$ We can show that
if we expand (\ref{Fx=0}) as $(x,z)\rightarrow(0,\gamma),$ with $\frac
{x}{z-\gamma}$ fixed, the result asymptotically matches to he corner layer
approximation in Theorem \ref{corner}. The latter must be expanded for
$\chi,l\rightarrow\infty$ with $\frac{l}{\chi}$ fixed.

\subsection{The corner layer at $\left(  0,1\right)  $ (Region VIII)}

When $z\rightarrow1,$ (\ref{Fx=0}) becomes singular. Therefore, we shall find
an approximation to \ $F_{k}(x)$ \ close to the point $(0,1).$ We first
observe that when
\[
k=N-j,\quad\vartheta=SN,\quad S>0,\quad j,\ S=O(1)
\]
we obtain from (\ref{h}), by a calculation similar to that in the proof of
Lemma \ref{lemma 8},%
\begin{equation}
h_{k}(\vartheta)\sim N^{2j}\frac{1}{j!}\left[  \frac{\left(  1-\gamma\right)
S}{\lambda}\right]  ^{j}. \label{h(0,1)}%
\end{equation}
Using (\ref{Paprox}) and (\ref{h(0,1)}) in (\ref{Fexact}) we have%
\begin{align}
F_{k}(x)  &  \sim N^{2j}\frac{1}{j!}\left(  \frac{\lambda}{\lambda+1}\right)
^{N}\frac{1}{\sqrt{2\pi N}}\sqrt{\frac{\rho}{\phi\gamma\left(  1-\gamma
\right)  }}\nonumber\\
&  \times\frac{1}{2\pi i}%
{\displaystyle\int\limits_{\mathrm{Br}^{+}}}
\exp\left\{  xSN-\left(  1-\gamma\right)  \ln\left[  \left(  1-\gamma\right)
SN\right]  N-\gamma\ln(\gamma)N\right\} \label{F8.2}\\
&  \times\frac{1}{S}\left(  \frac{\phi}{S}\right)  ^{\alpha}\Gamma\left(
\frac{\phi}{S}+1-\alpha\right)  \left[  \frac{\left(  1-\gamma\right)
S}{\lambda}\right]  ^{j}\nonumber\\
&  \times\exp\left\{  \frac{\phi}{S}\ln\left[  \frac{\gamma}{\phi\left(
1-\gamma\right)  N}\right]  +\frac{2\phi-\rho-1}{S}\right\}  dS.\nonumber
\end{align}
Using the saddle point method we find that the integrand of (\ref{F8.2}) has a
saddle at%
\[
S^{\ast}(x,1)=\frac{1-\gamma}{x}%
\]
where $S^{\ast}(x,z)$ was defined in (\ref{Sstar}). Hence, to leading order,
we obtain
\begin{align}
F_{k}(x)  &  \equiv F_{j}^{(8)}(x)\sim N^{2j}\frac{1}{j!}\left(  \frac
{\lambda}{\lambda+1}\right)  ^{N}\frac{1}{2\pi N}\sqrt{\frac{\rho}{\phi\gamma
}}\frac{1}{1-\gamma}\nonumber\\
&  \times\exp\left\{  \left(  1-\gamma\right)  \ln\left[  \frac{ex}{\left(
1-\gamma\right)  ^{2}N}\right]  N-\gamma\ln(\gamma)N\right\} \label{F8}\\
&  \times\left(  \frac{x\phi}{1-\gamma}\right)  ^{\alpha}\Gamma\left(
\frac{x\phi}{1-\gamma}+1-\alpha\right)  \left[  \frac{\left(  1-\gamma\right)
^{2}}{\lambda x}\right]  ^{j}\nonumber\\
&  \times\exp\left\{  \frac{x\phi}{1-\gamma}\ln\left[  \frac{\gamma}%
{\phi\left(  1-\gamma\right)  N}\right]  +\left(  3\lambda-1\right)
x\right\}  .\nonumber
\end{align}
Note that from (\ref{Fx=0}) and (\ref{F8}) we have%
\[
F_{k}(x)=O\left(  x^{\alpha+(z-\gamma)N}\right)  =O\left(  x^{k-\left\lfloor
c\right\rfloor }\right)  ,\quad x\rightarrow0,\quad\left\lfloor c\right\rfloor
+1\leq k\leq N
\]
and thus we recover the result found in \cite{MR84a:68020} (equation (38)
therein), that
\[
\frac{d^{n}F_{k}}{dx^{n}}(0)=0,\quad\left\lfloor c\right\rfloor +1+n\leq k\leq
N.
\]

\section{Summary and numerical studies}

\begin{figure}[t]
\begin{center}
\rotatebox{0} {\resizebox{15cm}{!}{\includegraphics{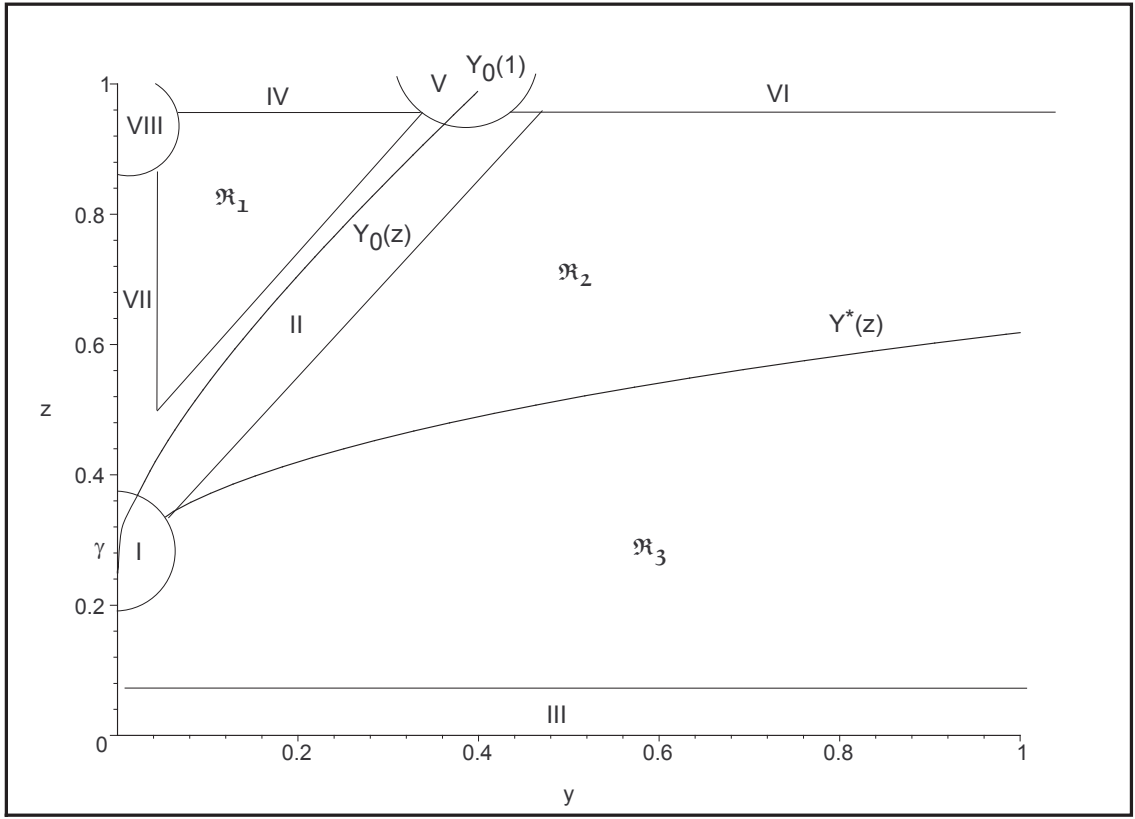}}}
\end{center}
\caption{A sketch of the different asymptotic regions.}%
\label{figregions1}%
\end{figure}

In most of the strip $\mathfrak{D=}\left\{  (y,z):y\geq0,\ 0\leq
z\leq1\right\}  ,$ the asymptotic expansion of $F_{k}(x)$ is given by Theorem
\ref{theorem1}. Below we summarize our results for the various boundary,
corner and transition layer regions, where Theorem \ref{theorem1} does not
apply. The paragraph number refers to the corresponding region (see Figure
\ref{figregions1}).

\begin{theorem}
As $N\rightarrow\infty$, asymptotic approximations of the joint distribution
$F_{k}(x)$ are as follows.

\begin{enumerate}
\item[I.] $k=c-\alpha+l=\left\lfloor c\right\rfloor +l,\quad-\infty
<l<\infty,\quad l=O(1)$ and $x=\frac{\chi}{N},\quad\chi=O(1),$ $F_{k}(x)\sim
F_{l}^{(1)}(\chi)$ as $N\rightarrow\infty$ where%
\begin{align}
F_{l}^{(1)}(\chi)  &  =\frac{1}{\sqrt{2\pi N}}\sqrt{\frac{\rho}{\gamma
\phi\left(  1-\gamma\right)  }}\left(  \frac{\beta}{2\gamma}\right)
^{l-\alpha}\exp\left[  N\Phi(\gamma)\right] \nonumber\\
&  \times\frac{1}{2\pi i}%
{\displaystyle\int\limits_{\mathrm{Br}^{+}}}
e^{\chi S}\frac{1}{S}\Gamma\left(  \frac{\phi}{S}+1-\alpha\right)
J_{l-\alpha+\frac{\phi}{S}}\left(  \frac{\beta}{S}\right) \label{I}\\
&  \times\exp\left\{  \frac{\phi}{S}\ln\left[  \frac{\beta S}{2\lambda
\phi\left(  1-\gamma\right)  }\right]  +\frac{\phi-\rho}{S}\right\}  \left(
\frac{\phi}{S}\right)  ^{\alpha}dS,\nonumber
\end{align}
$\mathrm{Br}^{+}$ is a vertical contour in the complex plane with
$\operatorname{Re}(S)>0,$ $J_{\cdot}(\cdot)$ is the Bessel function,
$\Gamma\left(  \cdot\right)  $ is the Gamma function ,%
\[
\beta=2\sqrt{\lambda\gamma\left(  1-\gamma\right)  }%
\]
and%
\begin{equation}
\Phi(z)=-z\ln(z)-(1-z)\ln\left(  1-z\right)  +z\ln\left(  \lambda\right)
-\ln\left(  \lambda+1\right)  .
\end{equation}

\item[II.] $y-Y_{0}(z)=O\left(  N^{-\frac{1}{2}}\right)  ,\quad\gamma<z<1$%
\begin{equation}
F_{k}(x)=F^{(2)}(V,z)\sim\frac{1}{\sqrt{2\pi N}}\sqrt{\frac{1}{z\left(
1-z\right)  }}\exp\left[  \Phi(z)N\right]  \frac{1}{\sqrt{2\pi}}%
{\displaystyle\int\limits_{-\infty}^{V}}
\exp\left(  -\frac{t^{2}}{2}\right)  dt \label{II}%
\end{equation}
with%
\[
V(y,z)=\frac{y-Y_{0}(z)}{\sqrt{Y_{2}(z)}}\sqrt{N}%
\]%
\[
Y_{0}(z)=\frac{z-\gamma}{\lambda+1}-\frac{\rho}{\left(  \lambda+1\right)
^{2}}\ln\left(  \frac{z\lambda+z-\lambda}{\rho}\right)  ,\quad\gamma<z<1
\]%
\begin{gather*}
Y_{2}(z)=\frac{2\zeta}{\left(  \lambda+1\right)  ^{4}}\ln\left(
\frac{z+z\lambda-\lambda}{\rho}\right) \\
-\frac{z-\gamma}{\left(  \lambda+1\right)  \left(  \lambda z+z-\lambda\right)
^{2}}\left[  \frac{2\zeta\rho}{\left(  \lambda+1\right)  ^{2}}+\frac{3\zeta
}{\left(  \lambda+1\right)  }\left(  z-\gamma\right)  +\left(  \lambda
-1\right)  \left(  z-\gamma\right)  ^{2}\right]
\end{gather*}%
\[
\zeta=2\lambda-\gamma+(\gamma-1)\lambda^{2}.
\]

\item[III.] $k=O(1)$%
\begin{gather}
F_{k}(x)-F_{k}(\infty)=F_{k}^{(3)}(y)-F_{k}(\infty)\sim\frac{N^{k}}{k!}\left(
\frac{\lambda}{\lambda+1}\right)  ^{N}\frac{1}{\sqrt{2\pi N}}\label{III}\\
\times\frac{1}{\Theta_{0}}\sqrt{\frac{-\theta_{0}}{\Theta_{0}-\theta_{0}}%
}\left(  \lambda-\gamma\Theta_{0}\right)  ^{k}\exp\left[  N\Psi_{0}\left(
y,\Theta_{0}\right)  \right]  \left[  \frac{\partial^{2}\Psi_{0}}%
{\partial\vartheta^{2}}\left(  y,\Theta_{0}\right)  \right]  ^{-\frac{1}{2}%
}\nonumber
\end{gather}
where
\[
\Psi_{0}\left(  y,\vartheta\right)  =y\vartheta+\mu(\vartheta)+V(\vartheta
)\ln\left[  -R_{1}(\vartheta)\right]  +\left[  1-V(\vartheta)\right]
\ln\left[  R_{2}(\vartheta)\right]
\]
and $\Theta_{0}(y)$ is the solution to the equation%
\begin{equation}
\frac{\partial\Psi_{0}}{\partial\vartheta}\left(  y,\Theta_{0}\right)  =0.
\end{equation}

\item[IV.] $k=N-j,\quad\ j=O(1),\quad0<y<Y_{0}(1)$%
\begin{align}
F_{k}(x)  &  =F_{j}^{(4)}(y)\sim\frac{N^{j}}{j!}\frac{1}{\sqrt{2\pi N}}\left(
\frac{\lambda}{\lambda+1}\right)  ^{N}\frac{1}{\Theta_{1}}\sqrt{\frac
{-\theta_{0}}{\Theta_{1}-\theta_{0}}}\left[  \frac{1+\left(  1-\gamma\right)
\Theta_{1}}{\lambda}\right]  ^{j}\label{IV}\\
&  \times\exp\left[  y\Theta_{1}N+\mu(\Theta_{1})N\right]  \frac{1}{\sqrt
{\mu^{\prime\prime}(\Theta_{1})}}\nonumber
\end{align}
where $\Theta_{1}(y)$ is the solution to%
\begin{equation}
y+\mu^{\prime}(\Theta_{1})=0.
\end{equation}

\item[V.] $k=N-j,\quad\ j=O(1),\quad y-Y_{0}(1)=O\left(  N^{-\frac{1}{2}%
}\right)  $%
\begin{equation}
F_{k}(x)=F_{j}^{\left(  5\right)  }(V)\sim\frac{1}{j!}\left(  \frac{N}%
{\lambda}\right)  ^{j}\left(  \frac{\lambda}{\lambda+1}\right)  ^{N}\frac
{1}{\sqrt{2\pi}}%
{\displaystyle\int\limits_{-\infty}^{V(y,1)}}
\exp\left(  -\frac{t^{2}}{2}\right)  dt \label{V(5)}%
\end{equation}
with%
\[
V(y,1)=\left[  y-Y_{0}(1)\right]  \sqrt{\frac{N}{Y_{2}(1)}}.
\]

\item[VI.] $k=N-j,\quad\ j=O(1),\quad y>Y_{0}(1)$%
\begin{equation}
F_{k}(x)-F_{k}(\infty)=F_{j}^{\left(  6\right)  }(y)-\left(  \frac{\lambda
}{1+\lambda}\right)  ^{N}\binom{N}{j}\lambda^{-j}\sim F_{j}^{(4)}(y).
\label{VI}%
\end{equation}
where $F_{j}^{\left(  4\right)  }(y)$ is as in item IV.

\item[VII.] $x=O(1),\quad\gamma<z<1$%
\begin{gather}
F_{k}(x)=F^{(7)}(x,z)\sim\left(  \frac{1}{2\pi N}\right)  ^{\frac{3}{2}}%
\sqrt{\frac{\rho}{\phi\gamma\left(  1-z\right)  }}\frac{1}{z-\gamma}%
\Gamma\left(  \frac{x\phi}{z-\gamma}+1-\alpha\right)  \left(  \frac{x\phi
}{z-\gamma}\right)  ^{\alpha}\nonumber\\
\times\exp\left\{  (z-\gamma)\ln\left[  \frac{xe}{\left(  z-\gamma\right)
^{2}N}\right]  N+z\ln(\lambda)N-\left(  1-z\right)  \ln\left(  1-z\right)
N-\ln(\lambda+1)N-\gamma\ln(\gamma)N\right\} \label{VII}\\
\times\exp\left\{  \frac{x\phi}{z-\gamma}\ln\left[  \frac{\gamma}{\phi\left(
z-\gamma\right)  N}\right]  +\left[  \frac{2\lambda\left(  1-\gamma\right)
}{z-\gamma}+\left(  \lambda-1\right)  \right]  x\right\}  .\nonumber
\end{gather}

\item[VIII.] $k=N-j,\quad\ j=O(1),\quad x=O(1)$%
\begin{align}
F_{k}(x)  &  =F_{j}^{(8)}(x)\sim N^{2j}\frac{1}{j!}\left(  \frac{\lambda
}{\lambda+1}\right)  ^{N}\frac{1}{2\pi N}\sqrt{\frac{\rho}{\phi\gamma}}%
\frac{1}{1-\gamma}\nonumber\\
&  \times\exp\left\{  \left(  1-\gamma\right)  \ln\left[  \frac{ex}{\left(
1-\gamma\right)  ^{2}N}\right]  N-\gamma\ln(\gamma)N\right\} \label{VIII}\\
&  \times\left(  \frac{x\phi}{1-\gamma}\right)  ^{\alpha}\Gamma\left(
\frac{x\phi}{1-\gamma}+1-\alpha\right)  \left[  \frac{\left(  1-\gamma\right)
^{2}}{\lambda x}\right]  ^{j}\nonumber\\
&  \times\exp\left\{  \frac{x\phi}{1-\gamma}\ln\left[  \frac{\gamma}%
{\phi\left(  1-\gamma\right)  N}\right]  +\left(  3\lambda-1\right)
x\right\}  .\nonumber
\end{align}

\end{enumerate}
\end{theorem}

From (\ref{Exact}) we see that the continuous part of the density (for $x>0)$
is given by
\begin{equation}
f_{k}(x)\equiv F_{k}^{\prime}(x)=\left(  \frac{\lambda}{\lambda+1}\right)
^{N}\frac{1}{2\pi i}%
{\displaystyle\int\limits_{\mathrm{Br}^{+}}}
e^{x\vartheta}\left[
{\displaystyle\prod\limits_{j=0}^{N-\left\lfloor c\right\rfloor -1}}
\frac{\theta_{j}}{\theta_{j}-\vartheta}\right]  h_{k}(\vartheta)d\vartheta.
\label{density}%
\end{equation}
Therefore, we can approximate the density by multiplying the asymptotic
approximations by the corresponding saddle point:
\begin{align*}
f_{k}(x)  &  \sim\Theta G_{1}(y,z)\\
f_{k}(x)  &  \sim\Theta^{+}G_{2}(y,z)\\
f_{k}(x)  &  \sim\Theta_{0}\left[  F_{k}^{(3)}(y)-F_{k}(\infty)\right] \\
f_{k}(x)  &  \sim\Theta_{1}F_{j}^{(4)}(y)\\
f_{k}(x)  &  \sim S^{\ast}(x,z)NF^{(7)}(x,z).
\end{align*}
Here the right sides are the expansions valid in the various regions.

\begin{figure}[ptb]
\begin{center}
\rotatebox{270} {\resizebox{11cm}{!}{\includegraphics{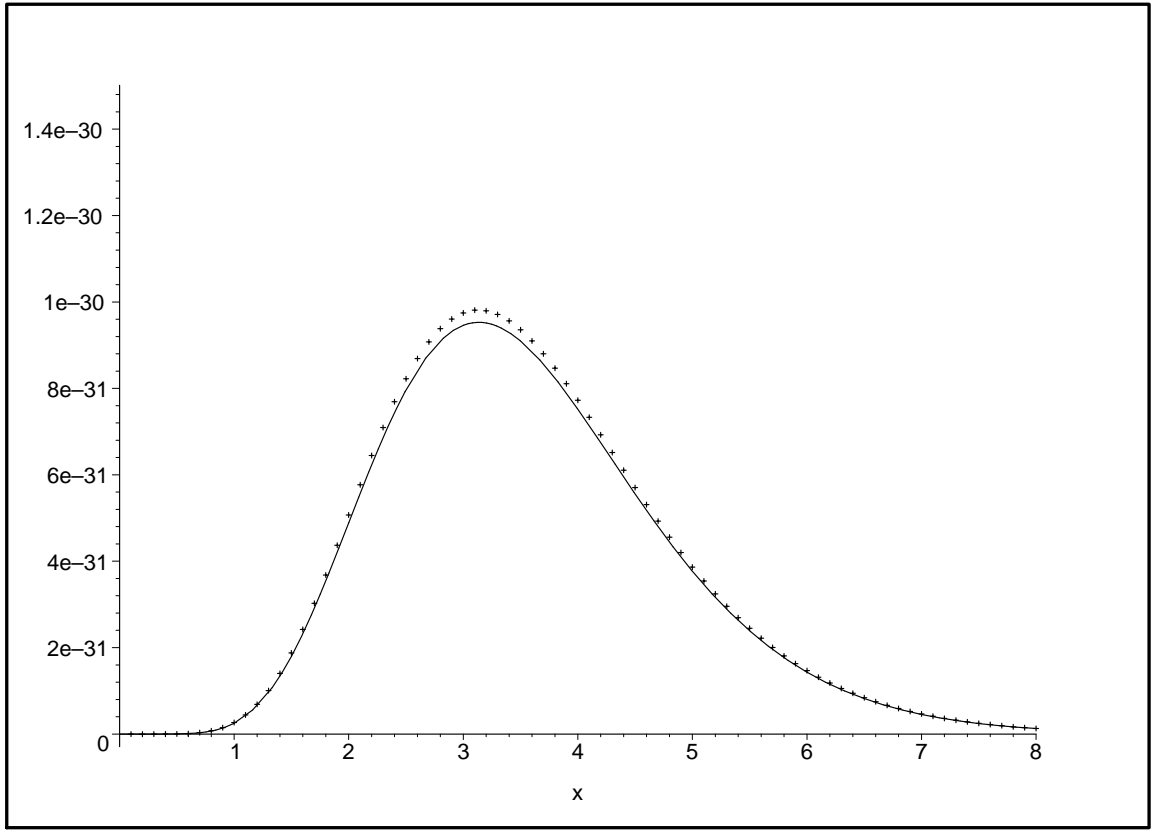}}}
\end{center}
\caption{A comparison of the exact (solid curve) and asymptotic (+++) values
of the density $f_{k}(x)$ for $k=17$ with $N=20$.}%
\label{figcompare1}%
\end{figure}

We compare our asymptotic approximations to the density with the exact results
obtained by evaluating (\ref{density}). We use $N=20,\ \gamma=0.37987897\ $and
$\lambda=0.0122448.$ In Figure \ref{figcompare1} we graph $f_{k}(x)$ and the
asymptotic approximation $\Theta G_{1}(y,z)$, for $k=17.$ We note that the
asymptotics predict that the maximum value will be achieved at $x=NY_{0}%
\left(  \frac{k}{N}\right)  \simeq3.052.$

\begin{figure}[ptb]
\begin{center}
\rotatebox{270} {\resizebox{11cm}{!}{\includegraphics{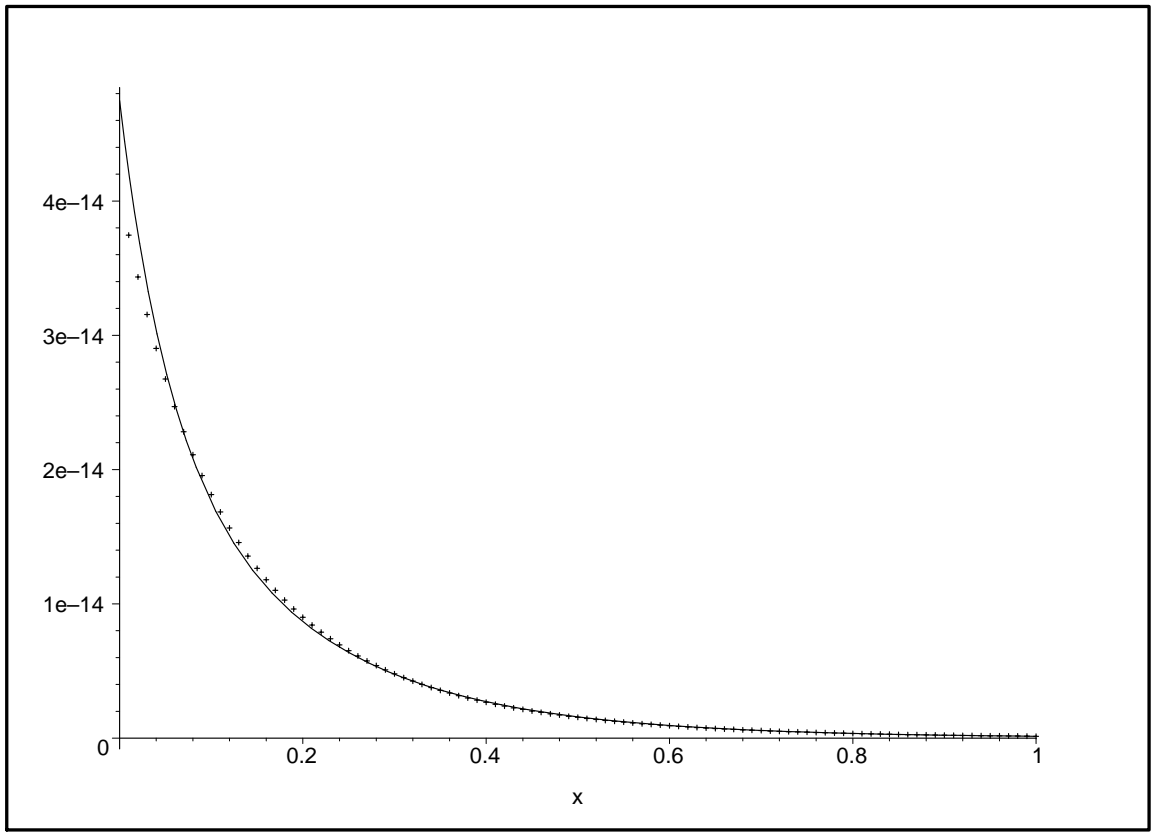}}}
\end{center}
\caption{A comparison of the exact (solid curve) and asymptotic (+++) values
of the density $f_{k}(x)$ for $k=3$ with $N=20$.}%
\label{figcompare3}%
\end{figure}

In Figure \ref{figcompare3} we graph $f_{3}(x)$ and the asymptotic
approximation $\Theta^{+}G_{2}(y,z).$ We note that the maximum is now achieved
at $x=0,$ which is consistent with the asymptotic analysis.

\begin{figure}[ptb]
\begin{center}
\rotatebox{270} {\resizebox{11cm}{!}{\includegraphics{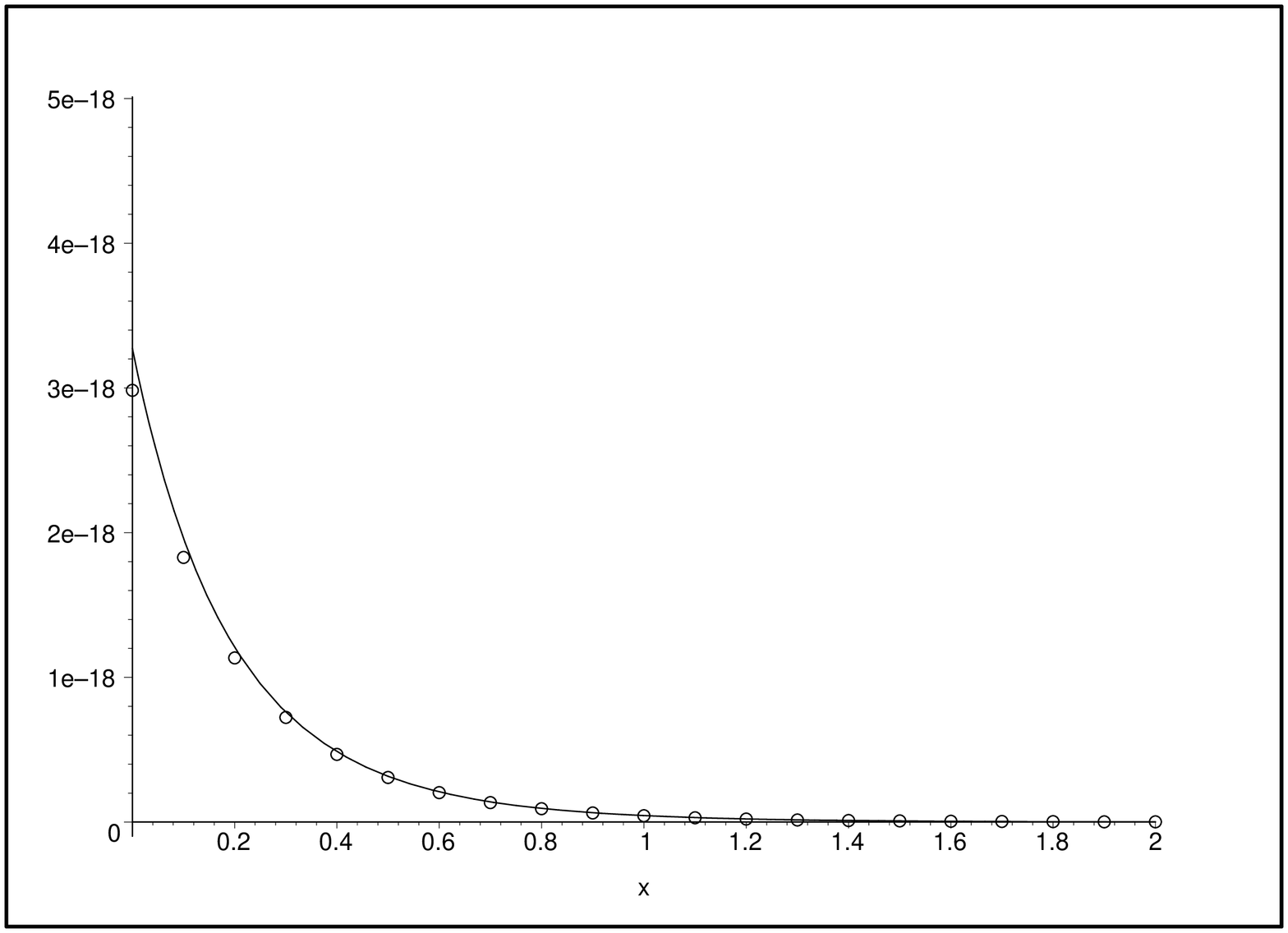}}}
\end{center}
\caption{A comparison of the exact (solid curve) and asymptotic ($\circ
\circ\circ$) values of the density $f_{k}(x)$ for $k=0$ with $N=20$.}%
\label{figcompare5}%
\end{figure}

In Figure \ref{figcompare5} we graph $f_{0}(x)$ and $\Theta_{0}\left[
F_{0}^{(3)}(y)-F_{0}(\infty)\right]  .$ This shows the importance of treating
Region III.

\begin{figure}[ptb]
\begin{center}
\rotatebox{270} {\resizebox{11cm}{!}{\includegraphics{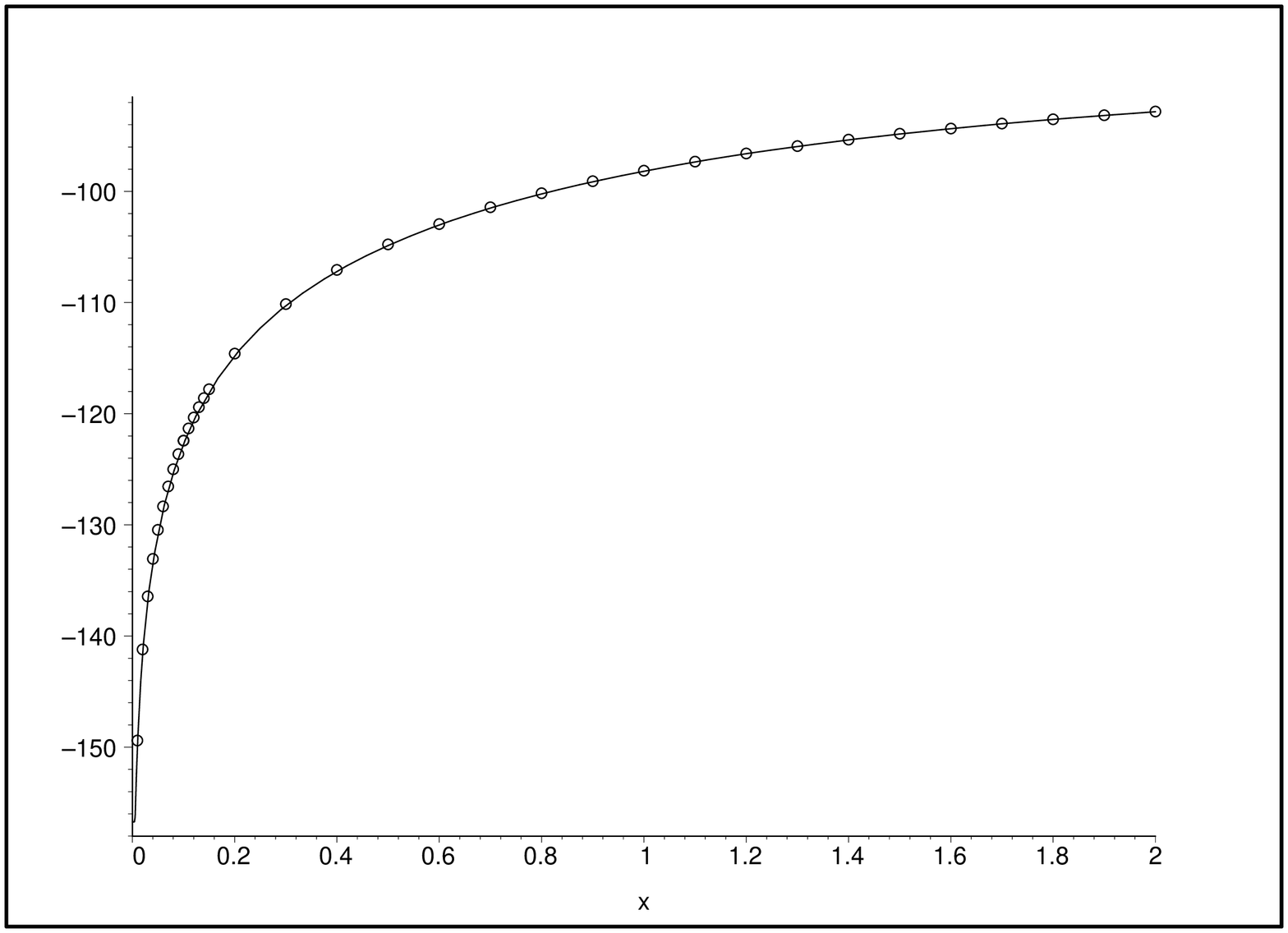}}}
\end{center}
\caption{A comparison of the exact (solid curve) and asymptotic ($\circ
\circ\circ$) values of $\ln\left[  f_{k}(x)\right]  $ for $z=1$ with $N=20$.}%
\label{figcompare6}%
\end{figure}

In Figure \ref{figcompare6} we graph $\ln\left[  f_{20}(x)\right]  $ and
$\ln\left[  \Theta_{1}F_{0}^{(4)}(y)\right]  ,$ the latter being the boundary
layer approximation to the density at $z=1.$

\begin{figure}[ptb]
\begin{center}
\rotatebox{270} {\resizebox{11cm}{!}{\includegraphics{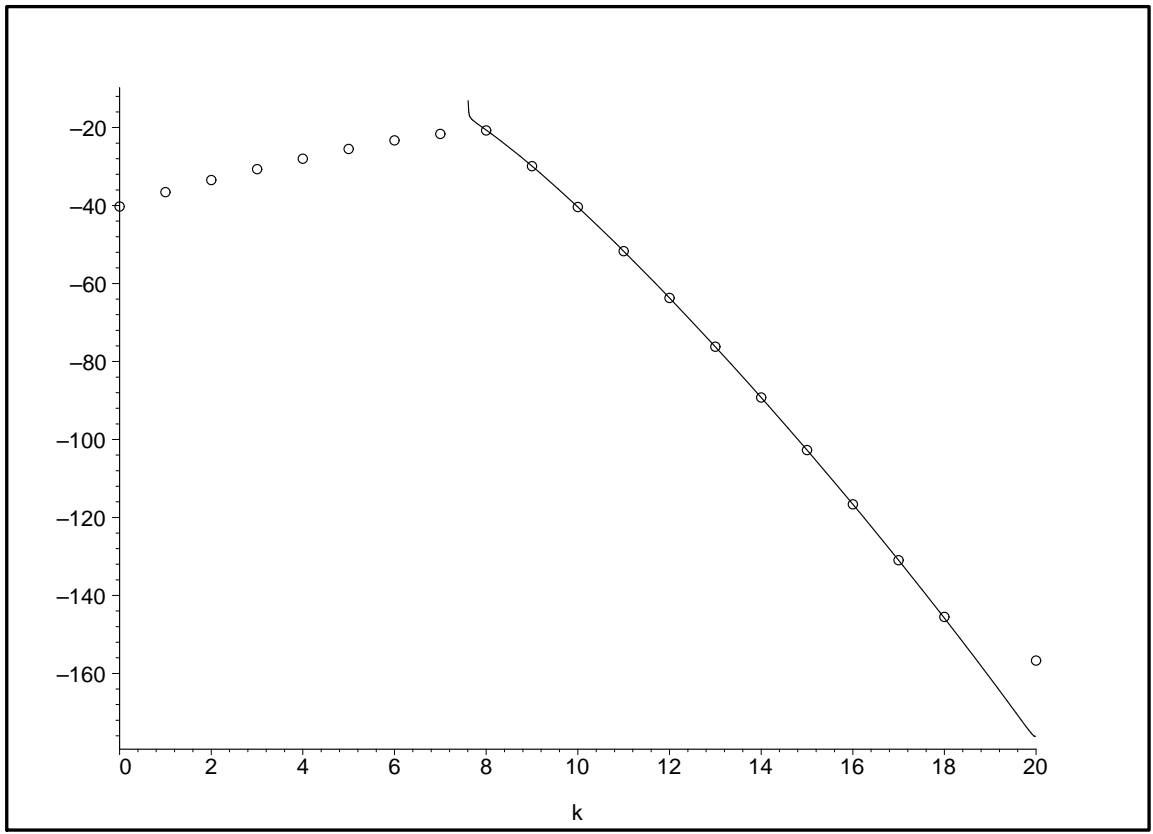}}}
\end{center}
\caption{A comparison of the exact ($\circ\circ\circ$) and asymptotic (solid
curve) values of $\ln\left[  f_{k}(x)\right]  $ for $x=0.001$ with $N=20$.}%
\label{figcompare7}%
\end{figure}

In Figure \ref{figcompare7} we graph $\ln\left[  f_{k}(x)\right]  $ and
$\ln\left[  S^{\ast}(x,z)NF^{(7)}(x,z)\right]  $ with $x=0.001.$ We see how
the asymptotic approximation breaks down for $z\approx\gamma$ and for
$z\approx1$ (it does not apply for $z<\gamma$).

\begin{figure}[ptb]
\begin{center}
\rotatebox{270} {\resizebox{11cm}{!}{\includegraphics{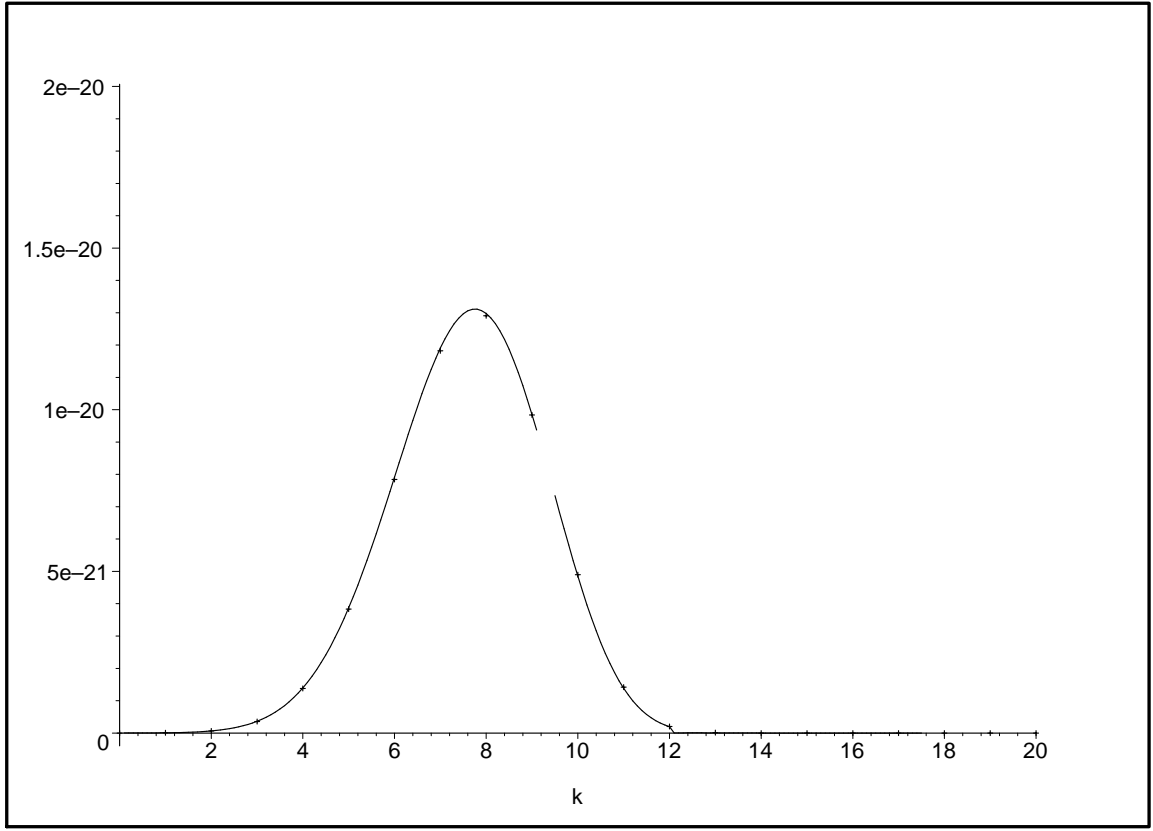}}}
\end{center}
\caption{A comparison of the exact (+++) and asymptotic (solid curve) values
of the density $f_{k}(x)$ for $x=5$ with $N=20$.}%
\label{figcompare4}%
\end{figure}

In \ref{figcompare4} we graph $f_{k}(x)$ and the asymptotic approximations
$\Theta G_{1}(y,z)$ and $\Theta^{+}G_{2}(y,z)$ with $x=5.$ To show in more
detail the different ranges, in Table 1 we compare the exact value of the
density with the approximations given by $\Theta_{0}F_{k}^{(3)}(y)$ and
$\Theta^{+}G_{2}(y,z)$, for $x=1$ ($y=0.05$) with $0\leq k\leq9.$ In Table 2
we compare $f_{k}(x)$ with the approximations given by $\Theta_{1}F_{j}%
^{(4)}(y)$ and $\Theta G_{1}(y,z)$, for $x=1$ ($y=0.05$) with $10\leq k\leq
N$. We note that for our choice of parameters $c\simeq7.598$ and
$Y_{0}(1)\simeq.4998.$ Also, when $y=1/N=0.05$, the equation $Y^{\ast}(z)=y$
has the solution $z\simeq.4811$, so that the transition between $\mathfrak{R}%
_{2}$ and $\mathfrak{R}_{3}$ occurs at $Nz\simeq9.622$.

\begin{table}[ptb]
\caption{}
\begin{center}
\medskip%
\begin{tabular}
[c]{|c|c|c|c|c|}\hline
$k$ & $\Theta$ & Exact & $\Theta_{0}F_{k}^{(3)}(y)$ & $\Theta^{+}G_{2}%
(y,z)$\\\hline
$0$ & $-3.418$ & $.445410^{-19}$ & $.432410^{-19}$ & \\
\ $1$ & $-3.593$ & $.125510^{-17}$ & $.113310^{-17}$ & $.133410^{-17}$\\
\ $2$ & $-3.785$ & $.168610^{-16}$ & $.148610^{-16}$ & $.172710^{-16}$\\
\ $3$ & $-3.994$ & $.142510^{-15}$ & $.129810^{-15}$ & $.141510^{-15}$\\
\ $4$ & $-4.216$ & $.835510^{-15}$ & $.850310^{-15}$ & $.818710^{-15}$\\
\ $5$ & $-4.444$ & $.352210^{-14}$ & $.445910^{-14}$ & $.350110^{-14}$\\
\ $6$ & $-4.660$ & $.106210^{-13}$ & $.194710^{-13}$ & $.101010^{-13}$\\
\ $7$ & $-4.825$ & $.219610^{-13}$ &  & $.207310^{-13}$\\
\ $8$ & $-4.838$ & $.283310^{-13}$ &  & $.270710^{-13}$\\
\ $9$ & $-4.411$ & $.182610^{-13}$ &  & $.164810^{-13}$\\\hline
\end{tabular}
\end{center}
\end{table}

\begin{table}[ptb]
\caption{}
\begin{center}
\medskip%
\begin{tabular}
[c]{|c|c|c|c|c|}\hline
$k$ & $\Theta$ & Exact & $\Theta G_{1}(y,z)$ & $\Theta_{1}F_{j}^{(4)}%
(y)$\\\hline
$10$ & $-3.177$ & $.238910^{-14}$ & $.229710^{-14}$ & \\
\ $11$ & $-1.650$ & $.726610^{-16}$ & $.740010^{-16}$ & \\
\ $12$ & $-.2175$ & $.857810^{-18}$ & $.863510^{-18}$ & \\
\ $13$ & $1.115$ & $.497210^{-20}$ & $.507410^{-20}$ & \\
\ $14$ & $2.378$ & $.159910^{-22}$ & $.164610^{-22}$ & \\
\ $15$ & $3.591$ & $.304410^{-25}$ & $.314510^{-25}$ & \\
\ $16$ & $4.769$ & $.354010^{-28}$ & $.374010^{-28}$ & $.146710^{-27}$\\
\ $17$ & $5.921$ & $.251910^{-31}$ & $.265010^{-31}$ & $.532610^{-31}$\\
\ $18$ & $7.052$ & $.106310^{-34}$ & $.114310^{-34}$ & $.144910^{-34}$\\
\ $19$ & $8.167$ & $.242210^{-38}$ & $.242210^{-38}$ & $.263010^{-38}$\\
\ $20$ & $9.269$ & $.228410^{-42}$ &  & $.238610^{-42}$\\\hline
\end{tabular}
\end{center}
\end{table}

\section{Conditional limit laws}

In this section we will study the probability of $k$ sources being
\textbf{on}, given that the buffer content is equal to $x>0$ and the
probability of the buffer having some non-zero content $x$, given that $k$
sources are \textbf{on,} for various ranges of $x$ and $k.$

\subsection{The distribution of active sources conditioned on $x$}

We define
\begin{equation}
M_{1}(x)\equiv\sum\limits_{k=0}^{N}f_{k}(x)=\underset{t\rightarrow\infty}%
{\lim}\frac{d}{dx}\Pr\left[  X(t)\leq x\right]  ,\quad x>0, \label{M1}%
\end{equation}%
\begin{equation}
M_{2}(k)\equiv\int\limits_{0}^{\infty}f_{k}(x)dx=\underset{t\rightarrow\infty
}{\lim}\Pr\left[  Z(t)=k,\ X(t)>0\right]  ,\quad0\leq k\leq N \label{M2}%
\end{equation}
and the conditional probabilities
\begin{equation}
f(k\ |\ x)=\frac{f_{k}(x)}{M_{1}(x)}=\underset{t\rightarrow\infty}{\lim}%
\Pr\left[  Z(t)=k\ |\ \ X(t)=x>0\right]  ,\quad x>0,\quad0\leq k\leq N
\label{fkx}%
\end{equation}%
\begin{equation}
f(x\ |\ k)=\frac{f_{k}(x)}{M_{2}(k)}=\underset{t\rightarrow\infty}{\lim}%
\frac{d}{dx}\Pr\left[  X(t)\leq x\ |\ Z(t)=k,\ X(t)>0\right]  ,\quad
x>0,\quad0\leq k\leq N. \label{fxk}%
\end{equation}
Note that
\begin{equation}
M_{2}(k)=F_{k}(\infty)-F_{k}(0),\quad0\leq k\leq N. \label{laplaceM2}%
\end{equation}

We first consider $x=O(N).$ In this range we compute $M_{1}(x)$ using the
approximation in $\mathfrak{R}_{3}$%
\begin{align*}
\sum\limits_{k=0}^{N}f_{k}(x) &  \sim\frac{1}{2\pi}\left(  \frac{\lambda
}{\lambda+1}\right)  ^{N}\int\limits_{0}^{1}\frac{1}{W_{+}(\Theta^{+},z)}%
\sqrt{\frac{-\theta_{0}}{\Theta^{+}-\theta_{0}}}\sqrt{\frac{1}{\eta_{ww}%
(W_{+},\Theta^{+},z)}}\\
&  \times\sqrt{\frac{1}{\Psi_{\vartheta\vartheta}(y,W_{+},\Theta^{+},z)}}%
\exp\left[  N\Psi(y,W_{+},\Theta^{+},z)\right]  dz.
\end{align*}
The main contribution to the integral comes from near $z=\gamma$ and by
Laplace's method we have%
\begin{align*}
\sum\limits_{k=0}^{N}f_{k}(x) &  \sim\frac{1}{\sqrt{2\pi N}}\left(
\frac{\lambda}{\lambda+1}\right)  ^{N}\sqrt{\frac{-\theta_{0}}{\Theta
^{+}\left(  y,\gamma\right)  -\theta_{0}}}\\
&  \times\sqrt{\frac{1}{\Psi_{\vartheta\vartheta}\left[  y,1,\Theta^{+}\left(
y,\gamma\right)  ,\gamma\right]  }}\exp\left\{  N\Psi\left[  y,1,\Theta
^{+}\left(  y,\gamma\right)  ,\gamma\right]  \right\}
\end{align*}
since%
\[
\frac{\partial\Psi}{\partial z}(y,W_{+},\Theta^{+},z)=0\text{ \ for
\ }z=\gamma
\]%
\[
\frac{\partial^{2}\Psi}{\partial z^{2}}(y,W_{+},\Theta^{+},z)=\frac{\Theta
^{+}}{\rho}\text{ \ for \ }z=\gamma
\]
and%
\[
W_{+}(\vartheta,\gamma)=1,\quad\eta_{ww}(1,\vartheta,\gamma)=-\frac{\rho
}{\vartheta}.
\]
Therefore,%
\begin{align*}
f(k\ |\ x) &  \sim\frac{1}{\sqrt{2\pi N}}\frac{1}{W_{+}(\Theta^{+},z)}%
\sqrt{\frac{\Theta^{+}\left(  y,\gamma\right)  -\theta_{0}}{\Theta^{+}%
-\theta_{0}}}\sqrt{\frac{1}{\eta_{ww}(W_{+},\Theta^{+},z)}}\\
\times &  \sqrt{\frac{\Psi_{\vartheta\vartheta}\left[  y,1,\Theta^{+}\left(
y,\gamma\right)  ,\gamma\right]  }{\Psi_{\vartheta\vartheta}(y,W_{+}%
,\Theta^{+},z)}}\exp\left\{  N\Psi(y,W_{+},\Theta^{+},z)-N\Psi\left[
y,1,\Theta^{+}\left(  y,\gamma\right)  ,\gamma\right]  \right\}
\end{align*}
and expanding this for $z\rightarrow\gamma$ yields the Gaussian form%
\begin{align*}
f(k\ |\ x) &  =f\left(  N\gamma+\sqrt{N}\Lambda\ |\ Ny\right)  \sim\frac
{1}{\sqrt{2\pi N}}\sqrt{\frac{-\Theta^{+}\left(  y,\gamma\right)  }{\rho}}%
\exp\left\{  \frac{1}{2}\frac{\partial^{2}\Psi}{\partial z^{2}}\left[
y,1,\Theta^{+}\left(  y,\gamma\right)  ,\gamma\right]  \left(  z-\gamma
\right)  ^{2}N\right\}  \\
&  =\frac{1}{\sqrt{2\pi N}}\sqrt{\frac{-\Theta^{+}\left(  y,\gamma\right)
}{\rho}}\exp\left\{  \frac{1}{2}\frac{\Theta^{+}\left(  y,\gamma\right)
}{\rho}\Lambda^{2}\right\}
\end{align*}
with $\Lambda=\sqrt{N}\left(  z-\gamma\right)  .$ This means that given a
buffer size $x=Ny,\ y>0,$ the most likely number of active sources will be
$c=N\gamma,$ with a Gaussian spread whose variance is obtained by solving
(\ref{Thetaplus}) with $z=\gamma.$ This is also illustrated in Figure
\ref{figcompare4}.

Next we assume a small buffer size, with $x=O(1/N).$ Again the main
contribution to the sum (\ref{M1}) comes from $k\approx c,$ but for $x=O(1/N)$
we must use the Region I approximation $f_{k}(x)\sim NF_{l}^{(1)\prime}%
(\chi).$ Thus,%
\begin{align*}
f_{k}(x)  &  \sim\frac{\sqrt{N}}{\sqrt{2\pi}}\sqrt{\frac{\rho}{\gamma
\phi\left(  1-\gamma\right)  }}\left(  \frac{\beta}{2\gamma}\right)
^{l-\alpha}\exp\left[  N\Phi(\gamma)\right] \\
&  \times\frac{1}{2\pi i}%
{\displaystyle\int\limits_{\mathrm{Br}^{+}}}
e^{\chi S}\Gamma\left(  \frac{\phi}{S}+1-\alpha\right)  J_{l-\alpha+\frac
{\phi}{S}}\left(  \frac{\beta}{S}\right) \\
&  \times\exp\left\{  \frac{\phi}{S}\ln\left[  \frac{\beta S}{2\lambda
\phi\left(  1-\gamma\right)  }\right]  +\frac{\phi-\rho}{S}\right\}  \left(
\frac{\phi}{S}\right)  ^{\alpha}dS
\end{align*}
which we can write as a residue series%
\begin{align}
f_{k}(x)  &  \sim\frac{\sqrt{N}}{\sqrt{2\pi}}\sqrt{\frac{\rho\phi}%
{\gamma\left(  1-\gamma\right)  }}\left(  \frac{\beta}{2\gamma}\right)
^{l-\alpha}\exp\left[  N\Phi(\gamma)\right] \nonumber\\
&  \times\sum\limits_{j=0}^{\infty}\exp\left[  -\frac{\phi}{j+1-\alpha}%
\chi+\frac{\left(  j+1-\alpha\right)  \left(  \rho-\phi\right)  }{\phi}\right]
\label{spectralcorner}\\
&  \times\frac{\left(  j+1-\alpha\right)  ^{j-1}}{j!}\left[  \frac
{2\lambda\left(  1-\gamma\right)  }{\beta}\right]  ^{j+1-\alpha}%
J_{l-j-1}\left[  \left(  \alpha-1-j\right)  \frac{\beta}{\phi}\right]
\nonumber
\end{align}
where we have used
\[
\operatorname{Re}\mathrm{s}\left[  \Gamma\left(  \frac{\phi}{S}+1-\alpha
\right)  ,\ S=\frac{\phi}{\alpha-1-j}\right]  =\frac{\left(  -1\right)
^{j+1}\phi}{j!\left(  j+1-\alpha\right)  ^{2}},\quad j\geq0.
\]

Using the generating function \cite{MR96i:33010}%
\[%
{\displaystyle\sum\limits_{j=-\infty}^{\infty}}
J_{j}(x)z^{j}=\exp\left[  \frac{x}{2}\left(  z-\frac{1}{z}\right)  \right]
\]
in (\ref{spectralcorner}) we obtain%
\[
\sum\limits_{l=-\infty}^{\infty}\left(  \frac{\beta}{2\gamma}\right)
^{l-\left(  j+1\right)  }J_{l-j-1}\left[  \left(  \alpha-1-j\right)
\frac{\beta}{\phi}\right]  =\exp\left[  \frac{\rho}{\phi}(j+1-\alpha)\right]
\]
and thus
\begin{align*}
M_{1}(x)  &  \sim\frac{\sqrt{N}}{\sqrt{2\pi}}\sqrt{\frac{\rho\phi}%
{\gamma\left(  1-\gamma\right)  }}\exp\left[  N\Phi(\gamma)\right]
\sum\limits_{j=0}^{\infty}\exp\left[  -\frac{\phi}{j+1-\alpha}\chi\right] \\
&  \times\frac{\left(  j+1-\alpha\right)  ^{j-1}}{j!}\exp\left[  \frac
{2\rho-\phi}{\phi}(j+1-\alpha)\right]  \left[  \frac{\lambda\left(
1-\gamma\right)  }{\gamma}\right]  ^{j+1-\alpha}.
\end{align*}
It follows that
\[
f(k\ |\ x)=f\left(  c+l-\alpha\ |\ \frac{\chi}{N}\right)  \sim
\]%
\[
\frac{\left(  \frac{\beta}{2\gamma}\right)  ^{l-\alpha}\sum\limits_{j=0}%
^{\infty}\exp\left[  -\frac{\phi}{j+1-\alpha}\chi+\frac{\left(  j+1-\alpha
\right)  \left(  \rho-\phi\right)  }{\phi}\right]  \frac{\left(
j+1-\alpha\right)  ^{j-1}}{j!}\left[  \frac{2\lambda\left(  1-\gamma\right)
}{\beta}\right]  ^{j+1-\alpha}J_{l-j-1}\left[  \left(  \alpha-1-j\right)
\frac{\beta}{\phi}\right]  }{\sum\limits_{j=0}^{\infty}\exp\left[  -\frac
{\phi}{j+1-\alpha}\chi+\frac{2\rho-\phi}{\phi}(j+1-\alpha)\right]
\frac{\left(  j+1-\alpha\right)  ^{j-1}}{j!}\left[  \frac{\lambda\left(
1-\gamma\right)  }{\gamma}\right]  ^{j+1-\alpha}}%
\]
which is a discrete distribution in $l.$

\subsection{The buffer density conditioned on $k$}

We first consider $\gamma<z\leq1.$ In this range $F_{k}(0)=0$ and from
(\ref{laplaceM2}) we have $M_{2}(k)=F_{k}(\infty).$ We have%
\[
\frac{\partial\Psi}{\partial y}\left[  y,W_{-},\Theta,z\right]  =\Theta(y,z)
\]
so that the maximum of the density occurs at $\Theta=0,$ which corresponds to
the transition curve
\[
y=Y_{0}(z)=\frac{z-\gamma}{\lambda+1}-\frac{\rho}{\left(  \lambda+1\right)
^{2}}\ln\left(  \frac{z\lambda+z-\lambda}{\rho}\right)  ,\quad\gamma<z<1.
\]
Therefore, we shall use the approximate solution $F^{(2)}(V,z),$ since Region
II corresponds to $\Theta\approx0.$ From (\ref{II}) we get, for $\gamma<z<1$%
\[
f_{k}(x)=\frac{1}{N}\frac{\partial F^{(2)}}{\partial V}\frac{\partial
V}{\partial y}\sim\frac{1}{2\pi N}\sqrt{\frac{1}{z\left(  1-z\right)  }}%
\exp\left[  \Phi(z)N\right]  \frac{1}{\sqrt{Y_{2}(z)}}\exp\left(  -\frac
{V^{2}}{2}\right)  ,
\]%
\begin{gather*}
Y_{2}(z)=\frac{2\zeta}{\left(  \lambda+1\right)  ^{4}}\ln\left(
\frac{z+z\lambda-\lambda}{\rho}\right) \\
-\frac{z-\gamma}{\left(  \lambda+1\right)  \left(  \lambda z+z-\lambda\right)
^{2}}\left[  \frac{2\zeta\rho}{\left(  \lambda+1\right)  ^{2}}+\frac{3\zeta
}{\left(  \lambda+1\right)  }\left(  z-\gamma\right)  +\left(  \lambda
-1\right)  \left(  z-\gamma\right)  ^{2}\right]  .
\end{gather*}
From (\ref{laplaceM2}) we obtain%
\[
M_{2}(k)=F_{k}(\infty)\sim\frac{1}{\sqrt{2\pi N}}\sqrt{\frac{1}{z\left(
1-z\right)  }}\exp\left[  \Phi(z)N\right]
\]
and we conclude that%
\[
f(x\ |\ k)=f\left[  NY_{0}(z)+\sqrt{N}\Omega\ |\ Nz\right]  \sim\frac{1}%
{\sqrt{2\pi}}\frac{1}{\sqrt{Y_{2}(z)N}}\exp\left\{  -\frac{1}{2}\frac
{\Omega^{2}}{Y_{2}(z)}\right\}
\]
with $\Omega=\left[  y-Y_{0}(z)\right]  \sqrt{N}.$ This shows that if the
number of active sources exceeds the service rate, the buffer content is
likely to be large, approximately $NY_{0}(z)$ and with a Gaussian spread
proportional to $Y_{2}(z)$ (see also Figure \ref{figcompare1}). While this
calculation assumed that $\gamma<z<1,$ we can show that the result remains
valid for $z\rightarrow1.$

Next we consider the range $0\leq z<\gamma.$ We use the asymptotic
approximation $F_{k}(x)-F_{k}(\infty)\sim G_{2}(y,z),$ or
\begin{align*}
f_{k}(x)  &  \sim\frac{1}{2\pi N}\left(  \frac{\lambda}{\lambda+1}\right)
^{N}\frac{1}{W_{+}(\Theta^{+},z)}\sqrt{\frac{-\theta_{0}}{\Theta^{+}%
-\theta_{0}}}\sqrt{\frac{1}{\eta_{ww}(W_{+},\Theta^{+},z)}}\\
\times &  \sqrt{\frac{1}{\Psi_{\vartheta\vartheta}(y,W_{+},\Theta^{+},z)}}%
\exp\left[  N\Psi(y,W_{+},\Theta^{+},z)\right]  .
\end{align*}
From (\ref{laplaceM2}) we have%
\[
M_{2}(k)=F_{k}(\infty)-F_{k}(0)\sim-G_{2}(0,z).
\]
Hence,%
\begin{align*}
f(x\ |\ k)  &  \sim-\Theta^{+}(0,z)\frac{W_{+}\left[  \Theta^{+}%
(0,z),z\right]  }{W_{+}(\Theta^{+},z)}\sqrt{\frac{\Theta^{+}(0,z)-\theta_{0}%
}{\Theta^{+}-\theta_{0}}}\sqrt{\frac{\Psi_{\vartheta\vartheta}\left[
0,W_{+},\Theta^{+}(0,z),z\right]  }{\Psi_{\vartheta\vartheta}(y,W_{+}%
,\Theta^{+},z)}}\\
&  \times\sqrt{\frac{\eta_{ww}\left[  W_{+},\Theta^{+}(0,z),z\right]  }%
{\eta_{ww}(W_{+},\Theta^{+},z)}}\exp\left\{  N\Psi\left[  y,W_{+},\Theta
^{+},z\right]  -N\Psi\left[  0,W_{+},\Theta^{+}(0,z),z\right]  \right\}
\end{align*}
and as $y\rightarrow0$%
\begin{align*}
f(x\ |\ k)  &  =f(x\ |\ Nz)\sim-\Theta^{+}(0,z)\exp\left[  \left\{
N\frac{\partial\Psi}{\partial y}\left[  0,W_{+},\Theta^{+}(0,z),z\right]
y\right\}  \right] \\
&  =-\Theta^{+}(0,z)\exp\left[  \Theta^{+}(0,z)x\right]  .
\end{align*}
This shows that if the number of \textbf{on} sources is less than the service
rate $c$ and the buffer is not empty, the conditional density is approximately
exponential in $x,$ with mean $\frac{1}{-\Theta^{+}(0,z)}.$ This exponential
shape is also illustrated in Figure \ref{figcompare3}.

Finally, we consider the case $k\approx c.$ We shall again use the
approximation in Region I and from (\ref{spectralcorner}) we get%
\begin{align}
M_{2}(k)  &  \sim\frac{1}{\sqrt{2\pi N}}\sqrt{\frac{\rho}{\gamma\phi\left(
1-\gamma\right)  }}\left(  \frac{\beta}{2\gamma}\right)  ^{l-\alpha}%
\exp\left[  N\Phi(\gamma)\right]  \sum\limits_{j=0}^{\infty}\exp\left[
\frac{\left(  j+1-\alpha\right)  \left(  \rho-\phi\right)  }{\phi}\right]
\nonumber\\
&  \times\frac{\left(  j+1-\alpha\right)  ^{j}}{j!}\left[  \frac
{2\lambda\left(  1-\gamma\right)  }{\beta}\right]  ^{j+1-\alpha}%
J_{l-j-1}\left[  \left(  \alpha-1-j\right)  \frac{\beta}{\phi}\right]
.\nonumber
\end{align}
Therefore,%
\[
f(x\ |\ k)=f\left(  \frac{\chi}{N}\ |\ c+l-\alpha\right)  \sim
\]%
\[
\frac{\phi N\sum\limits_{j=0}^{\infty}\exp\left[  -\frac{\phi}{j+1-\alpha}%
\chi+\frac{\left(  j+1-\alpha\right)  \left(  \rho-\phi\right)  }{\phi
}\right]  \frac{\left(  j+1-\alpha\right)  ^{j-1}}{j!}\left[  \frac
{2\lambda\left(  1-\gamma\right)  }{\beta}\right]  ^{j+1-\alpha}%
J_{l-j-1}\left[  \left(  \alpha-1-j\right)  \frac{\beta}{\phi}\right]  }%
{\sum\limits_{j=0}^{\infty}\exp\left[  \frac{\left(  j+1-\alpha\right)
\left(  \rho-\phi\right)  }{\phi}\right]  \frac{\left(  j+1-\alpha\right)
^{j}}{j!}\left[  \frac{2\lambda\left(  1-\gamma\right)  }{\beta}\right]
^{j+1-\alpha}J_{l-j-1}\left[  \left(  \alpha-1-j\right)  \frac{\beta}{\phi
}\right]  }%
\]
which is a density in $\chi,$ consisting of an infinite mixture of exponentials.

\section*{Acknowledgement}
The work of C. Knessl was partially supported by NSF Grants DMS 99-71656,
DMS 02-02815 and NSA Grant MDA 904-03-1-0036. The work of D. Dominici was
supported in part by NSF Grant 99-73231, provided by Professor Floyd Hanson.
We wish to thank him for his generous sponsorship.

\bibliographystyle{abbrv}
\bibliography{biblio}

\end{document}